\documentclass{amsart}
\usepackage{amsmath,amsthm,amscd,amssymb}
\usepackage{latexsym}
\usepackage[final]{graphicx}
\usepackage{epic, eepic}
\def\squarebox#1{\hbox to #1{\hfill\vbox to #1{\vfill}}}

\newcommand{\cz}{{\mathbb C}}
\newcommand{\nz}{{\mathbb N}}

\newcommand{\rz}{{\mathbb R}}
\newcommand{\zz}{{\mathbb Z}}

\def\tr{\mathop{\mathrm{tr}} \nolimits} 
\def\supp{\mathop{\mathrm{supp}} \nolimits} 
\def\qed{\hbox {\hskip 1pt \vrule width 4pt height 6pt depth 1.5pt
        \hskip 1pt}}

\def\canform{canonical transformation{}}

\def\fourior{Fourier integral operator{}}

\def\lhs{left hand side}

\def\neigh{neighborhood{}}

\def\pseudor{pseudodifferential operator{}}

\def\rhs{right hand side}

\def\wrt{with respect to}

\def\Re{{\rm Re\,}}
\def\Im{{\rm Im\,}}
\newtheorem{theorem}{Theorem}
\newtheorem{definition}{Definition}
\newtheorem{lemma}{Lemma}

\newtheorem{proposition}{Proposition}
\newtheorem{note}{Note}

\title{ Scattering poles near the real axis  for two strictly convex
obstacles. }
\author{Alexei Iantchenko}
\address{Malm\"o University\\
School of Technology and Society\\
SE-205 06 Malm\"o\\
Sweden
} \email{ ai@ts.mah.se}
\date{\today}
\begin{document}

   \begin{abstract}
To study the location of poles for the acoustic scattering matrix
for two strictly convex obstacles  with smooth boundaries, one uses
an approximation of the quantized billiard operator $M$ along the
trapped ray between the two obstacles. Using this method G{\'e}rard
({\slshape cf.} \cite{Gerard1988}) obtained complete asymptotic
expansions for the poles in a strip ${\rm Im}\, z\leq c$ as ${\rm
Re}\, z$ tends to infinity. He established the existence of parallel
rows of poles close to $\frac{\pi k}{d} +ij\delta ,$ $k\in\zz,$
$j\in\zz_+.$ Assuming that the boundaries are analytic and the
eigenvalues of Poincar{\'e} map are non-resonant we use  the
Birkhoff normal form for $M$ to improve his result and  to get the
 complete asymptotic
expansions for the poles in any logarithmic  neighborhood of real
axis.
\end{abstract}
\maketitle
\bibliographystyle{plain}
\pagestyle{headings}
\section{Some known results on the localization of resonances.}
\subsection{Introduction}
 Following \cite{SjostrandZworski2000} we denote by
$-\Delta_{\rz^{n+1}\setminus{\mathcal O}}$ the Dirichlet Laplacian
on a connected exterior domain $\rz^{n+1}\setminus{\mathcal O},$
where ${\mathcal O}$ is compact with a $C^\infty$ boundary. Then the
resolvent
$$R_{\mathcal O}(\lambda)\stackrel{\rm
def}{=}\left(-\Delta_{\rz^{n+1}\setminus{\mathcal O}}
-\lambda^2\right)^{-1}:\,\,L^2(\rz^{n+1}\setminus{\mathcal
O})\mapsto H^2(\rz^{n+1}\setminus{\mathcal O})\cap
H_0^1(\rz^{n+1}\setminus{\mathcal O}),\,\,{\rm Im}\,\lambda <0,$$
continues meromorphically across continuous spectrum ${\rm
Im}\,\lambda =0,$ to an operator
$$R_{\mathcal O}(\lambda):\,\,L_{\rm comp}
^2(\rz^{n+1}\setminus{\mathcal O})\mapsto H_{\rm loc}
^2(\rz^{n+1}\setminus{\mathcal O})\cap H_{0,{\rm
loc}}^1(\rz^{n+1}\setminus{\mathcal O}).$$ Here
$H^2(\rz^{n+1}\setminus{\mathcal O})$ is the standard Sobolev space,
$H_0^1(\rz^{n+1}\setminus{\mathcal O})$ is the closure of $C_{\rm
comp}^\infty(\rz^{n+1}\setminus{\mathcal O})$ in the $H^1-$norm and
by $L^2_{\rm comp},$ where by $L_{\rm comp}^2$ we mean the elements
from $L^2$ that are zero outside some bounded set, and by $H_ {\rm
loc}^2,$ $H_{0,{\rm loc}}^1$ functions that are locally in these
spaces.

 We recall
that $R_{\mathcal O}(\lambda)$ is globally meromorphic in
$\lambda\in\cz$ when $n+1$ is odd and in $\lambda\in\Lambda,$ the
logarithmic covering of the complex plane,
 when $n+1$ is even.
\begin{definition} The poles of $R_{\mathcal O}$ are called
resonances or scattering poles.\end{definition}
 Note that we adopted a
convention that the resonances $\lambda_j$ are in the upper half
plane ${\rm Im}\,\lambda >0.$

We cite some results following the review by Zworski
\cite{Zworski1994} (note his convention  that $\Im\lambda_j <0.$)

An obstacle ${\mathcal O}$ contained in a ball $B(0,R)$ is called
{\em non-trapping} if there exists $T>0$ such that every broken
characteristic ray (see e.g. \cite{HormanderLPDO}, Section 24.2)
starting in $(\rz^{n+1}\setminus {\mathcal O})\cap B(0,R)$ leaves
that set within time $T.$

Propagation of singularities for boundary value problems,
established by Melrose-Sj\"ostrand and Ivrii (see Chapter 24 of
\cite{HormanderLPDO} and the references given), gives the basic
results:
\begin{align}
  &{\mathcal O}\,\,\mbox{is non-trapping},\,\,\partial{\mathcal
  O}\,\,\mbox{is}\,\, C^\infty\,\,\Rightarrow\nonumber\\
&\mbox{there exist $a,b >0$ such that a region} \,\,\{ z;\,\,{\rm\,
Im}\,z\leq a\log(1+|z|)
+b\}\label{ln-omgivn}\\&\mbox{contains no poles;} \nonumber\\
&\forall N,\,\,\sharp\{\lambda_j:\,\, |\Im\lambda_j|\leq N\log
(1+|\lambda_j|)\}<\infty.\label{ln-number}
\end{align}

The situation is somewhat more satisfying in the analytic case where
one can apply the propagation of Gevrey 3 singularities due to
Lebeau \cite{Lebeau1984}. In fact, Bardos-Lebeau-Rauch
\cite{BardosLebeauRauch1987} (and also independently Popov
\cite{Popov1987}) used the work of Lebeau to show that
\begin{align}
  &{\mathcal O}\,\,\mbox{is non-trapping},\,\,\partial{\mathcal
  O}\,\,\mbox{is real analytic}\,\,\Rightarrow\nonumber\\
&\mbox{there exist $a,b >0$ such that a region} \,\,\{ z;\,\,{\rm\,
Im}\,z\leq a(1+|z|)^{1/3} +b\}\label{3rot-omgivn}\\&\mbox{contains
no poles.} \nonumber
\end{align}
Then clearly,
\begin{align}
&\exists B,C>0:\,\,\sharp\{\lambda_j:\,\, |\Im\lambda_j|\leq B |
\Re\lambda_j |^{1/3}\} < C.\label{3rot-number}
\end{align}

This is known to be optimal in the following sense
\cite{BardosLebeauRauch1987}: a non-degenerate isolated simple
geodesic $\gamma$ of length $d_\gamma$ on the boundary of an
analytic strictly convex obstacle in $\rz^{n+1},$ $n+1$ odd,
generates infinitely many poles in any region
\begin{equation*}
\{ z:\,\,\Im z <\frac{B}{d_\gamma} |z|^{\frac13}\},\,\, B
>\omega_\gamma
=2^{-\frac13}\zeta_1\cos\left(\frac{\pi}{6}\right)\cdot\int_0^{d_\gamma}\rho^{2/3}
(s)ds,
\end{equation*}
where $-\zeta_1$ is the first zero of the Airy function, $\rho$ is
the curvature of $\gamma$ in $\rz^{n+1}$ and $s$ the arc length
parameter.

For strictly convex obstacles the results on pole counting are
expected to be more precise and in particular the density of poles
near the real axis has already been estimated. To motivate them let
us consider (in \cite{Zworski1994} the convention was that
$\Im\lambda_j <0$)
$$N_\theta (r):=\sharp\{\lambda_j:\,\,\lambda_j\,\,\mbox{is a pole
of}\,\, R_{\mathcal O}(\lambda),\,\,|\lambda_j |\leq r,\,\,0 \leq
\arg\,\lambda_j < \theta,\,\,\pi-\theta <\arg\,\lambda_j \leq\pi
+\theta\}.$$ Then for strictly convex smooth obstacles  we get
\begin{equation}\label{3.19}
N_\theta (r) ={\mathcal O} (\theta^{3/2}) r^{n+1},\,\, r
>1.\end{equation}

By applying complex scaling ''all the way to the boundary'' which is
possible because of strict convexity, but still using a form of
functional calculus from \cite{SjostrandZworski1991}, (\ref{3.19})
was obtained in \cite{SjostrandZworski1993} for an arbitrary convex
obstacle provided $r >r(\theta).$ It was noted by Harg{\'e} and
Lebeau \cite{HargeLebeau1994} that a particular choice of the angle
of scaling yields, among other things, a cubic pole free region of
the form (\ref{3rot-number}) for strictly convex obstacles with
smooth boundaries. It had been widely believed in applied
mathematics and known already in dimensions two and three (see
\cite{BabichGrigoreva1974}).

\subsection{Trapping obstacles: results of Ikawa and Gerard }

We consider scattering by two strictly convex obstacles. Note that
scatterer is always trapping if it is not connected.

 We suppose that
\begin{align}
 &{\mathcal O}=\Omega_1\cup\Omega_2\subset\rz^{n+1},\,
 \,\mbox{with $C^\infty$
boundary}\,\, \partial\Omega_i,\,\, i=1,2,\nonumber\\
&\mbox{where}\,\, \Omega_1,\,\,\Omega_2\,\,\mbox{are compact and
{\em strictly convex}},\,\,
\Omega_1\cap\Omega_2=\emptyset.\label{hypotheseGerard}
\end{align}

  Let $d$ be the
distance between $\Omega_1$ and $\Omega_2$ and
$a_i\in\partial\Omega_i$ the points on the boundary such that
$|a_1-a_2|=d.$ Then there is one trapped broken characteristic ray
connecting $a_1$ and $a_2.$

 Under these assumptions on $\Omega$ Bardos, Guillot and
Ralston
 ({\slshape cf.} \cite{BardosGuillotRalston1982}, for $n+1$ odd) show the existence
of an infinite number of resonances in $\{ z;\,\,{\rm
Im}\,z\leq\epsilon\log|z|\}$ for any $\epsilon >0.$

 Thus their
result shows a difference in location of resonances between cases of
trapping obstacles and of non-trapping obstacles.

 The most complete results on location of poles were given by Ikawa and G\'erard. There results shows that
 uniformly in the strip $0 <{\rm Im}\, z<c,$ as ${\rm
Re}\,z$ goes to infinity the resonances are well approximated by
pseudo-poles of \cite{BardosGuillotRalston1982}:
\begin{equation}\label{pseudopoles}
\lambda_{\alpha, k}=k\frac{\pi}{d}
+\frac{i}{4d}\sum_{j=1}^{n}\ln\nu_j\cdot(2\alpha_j+1),\,\,(k,\alpha)\in\zz\times\nz^{n},\end{equation}
where $\nu_j$ are the eigenvalues of the linear Poincar{\'e} map
${\rm D}\kappa(a_1,0)$.

Though Ikawa and G\'erard only consider the case of odd dimensions
$n+1,$ their results are valid in even dimensions $n+1\geq 2$ as
well.

  In \cite{Ikawa} Ikawa obtained the first
 string of resonances closest to the real axis ($\alpha =0$) and in ({\slshape
cf.} \cite{Gerard1988}) G\'erard got the complete asymptotics of all
strings.

We start by defining the canonical transformation of the billiard
$\kappa:\,\,T^*\partial\Omega_1\mapsto T^*\partial\Omega_1$
following G\'erard \cite{Gerard1988}. For $\rho\in
T^*\partial\Omega_1$ in a \neigh{} of $(a_1,0)\in
T^*\partial\Omega_1$ we draw an outgoing half-ray $\gamma$ issued
from $\rho$ and we denote $\gamma'$ the half-ray reflected from
$\partial\Omega_2.$ If $\gamma'$ intersects $\partial\Omega_1$ then
we define $\kappa(\rho)$ the projection on $T^*\partial\Omega_1$ of
the point of intersection between $\gamma'$ and $\partial
T^*(\rz^{n+1}\setminus\Omega_1).$

For $\rho\in T^*\partial\Omega_1,$ close to $(a_1,0),$ we denote by
$\kappa_1(\rho)$ the projection on $T^*\partial\Omega_2$ of the
point of intersection of $\gamma$ and $\partial
T^*(\rz^{n+1}\setminus\Omega_2).$ $\kappa_1$ is the canonical
transformation from $T^*\partial\Omega_1$ to $T^*\partial\Omega_2.$
If we define in the similar way the \canform{}
$\kappa_2:\,\,T^*\partial\Omega_2\mapsto T^*\partial\Omega_1$ then
the billiard map is the \canform{}
\begin{equation}\label{defkappa}\kappa=\kappa_2\circ\kappa_1.\end{equation} According to Petkov
(\cite{Petkov}), \cite{PetkovStoyanov1992}, Corollary 2.3.3) and
Bardos, Guillot, Ralston (\cite{BardosGuillotRalston1982},
Proposition 3) using that $\Omega_1,$ $\Omega_2$ are strictly
convex, the eigenvalues of ${\rm D}\kappa (\rho_1,0)$ are positive
$\neq  1,$ thus $\kappa$ is of hyperbolic type and we have
$\kappa(a_1,0)=(a_1,0).$

Denote by $\nu_j$ the eigenvalues $>1$ of ${\rm D}\kappa(a_1,0)$
numerated such that $1<\nu_1\leq\nu_2\leq\ldots\leq\nu_n$ and put
$b_0=\Pi_{i=1}^n\nu_j^{-1/2}.$ For $\alpha\in\nz$ we put
$K_\alpha=b_0\nu^{-\alpha}.$ (It is possible that
$K_\alpha=K_{\alpha'}$ for $\alpha\neq\alpha'$). For any value of
$K_\alpha$ and $a_{\alpha ,j,l}\in\cz,$ we define
\begin{equation}\label{as-exp}\lambda_l(\alpha, k)= \lambda_{\alpha
,k} +\sum_{j=1}^\infty a_{\alpha ,j,l}\lambda_{\alpha
,k}^{-j/2al}\,\,\mbox{with}\,\,al\in\nz,\,\, l=1,\ldots,a
\end{equation}
 which correspond to the asymptotic expansion for the
eigenvalues of an $a\times a$ matrix, where
\begin{equation}\label{a}
a={\rm Card}\,\{\alpha'|\,\,K_{\alpha'}=K_\alpha\}.\end{equation}
Let  $p_l$ be the multiplicity of $\lambda_l(\alpha, k)$ as
asymptotic eigenvalue.
\begin{theorem}[G\'erard, $n+1$ odd] For all $A>0,$ there exists $C>0$ such that, if
$\lambda_{\alpha ,k}$ is given by (\ref{pseudopoles}), then there
exist coefficients $a_{\alpha,j,l}\in\cz$ in asymptotic expansions
(\ref{as-exp}) such that for all $N\in\nz$ there exist $k_N\in\nz$
and $c_N\in\rz$ such that the poles situated in the region
$$\{\lambda\in\cz;\,\,\Im\lambda\leq A,\,\, |\lambda|>C\}$$ are
all  in the balls
$$\left| \lambda -\left( \lambda_{\alpha ,k}
+\sum_{j=1}^{m_N}a_{\alpha ,j,l}\lambda_{\alpha
,k}^{-j/2al}\right)\right|\leq C_N |\lambda_{\alpha
,k}|^{-N},\,\,k\geq k_N\in\nz,$$ where $m_N$ is the largest $j$ such
that $\frac{j}{2al}<N,$ $a$ is defined in (\ref{a}) and in each ball
there are exactly $p_l$ poles with multiplicities.
\end{theorem}

In his article \cite{Gerard1988} G\'erard reduces the problem to the
problem on the boundary of the one obstacle, $\partial\Omega_1,$ by
introducing  the quantum billiard operator $M(\lambda)$ which
quantizes the non-linear Poincar{\'e} transform along the trapped
trajectory $\kappa$. The operator $M(\lambda)$  is defined in the
following way.

Let $H_{i,+}(\lambda):\,\, C^\infty(\partial\Omega_i)\mapsto
C^\infty(\rz^{n+1}\setminus\Omega_i),\,\,i=1,2,$ be the  outgoing
resolvent of the problem
\begin{equation}\label{lapl}
\left\{\begin{array}{l}
  (\Delta +\lambda^2) H_{i,+}(\lambda) v=0\,\,\mbox{in}\,\,\rz^{n+1}\setminus\Omega_i \\
  H_{i,+}(\lambda) v_{|\partial\Omega_i}=v
\end{array}\right.\end{equation} extended as an operator
$H^{1/2}(\partial\Omega_i)\mapsto H^1_{\rm
loc}(\complement\Omega_i)$
   and analytical for $\Im\lambda <0.$
It is known that $H_{i,+}(\lambda)$ extends  analytically  as a
bounded operator $H^{1/2}(\partial\Omega_i)\mapsto H^1_{\rm loc}
(\complement\Omega_i)$ to a domain of the form (\ref{ln-omgivn}) if
the boundary $\partial{\mathcal O}$  is $C^\infty$ and
(\ref{3rot-omgivn}) if the boundary is real analytic.

 We
define $H_i(\lambda)v=H_{i,+}(\lambda) v_{|\partial\Omega_{i+1}},$
where $\partial\Omega_3=\partial\Omega_1$ and
\begin{equation}\label{Mdef}
M(\lambda)=H_2(\lambda) H_1(\lambda)=\gamma_1H_{2,+}\gamma_2
H_{1,+}:\,\,H^1(\partial\Omega_1)\mapsto
H^1(\partial\Omega_1),\end{equation} where $\gamma_i:\,\,H_{\rm
loc}^1(\complement\Omega_i)\mapsto H^{1/2}(\partial\Omega_i)$ is the
operator of restriction to $\partial\Omega_i$.

The operator $M(\lambda)$ defined on $H^1(\partial\Omega_1)$ for
$\Im\lambda <0$ continues analytically as a bounded operator
$H^{1/2}(\partial\Omega_1)\mapsto H^{1/2}(\partial\Omega_1)$ in the
domain of the form (\ref{ln-omgivn}) if  $\partial{\mathcal O}$ is
$C^\infty$ and (\ref{3rot-omgivn}) if $\partial{\mathcal O}$ is real
analytic, and there satisfies the following estimate
$$\exists D>0,\,\,\exists C>0,\,\,\|M(\lambda)\|_{{\mathcal
L}(H^{1/2}(\partial{\Omega}_1))}\leq C|\lambda|^2
e^{D\Im\lambda^+},$$
 where
$\Im\lambda^+=\max (\Im\lambda ,0).$

 We have the following relation between the
outgoing resolvent $H_+(\lambda)$ in the exterior of ${\mathcal
O}=\Omega_1\cup\Omega_2$ ($u=H_+v$ satisfies (\ref{lapl}) with
$\partial\Omega$ instead of
 $\partial\Omega_i$)   and
the billiard operator $M,$ used by G\'erard in \cite{Gerard1988}, page 91:\\
if $(v_1,v_2)\in\ C^\infty (\partial\Omega_1)\times C^\infty
(\partial\Omega_2)$ then
$$ H_+(\lambda)(v_1,v_2)=(H_{1,+} -H_{2,+}H_1)(I-M)^{-1}v_1+(H_{2,+}
-H_{1,+}(I-M)^{-1} H_2 +H_{2,+} H_1(I-M)^{-1}H_2)v_2.$$  Using this
relation G\'erard proved in \cite{Gerard1988} that
\begin{lemma}
The scattering poles  counted  with their multiplicities coincide
with
\begin{equation}\label{M-problem}
\{\lambda,\,\,0\in\sigma(I-M(\lambda))\}.\end{equation}
\end{lemma}

Strictly speaking, G\'erard only considers the case  of $\lambda$ in
a strip
\begin{equation}\label{strip}
0\leq {\rm Im}\,\lambda\leq c_1,\,\,{\rm Re}\,\lambda\geq
c_2,\end{equation}
 for
$c_2$  sufficiently large and $n+1$ odd, but his proof also works in
more general domains (\ref{ln-omgivn}), (\ref{3rot-omgivn}).

 Microlocally near $a_1,$ G\'erard reduced the problem (\ref{M-problem}) to the problem
 of finding the points $\lambda$  in a strip (\ref{strip}) for which the operator ${\displaystyle
I-e^{-2id\lambda}M_0(\lambda)}$ considered in some appropriate space
has non-trivial kernel. Here $M_0(\lambda)$
   is a semi-classical Fourier integral operator
 associated to the
real hyperbolic canonical transformation $\kappa$ with a fixed point
at $(0,0).$ G\'erard obtained  his result by approximating of
$M_0(\lambda)$ by its linearization.

In the case when $\nu_j,$ $j=1,\ldots,n$ are linearly independent
over the field of integer numbers and the boundaries are analytic,
we will use the full Birkhoff normal form of $M_0$ (as in
\cite{IantchenkoSjostrand2002}) in order to get the explicit
formulas for the resonances in any logarithmic \neigh{} of the real
axis of the form (\ref{ln-omgivn}).

\section{Our result.}
We suppose in addition to (\ref{hypotheseGerard}) that the
boundaries $\partial\Omega_i$, $i=1,2,$ are real analytic. Suppose
that the eigenvalues of the Poincar{\'e} map
$$0< \nu_n^{-1}\leq\ldots\leq\nu_1^{-1}
<1<\nu_1\leq\ldots\leq\nu_n$$
 verify
 the
non-resonance condition: \begin{equation}\label{non-res-intro}
\sum_1^n k_j\ln\nu_j=0,\,\, k_j\in\zz\,\,\Longrightarrow\,\,
k_1=\ldots=k_n=0.
\end{equation}
In order to have good separation of the strings of resonances we
impose the Diophantine condition:
\begin{equation}\label{dioph0}
\alpha\neq\beta,\,\,|\alpha |,|\beta|\leq
m\,\,\Rightarrow\,\,|\ln\nu(\alpha -\beta) |\geq \frac{1}{C(D)}
e^{-D m},\,\,D>0.\end{equation}

Then microlocally near $a_1=0$ we can  transform $M_0(\lambda)$ into
its semiclassical Birkhoff normal form to any order $r$
$$M_{00}(\lambda):=e^{-i\lambda F^r(I_1,\ldots,I_n;1/\lambda )},\,\, I_j=\left(
x_j\partial_{x_j}+1/2\right)/i\lambda.$$
 Here for any $r\in\zz_+,$
and $\imath$ in a \neigh{} of $\imath=0$ in $\rz^n,$
 \begin{equation}\label{F}
 F^r(\imath ;1/\lambda)=
F_0(\imath)+\frac{1}{\lambda}F_1(\imath)+\frac{1}{\lambda^2}F_2(\imath)+\ldots
+\frac{1}{\lambda^r} F_r(\imath) \end{equation} with $F_j$
polynomial of degree $(r-j)$  in $\imath,$ with
\begin{align*}& F_0(\imath)=G(\imath,\mu)+H(\imath),\,\, G(\imath
,\mu)=\sum_{i=1}^n\mu_i\imath_i,\,\mu_i=\ln\nu_i,\,\,H(\imath)=\sum_{j=2}^r
h_j(\imath)=O(\imath^2),\\
&\mbox{where}\,\,h_j(\imath)=\sum_{|\alpha|=j}\frac{1}{\alpha
!}\partial^\alpha H(0)\imath^\alpha\,\,\mbox{is a homogeneous
polynomial of degree $j.$}\end{align*}  Here   $F_0(\imath)$ is real
when $\imath\in\rz^n.$

We have for any $r,r'\geq j:$
$$F_j^r(\imath)-F_j^{r'}(\imath)={\mathcal
O}\left( \imath^{\min\{r,r'\}-j+1}\right).$$ In other words,
$h_j(\imath)$ do not depend on $r.$

Thus $F^r, F_j$ are holomorphic in a \neigh{} of $\imath =0$ and
 $M_{00}(\lambda)$ is analytic in any logarithmic
\neigh{} of the real axis of the form
\begin{equation}\label{domain}\Lambda_{A,B}:=\{\lambda\in\cz\,\,\mbox{or}\,\,\Lambda;\,\,\Im\lambda
<A\ln\Re\lambda,\,\,\Re\lambda >B\}.\end{equation}

Let $$ P(\lambda):=I-e^{-2id\lambda}M_{00}(\lambda).$$

 As ${\displaystyle
I_jx_j^{\alpha_j}=\frac{1}{2i\lambda}(2\alpha_j+1)x_j^{\alpha_j}},$
the monomials $x^\alpha =x_1^{\alpha_1}\cdot
x_2^{\alpha_2}\cdot\ldots\cdot x_n^{\alpha_n}$ are formally the
eigenfunctions of the operator $P(\lambda).$ We have then
$$P(\lambda)x^\alpha=\left(1-e^{-2id\lambda}K_\alpha(\lambda)\right)x^\alpha,\,\,K_\alpha(\lambda):=e^{-i\lambda
F^r\left(\frac{2\alpha +1}{2i\lambda} ;1/\lambda\right)}.$$

We approximate the problem (\ref{M-problem})  with the following
model problem:
\begin{equation}\label{modelproblem}
1-e^{-2id\lambda}K_\alpha(\lambda) =0\,\,\Leftrightarrow\,\,
2d\lambda +\lambda F^r \left(\frac{2\alpha
+1}{2i\lambda};1/\lambda\right)=2\pi
k,\,\,\alpha\in\zz_+^n,\,\,k\in\zz,\,\,\lambda\in\Lambda_{A,B}.
\end{equation}
Note that by substituting $F^r$ by its leading term near the origin,
$\sum_{i=1}^n\ln\nu_i\cdot\imath_i,$ we recover the pseudopoles
(\ref{pseudopoles}).

In general in order to solve (\ref{modelproblem}) we observe that if
$F(\imath ;1/\lambda)=F_0(\imath)+F_1(\imath)/\lambda
+F_2(\imath)/\lambda^2+\ldots$ is either an asymptotic or a finite
sum, and $F,$ $F_j$ are holomorphic in a fixed \neigh{} of $\imath
=0,$ and $F_0(0)$ is small (in our case $F_0(0)=0$), then
\begin{align*}
&\frac{\partial}{\partial \lambda}\left(2d\lambda +\lambda F\left(
\frac{2\alpha +1}{2i\lambda};1/\lambda\right)\right)=\\
&=2d + F\left( \frac{2\alpha +1}{2i\lambda};1/\lambda\right)
-\frac{1}{\lambda}\left(\frac{\partial F\left( \frac{2\alpha
+1}{2i\lambda};1/\lambda\right)}{\partial\imath} \frac{2\alpha
+1}{2i}-\frac{\partial F\left( \frac{2\alpha
+1}{2i\lambda};1/\lambda\right)}{\partial h}\right)\approx
2d\end{align*} for $|\alpha|/|\lambda|$ small.

We omit $r$ and write equation (\ref{modelproblem}) in the form
\begin{equation}\label{modelproblem-bis}\frac{\lambda}{k}\left( 1+\frac{1}{2d}
F\left(\frac{(2\alpha +1)/k}{2i\lambda /k}
;\left(\frac{\lambda}{k}\right)^{-1}
\frac{1}{k}\right)\right)=\frac{\pi}{d}\end{equation} with $\lambda
/k$ as unknown variable.

It is then clear that equation (\ref{modelproblem-bis}) has the
solution
$$\frac{\lambda_0}{k}=g(\frac{\alpha}{k};\frac{1}{k}),\,\,g(\theta;\frac{1}{k})=\frac{\pi}{d}
+g_0(\theta) +g_1(\theta)\frac{1}{k} +\ldots,\,\,g_0(0)=0.$$

We get the solution in the form
$$\frac{\lambda_0}{k}= \frac{\pi}{d} +
\frac{a_1}{k} +\frac{a_2}{k^{2}} +\ldots =\frac{\pi}{d} -\frac{
F_1(0)}{2dk}+\frac{\Re (a_2)}{k^2}
 +\ldots
+i\left(\frac{1}{4dk}\sum_{j=1}^n\mu_j (2\alpha_j +1) +\frac{\Im
(a_2)}{k^2} +\ldots\right),$$ for $k>>1,\,\,a_j=O(|\alpha
 |^{j}),\,\,\alpha={\mathcal O}(\ln k),$ where $a_j$
  for $j=1,2,\ldots$ are polynomials in $\alpha$
(see  Theorem \ref{asymptotics} for details) and
 $F_1(0)$ is real constant as the operator $M$ is microlocally unitary (see Appendix \ref{ss-unitary}). The coefficient $a_j$ is only dependent on
 the coefficients in the Birkhoff normal form
of order $j,$ $F^j(\imath ;1/\lambda)=
F_0^j(\imath)+\frac{1}{\lambda}F_1^j(\imath)+\frac{1}{\lambda^2}F_2^j(\imath)+\ldots
+\frac{1}{\lambda^r} F_r^j(\imath),$ and  independent of
 $r$ for $j\leq r.$
 If $\lambda_0^r$ and $\lambda_0^{r'}$
 are two such solutions with $r<r'$ then
 $$\frac{\lambda_0^{r'}}{k}-\frac{\lambda_0^r}{k}=\frac{a_{r+1}^{r'}-a_{r+1}^{r}}{k^{r+1}}+{\mathcal O}\left(
 \frac{|\alpha|^{r+2}}{k^{r+2}}\right)={\mathcal O}\left(
 \frac{|\alpha|^{r+1}}{k^{r+1}}\right)$$

The main result of the paper is the following theorem:
\begin{theorem}\label{th-main-bis} Suppose (\ref{non-res-intro}) and (\ref{dioph0}) are satisfied and let $\Lambda_{A,B}$ be
 the
logarithmic \neigh{} of real axis as in (\ref{domain}). For any
$N\in\nz$ and for all $A>0$ there exists $B>0,$ and
$r=r(N)\in\nz$  such that there exist functions $a_j=a_j(\alpha),$
$j=1,\ldots, r,$ polynomial in $\alpha$ of order $\leq j$ (
$a_j={\mathcal O}\left( |\alpha|^{j}\right)$ for $|\alpha|\geq 1$)
and there exists bijection $b_{N}$ between the set of
$$\lambda_{\alpha,k}^r=k\left( \frac{\pi}{d} +\frac{a_1}{k}
+\frac{a_2}{k^2}+\ldots+\frac{a_{r+1}}{k^{r+1}}\right)\in\Lambda_{A,B}$$
and the set of resonances in $\Lambda_{A,B}$, where elements in each
set are counted with their multiplicity, such that
$$ b_{N}(\lambda) -\lambda ={\mathcal O}\left(
|\lambda|^{-N}\right).$$

\end{theorem}

Solutions to (\ref{M-problem})  can be given by means of the Grushin
problem for $I-M$ in some weighted Sj\"ostrand's space
$H(\Lambda_{tG})$ associated to the escape function $G.$ The results
of propagation of analytic singularities by Lebeau \cite{Lebeau1984}
 can be applied to localize the analytic singular support of $M.$
Outside the analytic singular support  the norm of $M$ is small. In
a \neigh{} of the analytic singular support we use the properties of
the escape function. Outside some small \neigh{} $V_0$ of
$a_1\in\partial\Omega_1$ we get that the operator $I-M(\lambda)$ is
invertible. In $V_0$ we reduce $M$ by means of an analytic
\fourior{} of Bargman type $U$ to the operator
$e^{-2id\lambda}M_0(\lambda )$ in the Birkhoff normal form such that
for any $N_0\in\nz,$
$$UM(\lambda)u=e^{-2id\lambda}M_0(\lambda )Uu+{\mathcal
O}\left(|\lambda|^{-N_0}\right) |u|_{\mathcal H}\,\,\mbox{in}\,\,
H_\Phi.$$  Here $M_0\equiv e^{-i\lambda F(I;1/\lambda)}$ in the
sense of \cite{IantchenkoSjostrand2002}. We use the truncated
Birkhoff normal form to the order $r,$ $e^{-i\lambda
F^r(I;1/\lambda)},$ which is an analytic \fourior{} and the Grushin
problem in order to show that the solutions of (\ref{M-problem}) can
be approximated by the solutions of the model problem
\begin{equation}\label{M00-problem}
\{\lambda\in\cz,\,\,0\in \sigma\left( I-e^{-2id\lambda}e^{-i\lambda
F^r(I;1/\lambda)}\right)\}\end{equation} which leads to
(\ref{modelproblem}).

 Note that as $R_{\mathcal
O}(-\overline{\lambda})=\overline{R_{\mathcal O}(\lambda)}$ it is
sufficient to consider the case $\Re\lambda
>0.$

The structure of the paper is the following:

In Section \ref{s-G} we construct the escape function.

In Section \ref{reveiw} we review the method of FBI-Bargman transforms.

In Section \ref{FBI-B} we apply the FBI-Bargman transform to the quantum billiard operator $M$ and consider the norm of the transformed operator $M_1.$

In Section \ref{s-qbnf} we reduce $M_1$ to the Birkhoff normal form up to some order $r$ at $a_1.$

In Section \ref{s-deformation} we deform the space near $a_1$ such that monomials form an almost orthonormal  basis in some \neigh of $a_1.$

In Section \ref{s-model} we derive the asymptotic expansions for the solutions of the model problem (\ref{M-problem}).

In the last sections \ref{s-gr1}-\ref{gr-global} we prove Theorem \ref{th-main-bis} by the routine method of Grushin
problems.

In Appendix \ref{s-M} we review the norm estimates of the quantum billiard operator $M$ in the Sobolev space and show that $M$ is microlocally unitary with respect to the flux norm.
\\ \\
\textsc{Acknowledgements.}
  The author is grateful to J{\"o}hannes Sj{\"o}strand   for  numerous discussions and constant support during the preparation of the manuscript.

\section{Construction of  the global escape function on $X=T^*\partial\Omega_1$ for the billiard
\canform{} $\kappa.$}\label{s-G}
\subsection{Definitions of the domains}

    Let
$\kappa:\,\,T^*\partial\Omega_1\mapsto T^*\partial\Omega_1$ be the
billiard map  defined in  (\ref{defkappa}). We have $\kappa
(a_1,0)= (a_1,0).$ For $\rho=(x,\xi)\in{\rm neigh}\, (a_1,0)\subset
T^*\partial\Omega_1$ we draw an outgoing half-ray $\gamma$ issued
from $\rho.$ If $\gamma$ intersects $\partial
T^*(\rz^{n+1}\setminus\Omega_2)$ then we
  denote
$\gamma'$ the reflected half-ray. If $\gamma'$ intersects $\partial
T^*(\rz^{n+1}\setminus\Omega_1)$ then  we define $\kappa(\rho)$ the
projection on $T^*\partial\Omega_1$ of the point of intersection
between $\gamma'$ and $\partial T^*(\rz^{n+1}\setminus\Omega_1).$
Let $D(\kappa)=\{\rho\in
T^*\partial\Omega_1;\,\,\exists\kappa(\rho)\in
T^*\partial\Omega_1\}$ be the domain of definition of $\kappa.$


 We
denote $X=\partial\Omega_1.$ The points in $T^*X$  can be divided
into the following $3$ regions:
\begin{align*}
&{\mathcal H}\,\,\mbox{hyperbolic region:}\,\, \{(y,\eta)\in
T^*X;\,\,|\eta | < 1\}\\&{\mathcal G}\,\,\mbox{glancing
region:}\,\, \{ (y,\eta)\in T^*X;\,\,|\eta | = 1\}\\
&{\mathcal E}\,\,\mbox{elliptic region:}\,\, \{ (y,\eta)\in
T^*X;\,\,|\eta |
> 1\}.
\end{align*}
We denote $\rho_1=(a_1,0)\in T^*X$ and $B(\rho_1, c):=\{\rho\in
T^*X;\,\,{\rm dist}(\rho,\rho_1)< c\}.$

 We use the following convention:  $h$ denotes
$1/\Re\lambda,$  $\lambda\in\Lambda_{A,B}.$

 Let $W_0\subset T^*X$ be a small
\neigh{} of $\rho_1$ of the size $\sqrt{h\ln(1/h)}$ given by
\begin{equation}\label{W0}W_0:= B(\rho_1,c_0\sqrt{h\ln(1/h)})\subset
T^*\partial\Omega_1,
\end{equation}
for some constant $c_0 >0.$

 Let $W_1\subset T^*X$  be a larger \neigh{} of $\rho_1$
independent of $h$ such that
$$W_0\subset\subset W_1\subset B(\rho_1,c_1)\subset D(\kappa)\subset{\mathcal H}$$
  for some constant $c_1>0$ (it will be defined later in
(\ref{W1})).

\begin{proposition}\label{l-G} There exists a real function
$G\in C^\infty(T^*\partial\Omega_1)$  with the following properties:
\begin{enumerate}
\item $G={\mathcal O}\left( h\ln^2(1/h)\right)),\,\,|\nabla G|={\mathcal O}\left(
\sqrt{h\ln (1/h)}\right),\,\,|\nabla^2 G|={\mathcal O}(1);$
\item for all ${\displaystyle \rho\in W_0,\,\,G(\kappa(\rho)) -G (\rho)\geq C\, {\rm
dist}(\rho_1,\rho)^2};$
\item for all ${\displaystyle \rho\in W_1\setminus W_0,\,\,
G(\kappa(\rho)) -G (\rho)\geq C\, h\ln(1/h);}$
\item  for ${\displaystyle \rho\in T^*\partial\Omega_1\setminus W_1,\,\,
G(\kappa(\rho)) -G (\rho)\geq C\, h\ln^2(1/h).}$
\end{enumerate} Here $C$ is some positive
constant.

 \end{proposition}
We call $G$ a (global) escape function for $\kappa.$ In the
following sections we prove Proposition \ref{l-G}.
\subsection{Construction of the local escape function $G_{\rm int}$}
 As $\kappa$ is of  hyperbolic type, by the stable
manifold theorem, in some \neigh{} $W$ of $(a_1,0),$ $\kappa$ has an
incoming stable manifold $\Lambda_-$ and an outgoing stable manifold
$\Lambda_+$ which are lagrangian manifolds, intersecting
transversally at $(a_1,0)$ characterized by
$$\Lambda_\pm=\{(x,\xi)\in W|\,\,\kappa^{-n}(x,\xi)\in
W\,\,\mbox{for all}\,\,n\in\zz_\pm\}$$ and if
$(x,\xi)\in\Lambda_\pm$ then $\kappa^{-n}(x,\xi)\mapsto (a_1,0)$
exponentially fast as $n\rightarrow\pm\infty.$

\begin{figure}[h]
\begin{center}
\includegraphics[width=50mm]{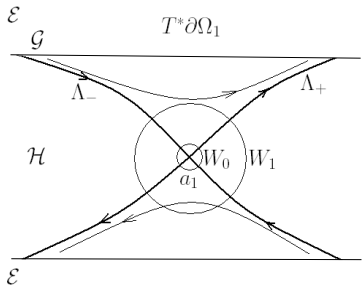}
\end{center}
\caption{{\em A part of $T^*\partial\Omega_1$ near
$(a_1,0)$.} } 
\end{figure}

There exists a local \canform{} $f:{\rm neigh}
((a_1,0),T^*X)\mapsto{\rm neigh} (0,\rz^{2n})$ such that $f
(\Lambda_+)=\{\xi =0\},$ $f(\Lambda_-)=\{ x =0\}$ in a \neigh{} $W$
of $(a_1,0)$ which corresponds to the symplectic change of
coordinates. In these new coordinates the differential of $\kappa$
at $(0,0)$ has the form
$${\rm D}\kappa(0,0)=\left(%
\begin{array}{cc}
  A & 0 \\
  0 & (A^T)^{-1} \\
\end{array}%
\right),$$ the  eigenvalues of $A$ are $\nu_1,\ldots,\nu_n.$

By the theorem of  Lewis and Sternberg ({\slshape cf.}
\cite{Sternberg1961}, \cite{Francoise1983}) there exists a
Hamiltonian $p(x,\xi)$ such that $\kappa\equiv\exp H_p$ in the sense
of formal Taylor series ({\slshape cf.} \cite{Iantchenko1998}
\cite{IantchenkoSjostrand2002}) and $p$ has the form
$p(x,\xi)=B(x,\xi)x\cdot\xi,$ where $B(0,0)=\ln A$ with eigenvalues
$\mu_j=\ln\nu_j,$ $j=1,\ldots,n.$

When $\kappa$ is of the form $\exp H_p$ near $(0,0)$ and $G\in
C^\infty$ with $dG(0,0)=0,$ then locally in a \neigh{} of $(0,0)$ we
get by Taylor expansion
$$ G(\exp H_p(\rho))=\sum\frac{1}{k!} H_p^k G (\rho)=G(\rho) +H_p
G+{\mathcal O}\left((x,\xi)^3\right)$$ and $G(\kappa(\rho))-G(\rho)=
H_pG+{\mathcal O}\left((x,\xi)^3\right).$ Thus  $G$ with properties
as in Proposition \ref{l-G} is  an escape function for $p$ in the
usual sense as defined in \cite{HelfferSjostrand1986}. This can be
used in order to construct $G$ locally near $0.$

 Suppose that the coordinates
are chosen such that
$$\kappa=\exp H_p,\,\,p=\sum_{i=1}^n\mu_ix_i\xi_i+{\mathcal O}\left((x,\xi)^3\right).$$
 It is enough to construct an escape function for
$p:=x\xi,\,\,x,\xi\in\rz.$

Inside the ball $x^2 +\xi^2 < {\mathcal O}(1) h\ln\left(1/h\right)$
we can take
$$G_0(x,\xi)=\frac12(x^2-\xi^2).$$ Then,
$H_pG_0=x^2+\xi^2={\mathcal O}\left( h\ln\left(1/h\right)\right)
,\,\,\nabla G_0={\mathcal O}\left(\sqrt{h\ln(1/h)}\right),\,\,
G_0={\mathcal O}\left( h\ln\left(1/h\right)\right).$

For $ c_1 h\ln(1/h)\leq x^2+\xi^2\leq c_2$ we use the following
ansatz: $G_1=f_h(x)-f_h(\xi).$

We have $H_p =x\partial_x -\xi\partial_\xi$ and we want that $H_p
G_1\sim h\ln(1/h)$ in this region.

We can for example choose  $f_h(x)$ satisfying
$$\left\{\begin{array}{cc}
  x\partial_x f_h (x)= h\ln(1/h), & \sqrt{h\ln(1/h)} <|x| <c \\
  f_h(\sqrt{h\ln(1/h)}) =h\ln(1/h)/2&  \\
\end{array}\right.$$ which gives
${\displaystyle
f_h(x)=h\ln(1/h)\ln\left(\frac{x}{\sqrt{h\ln(1/h)}}\right)+\frac{h\ln(1/h)}{2}}.$
And similarly for $f_h(\xi).$

With the above example in mind we define a local escape function in
the following way.

  Let
$f$ be a local real \canform{} $f(y,\eta)=(x,\xi)$
 such
that $\Lambda_-$ is transformed into $\{ x =0\}$ and $\Lambda_+$
$\mapsto$ $\{\xi =0\}$ and let $\widehat{W}_1$ be the image of
$W_1,$ $$ f(W_1)=\widehat{W}_1\subset{\rm neigh} (0,\rz^{2n}).$$ Put
$s=h\ln \left(1/h\right).$

 The local escape function $G_1$ on
$\widehat{W}_1$ is defined by $G_1=G_+ -G_-,$ where
$$G_+=\frac12 s\ln \left( s
+x^2\right),\,\,G_-=\frac12 s\ln \left( s +\xi^2\right).$$ Then
$$G_1=\frac12 s\ln\left(\frac{s + x^2}{s +\xi^2}
\right).$$

 If $x^2<<
s$ then
\begin{align*} G_+ &=\frac12s\left\{\ln(s)
+\ln\left(1+\frac{x^2}{s}\right)\right\} =\frac12 s\ln\left(
s\right) +\frac12 s\cdot\frac{x^2}{s} +s{\mathcal
O}\left(\left(\frac{x^2}{s}\right)^2\right).\end{align*} We have
then for $x^2+\xi^2<<s,$
\begin{equation}\label{as}
G_1=G_+ -G_-=\frac12 s\ln\left(\frac{s +x^2}{s +\xi^2}\right)
=\frac12(x^2 -\xi^2) + {\mathcal
O}\left(\frac{(x^2+\xi^2)^2}{s}\right).\end{equation} If $s<<x^2\leq
C,$ $\xi =0,$ we get
\begin{align*}
G_1 & =\frac12
s\ln\left(\frac{x^2(1+\frac{s}{x^2})}{s}\right)=\frac12 s\ln\left(
\frac{x^2}{s}\right)+{\mathcal O}\left( \frac{s^2}{x^2}\right)
=\frac12 h (\ln(1/h))^2 +{\mathcal O}\left( h \ln(1/h)\ln\ln
(1/h)\right).\end{align*} Thus the variation of $G_1$ is of order
$h\ln^2(1/h).$

For the derivatives we get
$$\partial_x G_1(x,\xi)=\frac{sx}{s+x^2},\,\,
\partial_\xi
G_1(x,\xi)=-\frac{s\xi}{s+\xi^2}.$$

As $$\partial_x^2 G_1(x,\xi)=
\frac{s}{s+x^2}\left(1-\frac{2x^2}{s+x^2}\right)$$ the maximum of
$\partial_x G_1$ is attained when $x^2=s$ with the maximum value
$\sqrt{s}/2.$ Using the similar estimate for $\partial_\xi G_1,$ we
get
$$  |\nabla G_1|={\mathcal O}\left( \sqrt{s}\right).$$

We have $$\partial^3_x G_1(x,\xi)=-\frac{2sx}{(s+x^2)^2}\left(
3-\frac{4x^2}{s+x^2}\right)$$ and we get that  $\partial_x^2 G_1$
attains the minimum at $x^2=3s$ with the minimum value $-1/8.$

By the similar estimate for $\partial_\xi^2 G_1,$ we get
$$  |\nabla^2 G_1|={\mathcal O}(1).$$

Similarly we get $|x\partial_x G_1|={\mathcal O}(s),$
$|\xi\partial_\xi G_1|={\mathcal O}(s).$

 We denote the billiard canonical transformation in
the new coordinates by $\hat{\kappa}:=f\circ\kappa\circ
f^{-1}:\widehat{W}_1\mapsto\widehat{W}_1.$ Denote $\rho=(x,\xi).$ We
have proved the first statement in
\begin{lemma}\label{l-G_1hat}
There exists a real function $G_1\in C^\infty(\widehat{W}_1)$ with
the following properties:
\begin{enumerate}
\item $G_1={\mathcal O}\left( h\ln^2(1/h)\right),\,\,|\nabla G_1|={\mathcal O}\left(
\sqrt{h\ln (1/h)}\right),\,\,  |\nabla^2 G_1|={\mathcal O}(1);$
\item for all $\rho\in\widehat{W}_1,$ ${\displaystyle |\rho|\leq
c\sqrt{h\ln (1/h)},\,\,G_1(\hat{\kappa}(\rho)) -G_1 (\rho)\geq C\,
|\rho |^2};$ \item for all $\rho\in\widehat{W}_1,$ ${\displaystyle
|\rho|\geq c\sqrt{h\ln (1/h)},\,\,G_1(\hat{\kappa}(\rho)) -G_1
(\rho)\geq C\, h\ln(1/h).}$
 \end{enumerate} Here $c,$ $C$ denote different
positive constants.
\end{lemma}
\underline{\bf Proof of (2) and (3):}  The canonical transformation in
the new coordinates is of the form
$\hat{\kappa}(x,\xi)=(Ax,(A^T)^{-1}\xi)+{\mathcal O}((x,\xi)^2) ,$
where $A$ is diagonal with entries $e^{\mu_i},$ $i=1,\ldots,n.$

We put $n=1,$ extension to the general case is straightforward.

Let $s:=h\ln(1/h).$ Then, choosing $\widehat{W}_1$ small enough, we
have for $\rho\in\widehat{W}_1,$
$$G_1(\rho)=\frac12 s\ln\left(\frac{s+x^2}{s+\xi^2}\right),\,\,
G_1(\hat{\kappa}(\rho))=\frac12
s\ln\left(\frac{s+e^{2\mu}x^2}{s+e^{-2\mu}\xi^2}\right)+{\mathcal
O}\left(s(x,\xi)\right).$$

 If $|\rho|^2=x^2+\xi^2 \leq c
s,$ then we can use the asymptotics (\ref{as}):
$$G_1(\hat{\kappa}(\rho)) -G_1(\rho)= (e^{2\mu} -1)x^2 -(e^{-2\mu}
-1)\xi^2+{\mathcal O}\left(\frac{|\rho|^4}{s}\right)\geq
2\mu|\rho|^2+{\mathcal O}\left(\frac{|\rho|^4}{s}\right).$$

If $cs\leq |\rho|^2=x^2+\xi^2 \leq \epsilon,$ for $\epsilon$ small
enough, we have
\begin{align} & G_1(\hat{\kappa}(\rho))-G_1(\rho)=\frac12
s\left(\ln\left(\frac{s+e^{2\mu}x^2}{s+x^2}\right)
+\ln\left(\frac{s+\xi^2}{s+e^{-2\mu}\xi^2}\right)\right)+{\mathcal
O}\left(s|\rho|\right)\geq
\label{G-1}\\
&\geq\frac12 s \left(\ln\left(\frac{1+e^{2\mu} c}{1+c}\right)
+\ln\left(\frac{1+c}{1+e^{-2\mu} c}\right)\right)+{\mathcal
O}\left(s|\rho|\right)\geq C_\epsilon s,\nonumber\end{align}using
that the \rhs{} of (\ref{G-1}) is increasing function in $x$ and in
$\xi.$ Here $C_\epsilon >0$ depends on the size of order $\epsilon$
of $\widehat{W}_1.$  \hfill\qed

In the original space (in coordinates $(y,\eta)$) in a \neigh{}
$W_1\subset T^*\partial\Omega_1$ we put
$$G_{\rm int}(y,\eta)=G_1(x,\xi),\,\, f(y,\eta)=(x,\xi)\,\,\Leftrightarrow\,\,
G_{\rm int}:=G_1\circ f^{-1}.$$

We call $G_{\rm int}$ the local escape function for $\kappa.$

\subsection{Construction of the exterior escape function $G_{\rm
ext}$}

As before $D(\kappa)$  denotes
 the domain of definition of
$\kappa$ and
$$R(\kappa):=\{\kappa(\rho),\,\,\rho\in D(\rho)\}$$ the image of
$\kappa.$
 On
$R(\kappa)$ we can define the inverse $\kappa^{-1}.$

 We define
\begin{align*}
&\widetilde{\Omega}_+(0)=\overline{\mathcal H}={\mathcal
H}\cup{\mathcal G},\,\,\widetilde{\Omega}_+(1)=D(\kappa),\,\,
\mbox{for}\,\,j\geq
1,\,\,\widetilde{\Omega}_+(j):= D(\kappa^j),\\
&\widetilde{\Omega}_+(0)\supset\widetilde{\Omega}_+(1)\supset\widetilde{\Omega}_+(2)\supset\ldots,\,\,
\bigcap_{j=0}^\infty\widetilde{\Omega}_+(j)= \Lambda_-^{\rm ext},\\
&\Lambda_-^{\rm ext}=\{(x,\xi)\in T^*X,\,\,\forall
j\in\nz,\,\,\exists\kappa^{j}(x,\xi)\in T^*X\},
\end{align*}
 where ${\rm ext}$ stands for ''extended'' by the iterated action of $\kappa.$ Note that as the obstacles
 are strictly convex  and $\rho\not\in\Lambda_-^{\rm ext}$ then
 there is $m\in\nz$ such that $\kappa^m (\rho)=\omega,$ infinite
 point.

 After $j$ reflections lie all reflected rays
except may be the last one in some \neigh{} of $(a_1,0).$ After
$\epsilon$-small perturbation of the outgoing ray it will have the
same number of reflections as an unperturbed one. Thus we have the
following lemma.
\begin{lemma}\label{lemma-distance}
For each $j=0,1,\ldots$ there exists $\epsilon (j)>0$ such that
$${\rm dist} (\widetilde{\Omega}_+(j+1),{\mathcal
H}\setminus\widetilde{\Omega}_+(j)) >\epsilon(j) >0.$$
\end{lemma}

Define
$\Omega_+(j):=\widetilde{\Omega}_+(j)\setminus\widetilde{\Omega}_+(j+1)\subset
T^*X$ such that the half-rays issued from the points of
$\Omega_+(j)$ come back exactly $j$ times and then disappears to
infinity. For example $\Omega_+(1)$ corresponds to the region in
${\mathcal H}$ such that if $\rho\in\Omega_+(1)$ than
$\kappa(\rho)=\rho '\in T^*X$ but $\kappa(\rho ')=\omega$ --
infinite point (only one reflection). We have \begin{align*}
&\Omega_+(0):={\mathcal H}\cup{\mathcal G}\setminus
D(\kappa),\,\,\Omega_+(1):=\kappa^{-1}(\Omega_+(0)),\,\,\ldots\\
&{\mathcal H}\cup{\mathcal
G}=\bigcup_{j=0}^\infty\Omega_+(j),\,\,\kappa:\,\,\Omega_+(j)\mapsto\Omega_+(j-1).
\end{align*}
Lemma \ref{lemma-distance} shows that approaching the glancing
surface (possibly tangentially) $\Omega_+(j),$ $j=0,1,\ldots,$ are
well separated.

We define also
\begin{align*} &\Omega_+(j+\frac12):=\{\rho\in\Omega_+(j);\,\,
{\rm dist}(\rho,\Omega_+(j+1) <{\rm dist}(\rho,\Omega_+(j-1))\},\\
&\Omega_+(j-\frac12):=\{\rho\in\Omega_+(j);\,\, {\rm
dist}(\rho,\Omega_+(j+1))\geq {\rm dist}(\rho,\Omega_+(j-1))\}.
\end{align*}

Thanks to Lemma \ref{lemma-distance} these domains extend naturally
to ${\mathcal E}.$ We denote such extended domains again by
$\Omega_+(j).$

Let $G_+^0$ be a step function which is  constant  on each
$\Omega_+(j),\,\,j=0,\frac12,1,\frac32,\ldots, N,$ defined as
follows \begin{align*}& G^0_{+|\cup_{j\geq
N}\Omega_+(j)}=0,\,\,G^0_{+|\Omega_+(N-1/2)}=\frac14 h\ln^2(1/h),\\
&G^0_{+|\Omega_+(N-1)}=\frac12
h\ln^2(1/h),\,\,\ldots,\,\,G^0_{+|\Omega_+(0)}=N\cdot\frac12
h\ln^2(1/h)\\
&G^0_{+|{\mathcal E}}=0
\end{align*}

Let $\omega\in C^\infty_0(\rz^n)$ be such that $\omega (x)\geq 0,$
$\omega (x)=0, $ for $|x| >1,$ $\int_{\rz^n}\omega (x) dx=1.$ Put
$\omega_\epsilon (x) =\epsilon^{-n}\omega \left(x /\epsilon\right).$
For $u\in C^\infty(\rz^n)$ we consider convolution
$$(\omega_\epsilon\ast u)(x)=\int\omega_\epsilon (y)
u(x-y)dy=\int\omega (y) u(x-\epsilon y)dy.$$ The mollifier operator
$u\mapsto\omega_\epsilon\ast u$ extends to the continuous linear
operator $D'\mapsto C^\infty$ with
$$\supp(\omega_\epsilon\ast u)\subset\supp u +\{
x\in\rz^n:\,\,|x|\leq C\epsilon\}$$ for some $C>0.$

Using the partition of identity on $T^*X$ and mollifiers on any
coordinate \neigh{} we can construct the regularization $G_+\in
C^\infty(T^*X)$ such that
for $\rho\in\bigcup_{j=1}^N\Omega_+(j)$ and $\epsilon$ small enough
$$G_+(\kappa(\rho))- G_+(\rho)\geq\frac14 h(\ln^2(1/h))> 0.$$

In the similar way we define $G_-$ associated to $\Omega_-(j)$
defined as before but for the incoming half-rays (by replacing
$\kappa$ with $\kappa^{-1}$):
\begin{align*}
&\widetilde{\Omega}_-(0)=\overline{\mathcal H}={\mathcal
H}\cup{\mathcal
G},\,\,\widetilde{\Omega}_-(1)=D(\kappa^{-1})=R(\kappa),\,\,
\mbox{for}\,\,j\geq
1,\,\,\widetilde{\Omega}_-(j):= D(\kappa^{-j})=R(\kappa^j),\\
&\widetilde{\Omega}_-(0)\supset\widetilde{\Omega}_-(1)\supset\widetilde{\Omega}_-(2)\supset\ldots,\,\,
\bigcap_{j=0}^\infty\widetilde{\Omega}_-(j)= \Lambda_+^{\rm ext},\\
&\Lambda_+^{\rm ext}=\{(x,\xi)\in T^*X,\,\,\forall
j\in\nz,\,\,\exists\kappa^{-j}(x,\xi)\in T^*X\},
\end{align*}
where $\Lambda_+^{\rm ext}$ denotes extension of $\Lambda_+$ by the
iterated action of $\kappa^{-1}.$ If $\rho\not\in\Lambda_+^{\rm
ext}$ then
 there is $m\in\nz$ such that $\kappa^{-m} (\rho)=\omega.$
We have an analogue of Lemma \ref{lemma-distance}.

 Define
$\Omega_-(j):=\widetilde{\Omega}_-(j)\setminus\widetilde{\Omega}_-(j+1).$
 We have \begin{align*}
&\Omega_-(0):={\mathcal H}\cup{\mathcal G}\setminus
D(\kappa^{-1}),\,\,\Omega_-(1):=\kappa(\Omega_-(0)),\,\,\ldots\\
&{\mathcal H}\cup{\mathcal
G}=\bigcup_{j=0}^\infty\Omega_-(j),\,\,\kappa^{-1}:\,\,\Omega_-(j)\mapsto\Omega_-(j-1).
\end{align*}
\begin{align*} &\Omega_-(j+\frac12):=\{\rho\in\Omega_-(j);\,\,
{\rm dist}(\rho,\Omega_-(j+1) <{\rm dist}(\rho,\Omega_-(j-1)\},\\
&\Omega_-(j-\frac12):=\{\rho\in\Omega_-(j);\,\, {\rm
dist}(\rho,\Omega_-(j+1))\geq {\rm dist}(\rho,\Omega_-(j-1))\}.
\end{align*}

We denote the extension of $\Omega_-(j)$ to ${\mathcal E}$ again by
$\Omega_-(j).$

Let $G_-^0$ be the step function:
 \begin{align*}& G^0_{-|\cup_{j\geq
N}\Omega_+-(j)}=0,\,\,G^0_{-|\Omega_-(N-1/2)}=\frac14 h\ln^2(1/h),\\
&G^0_{-|\Omega_-(N-1)}=\frac12
h\ln^2(1/h),\,\,\ldots,\,\,G^0_{-|\Omega_-(0)}=N\cdot\frac12
h\ln^2(1/h),\\
&G^0_{-|{\mathcal E}}=0.
\end{align*}

We define the regularization $ G_-\in C^\infty(T^*X)$ such that
\begin{align*}&\forall\,\,\rho\in\bigcup_{j=1}^N\Omega_-(j),\,\,G_-(\kappa^{-1}(\rho))
-G_-(\rho )\geq\frac14 h\ln^2(1/h)
>0\,\,\Leftrightarrow\\&\Leftrightarrow G_-(\kappa(\rho '))-G_-(\rho ')\leq
-\frac14 h\ln^2(1/h)
<0,\,\,\rho'=\kappa^{-1}(\rho)\in\bigcup_{j=0}^{N-1}\Omega_-(j)
.\end{align*}

 Let $G_{\rm ext}:=G_+-G_-$.

 We define \begin{equation}\label{W2}W_2=W_2(N):=
 \bigcup_{i=1}^N\Omega_+(i)\cup\bigcup_{j=0}^{N-1}\Omega_-(j).\end{equation}
\begin{lemma}\label{l-W2} For all $\rho\in
W_2(N)$ we have
$$ G_{\rm ext}(\kappa(\rho)) -G_{\rm ext}(\rho)\geq\frac12 h\ln^2(1/h).$$
\end{lemma}

\underline{\bf Proof:} If $\rho\in\Omega_+(i),$ $i\in [1,N],$ then
$\exists$ $j\in [0,N-1]$ such that $\rho\in\Omega_-(j),$ then we
have
\begin{align*} &G_{\rm ext}(\kappa(\rho)) -G_{\rm
ext}(\rho)=\left( G_+(\kappa(\rho)) -G_+(\rho)\right) -\left(
G_-(\kappa(\rho)) -G_-(\rho)\right)\geq\\
&\geq \frac14 h\ln^2(1/h)-(-\frac14
h\ln^2(1/h)).\end{align*}\hfill\qed

Thus we gain at least $\frac12 h\ln^2(1/h)$ on $W_2.$

We have that $\Lambda_-^{\rm
ext}\cup\bigcup_{i=N+1}^\infty\Omega_+(i)=\Omega_+(N+1)$ is a
\neigh{} of $\Lambda_-^{\rm ext},$ and $\Lambda_+^{\rm
ext}\cup\bigcup_{j=N}^\infty\Omega_-(j)=\Omega_-(N)$ is  a \neigh{}
of $\Lambda_+^{\rm ext}.$

Let \begin{equation}\label{W1} W_1=\Omega_+(N+1)\cap\Omega_-(N).
\end{equation}
If $$\rho\in\Lambda_-^{\rm
ext}\cup\bigcup_{i=N+1}^\infty\Omega_+(i)\setminus
W_1=\Omega_+(N+1)\setminus\left(\Omega_+(N+1)\cap\Omega_-(N)\right)={\rm
neigh}(\Lambda_-^{\rm ext})\setminus W_1$$ then $G_+(\rho)=0$ and
there is $j$ such that $\rho\in\bigcup_{j=0}^{N-1}\Omega_-(j)$ and
$G_{\rm ext}(\kappa(\rho))- G_{\rm ext}(\rho)=-G_-(\kappa(\rho))+ G_-(\rho)
\geq\frac14 h\ln^2(1/h).$ Similar for $\rho\in\Lambda_+^{\rm
ext}\cup\bigcup_{j=N}^\infty\Omega_-(j)\setminus W_1.$ Thus we gain
there too at least $\frac14 h\ln^2(1/h).$

Together with Lemma \ref{l-W2} we get that on $\overline{\mathcal
H}\setminus W_1$ we have $G_{\rm ext}(\kappa(\rho)) -G_{\rm
ext}(\rho)\geq\frac14 h\ln^2(1/h).$

  Inside $W_1$ we use
$G_{\rm int}:=G_1\circ f^{-1},$ the interior escape function defined
in the previous section, and apply Lemma \ref{l-G_1hat}.

 Let $C$ be a large constant and
$\chi\in C_0^\infty(T^*X)$ such that $\chi =1$ in $W_1.$
 We
define the global escape function
$$G(\rho)=\frac{1}{C}\chi (\rho) G_{\rm int}(\rho ) + G_{\rm
ext}(\rho)\in C^\infty(T^*X).$$

Then $$G(\kappa(\rho))=\frac{1}{C}\chi (\kappa(\rho)) G_{\rm
int}(\kappa (\rho)) +G_{\rm ext}(\kappa(\rho))$$ and
$G(\kappa(\rho)) -G(\rho)$ verifies the inequalities stated in Lemma
\ref{l-G}.\hfill\qed

\section{The FBI  transform of Helffer and Sj\"ostrand and the Bargman transform with the global choice of
phase.}\label{reveiw} In this section we remind the definition of the space
$H(\Lambda_{tG})$ as in \cite{HelfferSjostrand1986}. We follow the
presentation in \cite{Sjostrand1996}.
\subsection{Complex domains and symbol classes: generalities.}

Let $\tilde{X}$ denote a complex \neigh{} of $X.$ Let $\Lambda\in
T^*\widetilde{X}$ be a closed I--Lagrangian manifold which is close
to $T^*X$ in the $C^\infty$--sense and which coincides with this set
outside a compact set. Recall that ''I--Lagrangian'' means
Lagrangian for the real symplectic form $\Im\sigma,$ where
$\sigma=\sum d\alpha_{\xi_j}\wedge d\alpha_{x_j}$ is the standard
complex symplectic form. This means that if we choose (analytic)
coordinates $y$ in $X$ and let $(y,\eta)$ be the corresponding
canonical coordinates on $T^*X$ and $T^*\widetilde{X},$ then
$\Lambda$ is of the form ${\displaystyle \{(y,\eta)
+iH_{tG(y,\eta)};\,\,(y,\eta)\in T^*X\}}$ for some real-valued
smooth function $tG(y,\eta)$ which is close to $0$ in the
$C^\infty$--sense and has compact support in $\eta.$

 Here
$H_G$ denotes the Hamilton field of $G.$ Since $\Lambda$ is close to
$T^*X,$ it is also R--symplectic in the sense that the restriction
to $\Lambda$ of $\Re\sigma$ is non-degenerate. We say that $\Lambda$
is an IR--manifold.

 A smooth function $a(x,\xi;h),$
defined on $\Lambda$ or on a suitable neighborhood of $T^*X$ in
$T^*\widetilde{X}$ is said to be of class $S^{m,k},$ if
$\partial_x^p\partial_\xi^q a={\mathcal O} (1)
h^{-m}\langle\xi\rangle^{k-|q|}.$

A formal classical symbol $a\in S_{\rm cl}^{m,k}$ is of the form
$a\sim h^{-m} (a_0 +h a_1 +\ldots)$ where $a_j\in S^{0, k-j}$ is
independent of $h.$ Here and in the following, we let $0<h\leq h_0$
for some sufficiently small $h_0 >0.$ When the domain of definition
is real or equal to $\Lambda,$ we can find a realization of $a$ in
$S^{m,k}$ (denoted by the same letter $a$) so that
$$a-h^{-m}\sum_0^N h^j a_j\in S^{-(N+1) +m, k-(N+1)}. $$
When the domain of definition is a complex domain, we say that $a\in
S_{\rm cl}^{m,k}$ is a formal classical analytic symbol ($a\in
S_{\rm cla}^{m,k}$) if $a_j$ are holomorphic and satisfy
$|a_j|\leq C_0 C^j (j!) |\langle\xi\rangle |^{k-j}.$

It is then standard, that we can find a realization $a\in S^{m,k}$
(denoted by the same letter $a$) such that
\begin{align*}
&\partial_x^k\partial_\xi^l\overline{\partial}_{x,\xi} a={\mathcal
O} (1) e^{-|\langle\xi\rangle |/Ch},\,\,
\left| a-h^{-m}\sum_{0\leq j\leq |\langle\xi\rangle |/C_0 h}h^j a_j
\right|\leq{\mathcal O} (1) e^{-|\langle\xi \rangle
|/C_1h},\nonumber
\end{align*}
where in the last estimate $C_0 >0$ is sufficiently large and $C,C_1
>0$ depend on $C_0.$ We will denote by $S_{\rm cl}^{m,k}$ and
$S_{\rm cla}^{m,k}$ also the classes of realizations of classical
symbols. We say that a classical (analytic) symbol $a\sim h^{-m}
(a_0+ha_1 +\ldots)$ is elliptic, if $a_0$ is elliptic, so that
$a_0^{-1}\in S^{0,-k}.$

\subsection{The FBI transforms }

Let $X=\partial\Omega_1.$ Following  Sj\"ostrand
\cite{Sjostrand1996} and Zworski \cite{Zworski1999}  we introduce an
FBI transform which is a map $T^{HS}:C^\infty (X)\mapsto
C^\infty(T^*X)$ given by
\begin{equation}\label{HeSjFBI}
 T^{HS}u(\alpha ,h)=\int_X e^{i\phi (\alpha ,x)/h}a(\alpha,x
;h)\chi (\alpha ,x) u(x)dx,
\end{equation}
 where $\chi$ is smooth cut-off with
support close to the diagonal $$\Delta=\{ (\alpha,x)\in T^*X\times
X,\,\,\alpha_x=x\}$$ and equal to $1$ in a \neigh{}   of $\Delta,$
$a\in S^{\frac{3n}{4},\frac{n}{4}}$ is elliptic, and $\phi$ is an
admissible phase function. There exists $b(\alpha,x;h)\in
S^{\frac{3n}{4},\frac{n}{4}}$ such that if
\begin{equation}\label{inverseHeSjFBI} S^{HS}v(x ,h)=\int_{T^*X} e^{-i\phi^*
(\alpha ,x)/h}b(\alpha,x ;h)\chi (\alpha ,x)
v(x)dx,\,\,\phi^*(\alpha,x)=\overline{\phi
(\overline{x},\overline{\alpha})}, \end{equation} then
$S^{HS}T^{HS}u=u +Ru,$ where $R$ has a distributional kernel
$R(x,y;h)$ satisfying
$$|\partial_x^k\partial_y^l R| \leq C_{k,l} e^{-1/c_0h}.$$

Let $G\in C^\infty(T^*X)$ be the escape function constructed in the
previous sections. We define its $C^\infty$ extension to
$\widetilde{T^*X},$ the complex \neigh{} of $T^*X$,  with the
property that
$$(dG)_\rho|_{JT_\rho(T^*X)}=0,$$ where
$J:\,\,T_\rho(\widetilde{T^*X})\mapsto T_\rho(\widetilde{T^*X})$ is
the complex involution. Denoting the extension by the same symbol,
$G\in C^\infty(\widetilde{T^*X}),$ we now define a $C^\infty$
$I-$lagrangian, $R-$symplectic submanifold of $\widetilde{T^*X}:$
\begin{equation}\label{Lambda}
\Lambda_{tG}=\exp (tH_G^{\Im\sigma})(T^*X)\subset
\widetilde{T^*X},\,\,\Im\sigma
(x,H_G^{\Im\sigma})=dG(x),\,\,\sigma=d\alpha_\zeta\wedge d\alpha_z.
\end{equation}

Here by a $C^\infty$  manifold we mean a graph of a $C^\infty$
function. We remark that a $C^\infty$ graph in $\widetilde{T^*X}$
over $T^*X$ which is $I-$lagrangian (i.e. its almost everywhere
defined tangent plane is Lagrangian \wrt{} $\Im\sigma$) can be
locally written as a graph of a differential of a $C^\infty$
function, $\widetilde{G}$ on $T^*X\simeq T^*\rz^n$ (locally):
\begin{align*}&
\{\alpha=(\alpha_x,\alpha_\xi)\in T^*\cz^n:\,\,\Im
\alpha_x=\frac{\partial\widetilde{G}}{\partial\Re\alpha_\xi} (\Re
\alpha_x,\Re\alpha_\xi),\,\,\Im\alpha_\xi
=-\frac{\partial\widetilde{G}}{\partial\Re\alpha_x}
(\Re \alpha_x,\Re\alpha_\xi)\},\\
&\,\,(\Re \alpha_x,\Re\alpha_\xi)\in T^*\rz^n\simeq
T^*X,\,\,\widetilde{G}=tG.\end{align*} The form $-\Im\zeta
dz|_{\Lambda_{tG}}$ is (formally) closed and, as $\Lambda_{tG}$ is
close to $T^*X,$ it has a primitive which is a $C^\infty$ function
on $\Lambda_{tG}.$ We denote it by $H\in C^\infty(\Lambda_{tG}
;\rz).$ We normalize $H$ by demanding that it is equal to $0$ for
$|\alpha_\zeta |$ large enough.

Parametrizing $\Lambda_{tG}$ by $T^*X:$
$$\Lambda_{tG}=\{ (x+it\frac{\partial G}{\partial \xi},\xi
-it\frac{\partial G}{\partial x});\,\,(x,\xi)\in T^*X\},$$
 we get
\begin{equation}\label{identification} H_t(\alpha_x,\alpha_\xi)=-\xi \cdot t\partial_\xi G(x,\xi)
+tG(x,\xi),\,\,(\alpha_x,\alpha_\xi)=(x,\xi)+itH_G(x,\xi),\,\,(x,\xi)\in
T^*X.
\end{equation}

The weight function $H_t(\alpha_x,\alpha_\xi)$ have the same
properties as $tG$ due to Proposition \ref{l-G}.
 Thanks to the global properties of
the amplitude and the phase function of $T^{HS}$ we can continue it
in $\alpha$ to a \neigh{} of $T^*X,$ $\widetilde{T^*X},$ in
$T^*\tilde{X}$ and in particular we can define
$$T_{\Lambda_{tG}} u(\alpha ;h)=Tu |_{\Lambda_{tG}}(\alpha
;h),\,\,\alpha\in\Lambda_{tG}.$$ A deformation argument as in
Sj\"ostrand gives an approximate inverse $S_{\Lambda_{tG}},$ with
the same properties as $S$ above, defined by (\ref{inverseHeSjFBI})
with $T^*X$ replaced by $\Lambda_{tG}.$ We put for $u\in C^\infty$
\begin{equation}\label{normHS}\|u\|_{tG}^2=\|
u\|^2_{H(\Lambda_{tG},\langle\alpha_\xi\rangle^m)}=\|T_{\Lambda_{tG}}
u\|^2_{L^2(\Lambda_{tG},\langle\alpha_\xi\rangle^m
e^{-H_t/h})}=\int_{\Lambda_{tG}} |T_{\Lambda_{tG}} u(\alpha
;h)|^2|\langle\alpha_\xi\rangle |^{2m} e^{-2H_t(\alpha
)/h}d\alpha.\end{equation} Using $S_{\Lambda_{tG}}$ we can show the
independence of this norm of the choice of a specific phase
function.

From \cite{Sjostrand1996} and \cite{Zworski1999} it follows also
that the kernel of $T_{\Lambda_{tG}}S_{\Lambda_{tG}}$ has the form
$$T_{\Lambda_{tG}}S_{\Lambda_{tG}}(\alpha,\beta)=c(\alpha,\beta)e^{-\frac{i}{h}\psi(\alpha,\beta)}
+r(\alpha,\beta;h),$$ where
$\psi(\alpha,\beta)=v.c._z(\phi(\alpha,z)-\phi^*(\beta,z)),$
$c(\alpha,\beta;h)={\mathcal O}(h^{-n})$ is supported near
$\alpha=\beta,$ $$r(\alpha,\beta;h)={\mathcal O}_N(1)
e^{-1/Ch}|\langle\beta_\zeta\rangle|^{-N }$$ and
$$-H(\alpha)-\Im\psi(\alpha,\beta) +H(\beta)\leq
-\frac{1}{c}d(\alpha,\beta)^2.$$ Here $d(\alpha,\beta)$ is any
non-degenerate distance on $\Lambda_{tG}.$ We used that
 $|\nabla^2 G|={\mathcal O}(1).$

\subsection{Generalities on the Bargman transforms with the global
choice of phase and choice of the  norm.}\label{ss-Bargman} Following
\cite{Sjostrand1996} we replace the norm $\|\cdot\|_{tG}$ by an
equivalent norm obtained by decomposition of $T^*X.$

 In the region where
$\xi$ is bounded, it will be convenient to work with transforms
which are holomorphic up to exponentially small errors, and for that
we make a different choice of $T,$ and take an FBI--transform of
Bargman type with a global choice of phase.

 The standard  Bargman transform in $\rz^n$ is given by  $$T_0
u(z,h)=h^{-\frac{3n}{4}}\int e^{-(z-x)^2/2h} u(x) dx,\,\, T_0:\,
L^2(\rz^n)\mapsto H_{\Phi_0}(\cz^n),\,\,\Phi_0(z)=(\Im z)^2/2,$$
where
$H_\Phi:=\{u\,\,\mbox{holomorphic},\,\,\int|u|^2e^{-2\Phi/h}dx<\infty\},$
associated to the complex canonical transformation $$\kappa_{T_0}
(x,\xi)=(x-i\xi
,\xi),\,\,\kappa(\rz^{2n})=\Lambda_{\Phi_0}:=\{z=\frac{2}{i}\frac{\partial\Phi_0(z)}{\partial
z}\}.$$

Let $\phi_0(\alpha,x)=i(\Im z)^2/2 +i(z-x)^2/2,$ $z=\alpha_x
-i\alpha_\xi.$ Then $\phi_0(\alpha ,x)$ is admissible phase function
and thus $e^{-2\Phi_0 (z)/h}T_0$ is local FBI transform of the type
described in the previous section:
$$e^{-2\Phi_0(z)/h}T_0 u(\alpha)= h^{-\frac{3n}{4}}\int
e^{\frac{i}{h}((\alpha_x-x)\alpha_\xi +i(x-\alpha_x)^2/2)}u(x)dx.$$
$$\kappa_{T_0} (\Lambda_{tG})=\Lambda_{\Phi_t},\,\,\Phi_t(z)=\Phi_0
(z) +tG(\Re z,-\Im z) +{\mathcal O}(t^2).$$

We need the global version  of Bargman transform on $T^*X.$ We equip
$X$ with some analytic Riemannian metric so that we have a distance
${\rm dist}$ and a volume density $dy.$

The function ${\rm dist}(x,y)^2$ is analytic in a \neigh{} of the
diagonal in $X\times X,$ so we can consider it as a holomorphic
function in a region
$$\left\{ (x,y)\in \widetilde{X}\times\widetilde{X};\,\,{\rm
dist}\,(x,y) <\frac{1}{C},\,\,|\Im x|,|\Im y|<\frac{1}{C}\right\}.$$
Put \begin{equation}\label{A.22} \phi(x,y)=i\lambda_0 {\rm
dist}(x,y)^2,
\end{equation}
where $\lambda_0 >0$ is a constant that we choose large enough,
depending on the size of the \neigh{} of the zero section in $T^*X,$
that we wish to cover.

For $x\in\widetilde{X},$ $|\Im x |<1/C,$ put
\begin{equation}\label{A.23}
T u(x;h) =h^{\frac{-3n}{4}}\int e^{\frac{i}{h} \phi (x,y)}\chi (x,y)
u(y) dy,\,\, u\in {\mathcal D}'(X),
\end{equation}
where $\chi$ is a smooth cut-off function with support in
$$\{(x,y)\in\widetilde{X}\times X;\,\,|\Im x| <1/C,\,\, {\rm dist}(y,y(x))
<1/C\}$$ containing ${\mathcal H}\cup {\mathcal G}.$ Here $y(x)\in
X$ is the point close to $x$, where $X\ni y\mapsto -\Im\phi(x,y)$
attains its non-degenerate maximum. We collect
 the following facts from \cite{Sjostrand1996}:\\
\begin{itemize}
\item The function $\Phi_0 (x)=-\Im\phi (x,y(x)),$
$x\in\widetilde{X},$ $|\Im x|<1/C,$ is strictly plurisubharmonic and
is of the order of magnitude $\sim |\Im x |^2.$ \item
$\Lambda_{\Phi_0}:=\{ (x,\frac{2}{i}\partial\Phi_0)\in
T^*\widetilde{X}\}$ is an IR--manifold given by  $\Lambda_{\Phi_0}
=\kappa_ T(T^* X),$ where $\kappa_T$ is the complex canonical
transform associated to $T,$ given by
$(y,-\phi_y'(x,y))\mapsto(x,\phi_x'(x,y)).$ Here and in the
following, we identify $\widetilde{X}$ with its intersection with a
tubular \neigh{} of $X$ which is independent of the choice of
$\lambda_0$ in (\ref{A.22}). \item If ${\displaystyle L_{\Phi_0}^2=
L^2(\widetilde{X};\,\, e^{-2\Phi_0/h} L(dx)),}$ for $L(dx)$ denoting
a choice of Lebesgue measure (up to a non-vanishing continuous
factor), then $T={\mathcal O}(1):\,\,L^2(X)\mapsto L_{\Phi_0}^2,$
$\overline{\partial}_x T={\mathcal O} (e^{-1/Ch}):\,\,L^2(X)\mapsto
L_{\Phi_0}^2.$ This means that up to an exponentially small error $
Tu$ is holomorphic for $u\in L^2(X)$ (and even for $u\in{\mathcal
D}'(X)$).  \item Unitarity: Modulo exponentially small errors and
microlocally, $ T$ is unitary $${\displaystyle L^2(X)\mapsto
L^2(\widetilde{X}; a_0e^{-2\Phi_0/h} L(dx))},$$ where $L(dx)$ is
chosen as indicated above, and $a_0 (x;h)$ is a positive elliptic
analytic symbol of order $0.$ \item Let $\Lambda_{tG}\subset
T^*\widetilde{X}$ be an IR--manifold as before. Then
$\kappa_T(\Lambda_{tG} )=\Lambda_{\Phi_t},$ where $\Phi_t$ can be
normalized by the requirement that $\Phi_t=\Phi_0$ near the boundary
of $\widetilde{X}.$ (Here is where we have to choose $\lambda_0$
large enough, depending on $\Lambda_{tG}.$ In the applications, for
a given elliptic operator, $\Lambda_{tG}$ and $T^*X$ will coincide
outside a fixed compact \neigh{} of the zero section, and the whole
study will be carried out with a fixed $\lambda_0.$)
\end{itemize}

  Locally we can identify $T^*X\simeq T^*\rz^n$  and take $\phi
(x,y)=i(x-y)^2/2.$ Then locally $\kappa_T:\,\,(y,\eta)\mapsto
(y-i\eta,\eta)=:\Lambda_{\Phi_0}$ with $\Phi_0=(\Im x)^2/2.$ Using
the local formula
$$H_t(y,\eta)=-\Re\eta\cdot\Im y +tG(\Re y,\Re \eta)$$ we get
\begin{align}&\Phi_t(x)=v.c._{(y,\Re \eta)\in\cz^n\times\rz^n}\left(-\Im\phi(x,y) -\Re\eta\cdot\Im
y +tG(\Re y,\Re\eta)\right)\,\,\Leftrightarrow\nonumber\\
&\Phi_t(x)=\Phi_0(x)+tG(x)+{\mathcal O} (t^2|x|^2),\,\,\nonumber\\
&\mbox{where}\,\, G(x)=G(y(x),\eta (x)),\,\,(y(x),\eta
(x))=\kappa_T^{-1}\left( x,\frac{2}{i}\frac{\partial \Phi_0
(x)}{\partial x}\right).\label{G(x)}
\end{align}

Let $W$ be a \neigh{} of $\rho_1\in T^*X,$ such that
$W_1\subset\subset W.$ We identify $W$ with $\kappa_T(W).$

 It follows from \cite{HelfferSjostrand1986}, Section 9, that for
$u\in H(\Lambda;\langle\alpha_\xi\rangle^m)$ we have
$$\|Tu\|_{H_{\Phi_t}(W)}\leq c\|u\|_{tG}.$$ Let $W'\subset\subset
W$ and $W_1\subset\subset W'$ be another \neigh{} of $a_1$ which is
identified with $\{\alpha\in\Lambda_{tG};\,\,\Re\alpha\in W'\}.$ One
can show, using methods from \cite{HelfferSjostrand1986},
\cite{Sjostrand1982}, that for $0\leq t <t_0$ with $t_0$ small
enough $\exists$ $\delta_0>0$ such that
$$\| u\|_{tG,W'}\leq {\mathcal O}(1)\left( \|Tu\|_{H_{\Phi_t}( W)}
+e^{-\delta_0/h}\|u\|_{tG,\complement W'}\right),$$ where in general
the seminorm $\|u\|_{tG,W}$ is defined as in (\ref{normHS}) by
integration over $W\subset\Lambda_{tG}$ only.

  We define a new norm
uniformly equivalent to $\|\cdot\|_{tG}$ as $h\rightarrow 0$ by
$$\| u\|_{tG,\complement W'} +\| Tu\|_{H_{\Phi_t}(W)}.$$

\section{The FBI-Bargman transform of the quantum billiard operator $M.$ }\label{FBI-B}

 The operator $M$ was defined in the Introduction,
formula (\ref{Mdef}). Here we  give an alternative  definition (see
Burq \cite{Burq1997}).

 Let
$h=1/\lambda\in\rz_+.$
 If $H_{1,+}$ is the outgoing resolvent of the
problem
\begin{equation}\label{H1+}
\left\{\begin{array}{rl}
  (-h^2\Delta -1) H_{1,+} u &  =0,\,\,\mbox{in}\,\,\complement{\Omega}_1 \\
  H_{1,+}u_{|_{\partial\Omega_1}} & =u,  \\
\end{array}
\right.\end{equation} we consider  the outgoing resolvent of the
problem
\begin{equation}\label{H21}
\left\{\begin{array}{rl}
  (-h^2\Delta -1) H_{21} u &  =0,\,\,\mbox{in}\,\,\complement{\Omega}_2 \\
  H_{21}u_{|_{\partial\Omega_2}} & =H_{1,+} u_{|_{\partial\Omega_2}}.  \\
\end{array}
\right.\end{equation} The operator $M$ is defined as
\begin{equation}\label{Mdef2}
Mu=H_{21}u_{|\partial\Omega_1}.\end{equation} Some known facts about
$M$ are collected in  Appendix \ref{s-M}.

Let $W_1$ be the \neigh{} of $\rho_1=(a_1,0)\in T^*X$ independent of
$h$ and  defined in (\ref{W1}): $W_1=\Omega_+(N+1)\cap\Omega_-(N).$
Let $G$ be the global escape function as in Proposition \ref{l-G}.

\subsection{Properties of the kernel  of $M$ on $\complement W_1$ }

Let $T$  be the FBI transform of Helffer and Sj\"ostrand as in
(\ref{HeSjFBI}) with admissible phase function $\phi$ in a complex
\neigh{} $U\times V$ of $(\rho_1,a_1)$  and let $S$ be an
approximate inverse to $T$ as in (\ref{inverseHeSjFBI}). Let
$\mu,\rho\in T^*X$ be arbitrary real for the moment. As $M$ is
bounded operator
 $L^2(X)\mapsto L^2(X)$ and $T,S$ are bounded $L^2(X)\mapsto
 L^2(T^*X)$ respectively $L^2(T^*X)\mapsto L^2(X)$ then
 the distributional kernel $K$ of $TMS$ given by \begin{equation}\label{K1} TMSTu(\alpha,h)=\int_{{\rm neigh}\, (\rho_1)}
K(\alpha,\beta ;h)Tu(\beta)d\beta,\,\,\alpha\in{\rm neigh}\,
(\rho_1),\end{equation}
 is well defined.

 By the results on the  propagation of
Gevrey 3 singularities due to  Lebeau \cite{Lebeau1984} (see also
Burq \cite{Burq1997}) it is known that on the complement of the
${\rm graph}\,(\kappa):$
\begin{equation}\label{graph}
\{(\alpha,\beta)\in T^*X\times
T^*X,\,\,\alpha=\kappa(\beta)\},\end{equation} $K(\alpha,\beta ;h)$
is small:
\begin{equation}\label{Lebeau}
K(\alpha,\beta ;h)= {\mathcal O}_N(1) e^{-1/ch^{1/3}}\left(\max
(|\langle \alpha_\xi\rangle|,|\langle\beta_\xi\rangle|)\right)^{-N
}\,\,\mbox{if}\,\,\alpha\neq\kappa(\beta),\,\,\mbox{for any}\,\,
N>0,
\end{equation} if $\kappa(\rho)$ is defined, otherwise no condition.
Since $H(\alpha)$ has compact support and is small, we have the same
estimate for the reduced kernel: $e^{(-H(\alpha) +H(\beta))/h}
K(\alpha,\beta;h).$

In a \neigh{} $ U\times U$ of the graph of $\kappa$ (\ref{graph}),
from the general principles, as
  $M$
is bounded, we know only that $K(\alpha,\beta ;h)$ is of order
${\mathcal O}(h^{-N'})$ for some $N'>0.$ If
$\alpha,\beta\in\complement W_1$ we can get a better bound on the
reduced kernel $e^{(-H(\alpha) +H(\beta))/h} K(\alpha,\beta;h)$ by
using the properties of the escape function $G,$ Lemma \ref{l-G}: if
$\alpha=\kappa(\beta)\in T^*X\setminus W_1$  then
 $G(\alpha)-G(\beta)\geq C h\ln^2(1/h)$ and the bound is
still valid in the whole \neigh{} of $\{(\alpha,\beta)\in
T^*X\setminus W_1;\,\,\alpha=\kappa(\beta)\}.$

Let $\Lambda_{tG}$ be the IR-submanifold of the complex \neigh{} of
$T^*X$ defined in (\ref{Lambda}). Let $K(\alpha,\beta;h)$ denote
also the extension of (\ref{K1}) to the complex \neigh{} of
$(\rho_1,\rho_1)$ of the form $(U\cap\Lambda_{tG})\times
(U\cap\Lambda_{tG}).$ Let $H\in C^\infty(\Lambda_{tG};\rz)$ satisfy
$d_\alpha H_{|\Lambda_{tG}}=-\Im\alpha_\xi
d{\alpha_x}_{|\Lambda_{tG}}$ and be equal to zero for $|\alpha_\xi
|$ large enough. By (\ref{identification})  we have locally
$$H_t(\alpha_x,\alpha_\xi)=-\Re\alpha_\xi \cdot\Im\alpha_\xi
+tG(\Re\alpha_x,\Re\alpha_\xi).$$

Then for $(\alpha,\beta)$ in a \neigh{} of $\{(\alpha,\beta)\in
\complement W_1;\,\,\Re(\alpha)=\kappa(\Re\beta)\}$ we can arrange
that
\begin{equation}\label{onchar1}
 e^{-H(\alpha)/h}
K(\alpha,\beta;h)e^{H(\beta)/h}= e^{-\frac{1}{Ch}\min(
h\ln^2(1/h),|\alpha-\kappa(\beta)|^2)}\left(\max (|\langle
\alpha_\xi\rangle|,|\langle\beta_\xi\rangle|)\right)^{-N }
\end{equation}
for some $N_0 >0$ and for any $N>0.$

Then by the general theory of \cite{HelfferSjostrand1986} and
\cite{Sjostrand1996} the estimate (\ref{onchar1}) implies that for
$u\in H(\Lambda_{tG})$ we have the bound:
\begin{equation}\label{boundCW} \|Mu\|_{tG,\complement W'}\leq {\mathcal
O}(h^{N_0})\| u\|_{tG}\end{equation} for some $N_0 >0.$ Here $W'$
 is any $h-$independent \neigh{} of $\rho_1,$ $W_1\subset\subset W'\subset\subset W$
identified with $\{\alpha\in\Lambda_{tG};\,\,\Re\alpha\in W'\}$ and
in general the seminorm $\|u\|_{tG,W}$ is defined as in
(\ref{normHS}) by integration over $W\subset\Lambda_{tG}$ only (see
the end of Section \ref{ss-Bargman}). Bound (\ref{boundCW}) extends to complex $\lambda\in\Lambda_{A,B}$ with some
different power $N_0>0.$

\subsection{Local form of $M$}

Microlocally near $(a_1,0)\in T^*\partial\Omega_1,$ the asymptotic
solution to (\ref{H21}) is given by means of WKB construction.

In the hyperbolic zone ${\mathcal H}$, $M$  is approximated by
\fourior{} $H \in I_h^0(X,X; \kappa'),$  where
 $\kappa'=\{ (x,\xi;y,-\eta):\,\,(x,\xi)=\kappa(y,\eta)\},$  associated to the real
\canform{} of billiard $\kappa,$ with the real phase and can be
taken in the form (see \cite{Gerard1988}):
\begin{equation}\label{H}
Hu(u,\lambda)=\left(\frac{\lambda}{2\pi}\right)^n\int e^{-i\lambda(
s (x,\theta) -y\theta +2d)}b(x,y,\theta ;\lambda)u(y)dyd\theta,
\end{equation} where $s$ solves the eikonal equation,
$s (x,\theta) -y\theta$ parameterizes the \canform{} of billiard
$\kappa,$ $s(0,0)=0.$

If $\lambda$ is complex, $\lambda\in\Lambda_{A,B},$ then
 we denote $\lambda_1=\Re\lambda,$ $\lambda_2=\Im\lambda$ and
 write
\begin{equation*}
H u(x,\lambda)=\left(\frac{\lambda_1}{2\pi}\right)^n\int\!\!\int
e^{-i\lambda_1 (s(x,\theta)-y\theta
+2d)}\widetilde{b}(x,y,\theta;\lambda_1) u(y,\lambda) dy
d\theta\end{equation*} with
$$\tilde{b}(x,y,\theta;\lambda)=( 1+i\tan\arg\lambda)^{n} e^{\lambda_2
(s(x,\theta)-y\theta +2d)}b(x,y,\theta;\lambda_1).$$

Let $T: L^2(X)\mapsto  H_{\Phi_0},\,\, H(\Lambda_t)\mapsto
H_{\Phi_t} $ (modulo exponentially small errors) be as in
(\ref{A.23}) associated to the complex \canform{} $\kappa_T$ such
that $\kappa_T(T^*X)=\Lambda_{\Phi_0}$ and
$\kappa_T(\Lambda_{tG})\mapsto\Lambda_{\Phi_t}.$

 For $\Phi$ denoting either $\Phi_0$ or
$\Phi_t,$ the \fourior{} $\hat{H}:H_\Phi\mapsto H_\Phi$ is
associated to the real \canform{} $\hat{\kappa}$ in the sense that
locally $\hat{\kappa}(\Lambda_\Phi)=\Lambda_\Phi.$

Let $V_1=\pi_x\kappa_T W_1$ and for any $x\in V_1$ let
$\gamma(x)\in\Lambda_{tG}$ be such that
$\pi_x\kappa_T(\gamma(x))=x.$
 In
$W_1\subset{\mathcal H}$ we apply the Bargman transform and denote
$\hat{H}$ the transformed operator associated to the \canform{}
$\hat{\kappa}=\kappa_{T}\circ\kappa\circ\kappa_{T}^{-1},$ satisfying
for $u\in L^2(W_1),$
$$TMu=\hat{H}Tu +{\mathcal O}(h^{\infty})\|u\|\,\,\mbox{in}\,\,H_{\Phi}.$$

As $Hu$ is asymptotic solution to (\ref{H21})
 and obstacles are analytic then $\hat{H}$ can be taken analytic of the form
\begin{equation}\label{kertra}
\hat{H} u(x,\lambda)=\left(\frac{\lambda}{2\pi}\right)^n\int\!\!\int
e^{-i\lambda \phi(x,y,\theta)}a(x,y,\theta;\lambda) u(y,\lambda)
d\theta dy, \,\,u\in H_{\Phi}.\end{equation}  with analytic phase and
an amplitude which is an analytic symbol of order $0$ realized with
some suitable contour $\Gamma (x)\subset W_1$ passing through the
critical point, and introducing some cut-off functions.

 It is well-known that, by
choosing an appropriate integration contour, $\hat{H}$ is associated
with a kernel $K_{\hat{H}}\in H_{\Phi(x)+\Phi(\overline{y})}$ such
that, formally:
$$\hat{H}u(x,h)=\int
K_{\hat{H}}(x,\overline{y},h)u(y,h)e^{-2\Phi(y)/h} L(dy),\,\,u\in
H_\Phi.$$  The kernel $K_{\hat{H}}$ is uniquely defined as an
element of  $H_{\Phi(x)+\Phi(\overline{y})}$  through the data of
$\hat{H},$ but not as a function.

Denote $$k(x,y;h)=K_{\hat{H}}(x,\overline{y},h)e^{-2\Phi(y)/h}.$$ If
$\Re\gamma(x)\neq\kappa(\Re\gamma(y))\not\in W_0$  then by the
results on the propagation of analytic singularities of Lebeau
\cite{Lebeau1984} (see also Burq \cite{Burq1997}) we have
\begin{equation}\label{estoff}
e^{-\Phi(x)/h}k(x,y;h)e^{\Phi(y)/h}={\mathcal O}\left(e^{-\epsilon
/h}\right),\end{equation} where $\Phi$ is either $\Phi_0$ or
$\Phi_t=\Phi_0(x) +tG(x) +{\mathcal O}(t^2|x|^2),$ where
$G(x):=G(\Re\gamma(x)),$ $\pi_x\kappa_T(\gamma(x))=x.$ Here we used
that $G={\mathcal O}(h\ln(1/h)).$

If $\Re\gamma(x)=\kappa(\Re\gamma(y))$ then we use the properties of
the escape function:
$$-G(\kappa(\Re\gamma(y))) +G(\Re\gamma(y))\leq -C h\ln (1/h),\,\,c>0,\,\,
\mbox{if}\,\, \Re\gamma\in W_1\setminus W_0.$$

Choosing $W_0$ large enough ($c_0$ in the definition (\ref{W0}) is
sufficiently large) we can arrange that
\begin{equation}\label{eston}e^{-\Phi_t(x)/h}k(x,y;h)e^{\Phi_t(y)/h}={\mathcal O}\left(h^{N_0}
\right)\,\,\mbox{if}\,\,\Re\gamma(x)=\kappa(\Re\gamma(y))\cap
V_0=\emptyset,\end{equation} for some $N_0>0$ and $V_0=\pi_x\kappa_T
W_0.$

\subsection{The Bargman transform  reducing the stable manifolds to $x=0,$ $\xi
=0$}\label{s-intr-tr}
 The outgoing $\Lambda_+$, incoming $\Lambda_-$ stable manifolds for $\kappa$
  are lagrangian, intersecting transversally at
$(0,0).$ We can introduce real symplectic coordinates $(x,\xi)$ such
that $\Lambda_+:\,\,\xi =0,\,\,\Lambda_-:\,\, x=0.$

 We consider the image of
$\Lambda_\pm$ under the complex canonical transformation
$\kappa_T:\,\,T^*X\mapsto \Lambda_{\Phi_0}$ associated to the
Bargman transform $T.$ Transformation $\kappa_T$ preserves the
properties of $\Lambda_\pm.$ We will often write $\Lambda_\pm$
instead of $\kappa_T(\Lambda_\pm)$.

We know that the fiber $\{ x=0\}$ is a Lagrangian manifold, which is
strictly negative with respect to $\Lambda_{\Phi_0}.$
 It follows that $\{x=0\}$ and $\kappa(\Lambda_-)$ are transversal.
 Then $\Lambda_\pm$ are given by
$\xi=\partial_x\phi_\pm,$ where $\phi_\pm$ are holomorphic, and
$-\Im\phi_+-\Phi_0(x)\leq 0.$ The image of $\Lambda_{tG}$ is of the
form $\Lambda_{\Phi_t},$ where $\Phi_t(x)=\Phi_0(x) +tG(x)
+{\mathcal O}(t^2|x|^2),$ where $G$ is considered also as a function
on $\Lambda_{\Phi_0}\simeq \cz_x^n$ (see (\ref{G(x)})):
 \begin{equation}\label{GX}
G(x):=G\circ(\pi_x\circ\kappa_{T|T^*X})^{-1}=G(y,\eta)\,\,
\mbox{with}\,\, (y,\eta)\in\Lambda_{tG}\,\, \mbox{given by}\,\,
\pi_x(\kappa_T(y,\eta))=x.\end{equation}

(If $\phi=i(x-y)^2/2,$ and $(y,\eta)\in\rz^{2n}$ (instead of complex
$\Lambda_{tG}$), then $x=y-i\eta$ and $G(x)=G(\Re x,-\Im x).$)

 The
strict positivity of $\Lambda_+$ with respect to $\Lambda_{tG}$ then
implies that
$$\Im\phi_+ +\Phi_t\sim |x|^2.$$ Using the strict negativity of
$\Lambda_-$ with respect to $\Lambda_{tG}$ we have (see
\cite{Sjostrand1986})
\begin{lemma}\label{l-Sj}There is totally real linear space, $L\subset\cz^n,$
of real dimension $n,$ such that,
$$\Phi_t +\Im\phi_-\sim -|x|^2\,\,\mbox{on}\,\, L.$$
\end{lemma}

Considering still the situation after application of $\kappa_T,$ we
now introduce a new \canform{} $\kappa_F$ which maps $\Lambda_+$ to
$\{\xi =0\}$ and $\Lambda_-$ to $\{ x=0\}.$ Then it is easy to see
that $\kappa_F$ is given by
$\kappa_F:\,\,(y,-\phi_y'(x,y))\mapsto (x,\phi_x'(x,y)),$ where the generating
function $\phi(x,y)$ verifies: $\det\phi_{x,y}''\neq 0,\,\,\phi
(0,y)=-\phi_-(y).$

Let $f(x,y;h)$ be classical analytic symbol defined near
$(x,y)=(0,0).$ Using the Lemma \ref{l-Sj}, we see that if $u\in
H_{\Phi_t },$
then we can define $Fu\in H_{\hat{\Phi}_t },$
by
choosing a nice contour for the integral expression,
$$ Fu(x;h)=\int e^{i\phi (x,y)/h} f(x,y;h) u(y)dy.$$
Here, we use terminology of \cite{Sjostrand1982}, $\hat{\Phi}_t$ is
a new strictly pl.s.h. function, determined up to a constant by the
relation $\Lambda_{\hat{\Phi}_t}=\kappa_F(\Lambda_{\Phi_t}).$

As $\Lambda_+$ is strictly positive \wrt{} $\Lambda_{\Phi_t}$ we get
that $\kappa_F(\Lambda_+)$ is strictly positive \wrt{}
$\Lambda_{\hat{\Phi}_t}$ and
  the strict positivity of $\kappa_F(\Lambda_+):\,\,\xi =0$ (that we
shall from now on denote simply by ''$\Lambda_+$'') means
that,$$\hat{\Phi}_t\sim |x|^2.$$

Up to exponentially small errors (modulo equivalence in the spaces
$H_{\Phi_t},$ $H_{\hat{\Phi}_t},$) we can invert $F$ by an operator,
$$Gv(y;h)=\int e^{-i\phi (x,y)/h} g(x,y;h)v(x) dx.$$

The \canform{} of billiard $\kappa$   is transformed to
$$\kappa\equiv\exp H_p,\,\,{\rm D}\kappa (0,0)=\begin{pmatrix}
  A & 0 \\
  0& (A^T)^{-1}
\end{pmatrix},\,\,p=B(x,\xi)x\xi,\,\,B(0,0)={\rm
diag}(\mu_1,\ldots,\mu_n),$$ where we have used the Lewis-Sternberg
theorem (see {\cite{IantchenkoSjostrand2002}) and $\equiv$ denote
the equivalence relation for formal Taylor series at $(0,0).$

We denote $\widehat{W}_i:=\kappa_F(W_i)$   and
 $\hat{V}_i:=\pi(\widehat{W}_i)$ for $i=0,1.$

Composing the Bargman transform $T$ with the integral transform $F$
we get an operator of norm ${\mathcal O}(1),$
\begin{equation}\label{4.1}
FT:\,\,H(\Lambda_{tG})\mapsto H_\Phi  (V_1).\end{equation} Here
$t>0$ is small and fixed, $G$ is the escape function introduced in
Section \ref{s-G}, $\Phi$ denotes a function having all the
properties of $\hat{\Phi}_t$ and $V_1$ is small open \neigh{} of $0$
in $\cz^n.$

The direct definition of $T$ and $F$ only gives that
$\overline{\partial} (FT)$ is exponentially small, but we can
correct this by solving a $\overline{\partial}$-problem, using the
fact that $\Phi$ is strictly plurisubharmonic.

Composed operator $FT$ has microlocal inverse $SG$ of norm
${\mathcal O}(1):\,\, H_\Phi (V_1)\mapsto H(\Lambda_{tG})$ with the
properties:
\begin{itemize}
\item $SGFT$ is a \pseudor{} of order $0$ adopted to $\Lambda_{tG}$
in the sense of \cite{HelfferSjostrand1986}, which has compactly
supported symbol and which realizes the identity microlocally near
$\rho_1=(a_1,0).$
\item We have $\|FTSG u-u\|_{H_\Phi(\tilde{V}_1)}={\mathcal O}
(h^\infty)\| u\|_{H_\Phi(V_1)},$ where $\tilde{V}_1\subset\subset
V_1$ has the same properties as $V_1.$
\end{itemize}
We have
\begin{theorem}\label{Th.2.6} The transformed operator $FTHSG=e^{-i\lambda 2d}
M_1$ is an analytic \fourior{} given formally by
\begin{equation}\label{formofM1}
M_1 u(x,\lambda)=\left(\frac{\lambda}{2\pi}\right)^n\int
e^{-i\lambda (\varphi (x,\theta) -y\cdot\theta)}
b(x,y,\theta;\lambda ) u(y,\lambda )dy d\theta.\end{equation} Here
$b\in  S^{0,0}_{\rm cla},$   $b(0,0;\lambda
)=(\nu_1\cdot\ldots\cdot\nu_n )^{-1/2} +O(|\lambda |^{-1}),$ where
$\nu_1,\ldots,\nu_n$ are the eigenvalues $>1$ of ${\rm D}\kappa
(0,0).$ We have $\varphi(x,\theta)=A^{-1}x\cdot\theta +{\mathcal
O}(|(x,\theta)|^3),$ where  $A$ is diagonal with eigenvalues
$\nu_1,\ldots,\nu_n,$ $\nu_i=e^{\mu_i}.$

$M_1$ can be realized as bounded operator $
H_{\Phi}(\Omega_1)\mapsto H_{\Phi}(\Omega_2)$ with domains
$\Omega_2\subset\subset\Omega_1\subset\subset V_1,$ where we
identify $\Omega_1$ with $\pi\kappa \left(
\Lambda_\Phi\cap\pi^{-1}\Omega_1\right).$

  For any $\chi\in C_0^\infty(\cz^n)$ with
$\supp\chi\subset V_1\setminus V_0$ and some $N_0 >0,$ we have
\begin{equation}\label{est}\|\chi M_1 u\|_{L^2_{\Phi}(\Omega_2)}={\mathcal
O}(h^{N_0})\|u\|_{L^2_{\Phi}(\Omega_1)}.\end{equation}

\end{theorem}
 The form of $M_1$ in (\ref{formofM1}) follows as in
 \cite{Gerard1988}.
 The estimate (\ref{est}) follows from (\ref{estoff}) and (\ref{eston}).

Let $\Phi_t$ denote the new $\hat{\Phi}_t$ after the transformation
$F.$

\section{Semiclassical quantum Birkhoff normal form.}\label{s-qbnf}  In this section we perform an analytic
Birkhoff transform $\kappa_r$ up to some fix order $r,$
$$\kappa_r:\,\,\Lambda_{\Phi_t}\mapsto\Lambda_{\hat{\Phi}_t}.$$
Let  $\lambda\in\Lambda_{A,B}$ and we put  $h=1/\Re\lambda.$ We
need the following notion of equivalence used in papers
\cite{IantchenkoSjostrand2002},
\cite{IantchenkoSjostrandZworski2002}.
 Here  $\rho$ denotes either
pair $(x,\xi)$ or triple $(x,y,\theta).$
\begin{definition}\label{def-equiv}
 Let $U\in I^0(X,X,\kappa ')$
and $\tilde{U}\in I^0(X,X,\tilde{\kappa}')$ be two \fourior{}s. Then
$U\equiv\tilde{U}$ to the order $r$ if $\kappa
(0,0)=\tilde{\kappa}(0,0)=(0,0),$ $\kappa,$ $\tilde{\kappa}$ agree
to the order $\rho^{2r+1}$ at $(0,0)$, the terms with number $j,$
with $0\leq j\leq r,$ in the asymptotic expansions (in powers of
$h$, corresponding to $h^j$) of $\phi,\tilde{\phi},a,\tilde{a}$
(phase, respectively, amplitudes) agree to the order
$\rho^{2(r-j)+1}$ at $(0,0,0).$
\end{definition}
We write $a\equiv{\mathcal O}(|\rho|^{2r+1}+h^{r+1})$ to denote that
$a$ is equivalent  to zero to the order $r,$ and use $\equiv$ to
denote the equivalence to infinite order.

 The theorem proven in \cite{IantchenkoSjostrand2002} for real hyperbolic $\kappa$ says the following:
\begin{theorem}\label{th-log} Let $M_1$ be a \fourior{} as in Theorem
\ref{Th.2.6} which quantizes an analytic \canform{} $\kappa:\,\,{\rm
neigh\,}(0,{\cz }^{2n})\to {\rm neigh\,}(0,{\cz }^{2n}).$ Assume
that the eigenvalues of $d\kappa(0)$ satisfy
$$0< \nu_n^{-1}\leq\ldots\leq\nu_1^{-1}
<1<\nu_1\leq\ldots\leq\nu_n.$$ Let $\mu_j=\log\nu_j.$

Then there exists a
 \pseudor{}
$P(x,\lambda^{-1}D_x;h)$ with symbol $P(\rho;h)\sim p(\rho)
+hp_1(\rho)+\ldots,$ such that
\begin{equation}\label{3.23}
M_1\equiv e^{-i\lambda P}. \end{equation} $P$  is uniquely
determined modulo ``$\equiv$'' and up to an integer multiple of
$\frac{2\pi}{\lambda}$ by (\ref{3.23}) and the choice of $p_0$ such
that $p(\rho)-p_0(\rho)={\mathcal O}(|\rho|^3).$

Suppose \begin{equation}\label{non-res} \sum_1^n k_j\mu_j=0,\,\,
k_j\in\zz\,\,\Longrightarrow\,\, k_1=\ldots=k_n=0. \end{equation}

Then there exists elliptic \fourior{} $B$ and a classical symbol of
order $0,$
\begin{align*}
&F(\imath ; 1/\lambda) =F_0(\imath) +\lambda^{-1} F_1(\imath)
+\lambda^{-2}
F_2(\imath)+\ldots,\,\,\imath=(\imath_1,\ldots,\imath_n),\,\,\imath_j=\xi_j
x_j,\\
&F_0(\imath)=\sum_{j=1}^n\mu_j\imath_j +R(\imath),\,\, R(\imath)
={\mathcal O} (\imath^2),
\end{align*}
such that $$P\equiv B^{-1}F(I;1/\lambda)B,\,\,M_1\equiv B^{-1}
e^{-i\lambda F(I; 1/\lambda )}
B,\,\,I=(I_1,\ldots,I_n),\,\,I_j=\frac{1}{2i\lambda}(x_j\partial_{x_j}
+\partial_{x_j} x_j).$$ When $M_1$ is microlocally unitary near
$(0,0),$ $F(I; 1/\lambda)$ can be chosen to be microlocally
self-adjoint and unitary, respectively, at $(0,0,0).$

\end{theorem}

Fix some order $r\geq 1$ and let
$$F^r(\imath ;h):=F_0^r(\imath) +\lambda^{-1} F_1^{r-1}(\imath)
+\lambda^{-2} F_2^{r-2}(\imath)+\ldots +\lambda^{-r} F_r^0,$$ where
$F_j^{r-j}$ is polynomial of degree $(r-j)$ in $\imath$ ($F_r^0$ is
constant) such that
$$F_j(\imath) -F_j^{r-j}(\imath)={\mathcal O}(\imath^{r-j +1}),$$
($F_j^{r-j}(\xi\cdot x)$ is resonant up to order $2(r-j)$). For any
$r'>r$ the terms  in expansion of $F^{r'}(\imath ;1/\lambda)$ and
$F^{r'}$ with the same indexes coincide, thus the coefficients in
expansion $F^r(\imath ; 1/\lambda)$ are independent of $r.$

 Then
 $F^r(I; 1/\lambda)$ is {\em analytic} \pseudor{}. From the
construction of Birkhoff transform $B$ for the {\em Hamiltonian} it
is easy to see that there exists an {\em analytic} \fourior{} $B_r$
associated to an analytic $\kappa_r$ such that
$$P-B_r^{-1}F^r(I;1/\lambda)B_r=Q^r,$$ where $Q^r$ is \pseudor{} with the
symbol $q^r(\rho;1/\lambda)\sim q_0(\rho) +\lambda^{-1}
q_1(\rho)+\ldots,$ with $q_j(\rho)={\mathcal O}(|\rho|^{2(r-j)
+1}+h^{r-j+1}),$ if $0\leq j\leq r$ and $q_j=p_j$ for $j\geq r+1.$
Here $B_r^{-1}$ is analytic \fourior{} satisfying
$B_r^{-1}B_r-I=\tilde{Q}^r$ and $\tilde{Q}^r$ of the same type as
$Q_r.$

Operator $e^{-i\lambda F^r(I; 1/\lambda)}$ is analytic \fourior{}
which satisfies
$$M_1 - B_r^{-1} e^{-i\lambda F^r(I; 1/\lambda )} B_r\equiv{\mathcal
O}(|\rho|^{2r+1}+ h^{r+1}).$$

We denote $$M_0:=B_rM_1B_r^{-1},\,\,R^r:=M_0-e^{-i\lambda F^r(I;
1/\lambda)}.$$

Operators $M_0$ and $R^r$  have the form (\ref{formofM1}). Let
$$R^ru(y,h)=\lambda^{n}\int\!\!\int e^{-i\lambda\varphi^r
(x,y,\theta)} b^r(x,y,\theta,\lambda)dyd\theta$$ and
$b\sim\sum\lambda^j b_j,$ $b_j={\mathcal O}((|x|+|y|+|\theta
|)^{2(r-j)_++1}).$ The integration contour is chosen such that
$$-\Phi (x) +\Im\varphi^r(x,y,\theta)+\Phi (y)\sim -|x|^2
-|y|^2-|\theta |^2.$$
 We have
\begin{equation}\label{Rr}
e^{-\Phi/h}R^re^{\Phi/h}\equiv{\mathcal O}(|\rho|^{2r+1}+
h^{r+1}).\end{equation}

 Relation extends to
$\lambda\in \Lambda_{A,B}$  with the bound ${\mathcal O}\left(
h^{-C}(|\rho|^{2r+1}+ h^{r+1})\right)$ for some $C>0$ and usual
convention  $h=1/\lambda_1 ,$ $\lambda_1=\Re\lambda.$

Using (\ref{Rr}) and that for $\lambda\in\Lambda_{A,B},$ we have
$|e^{-i2d\lambda}|\leq {\mathcal O}(1)\lambda_1^m$ for some
(positive) $m,$ we get
\begin{lemma}\label{l-R}
Let $W_0=B(0,c_{0}\sqrt{s}),$ $s=h \ln(1/h),$ $c_{0}>0$ as in
(\ref{W0}) and $V_0$ such that $\pi(W_0)=V_0.$ Suppose
$\lambda\in\Lambda_{A,B}.$ For any $N\in\nz$ there is
$r=r(N)\in\nz$ such that
$$\|e^{-i2d\lambda}R^ru\|_{L^2_{\Phi}(\Omega_2)}\leq {\mathcal O}_r(
h^{N})\|u\|_{L_{\Phi}^2(\Omega_1)},$$ where
$\Omega_2\subset\subset\Omega_1\subset\subset V_0$ and $\Phi$ is
either $\Phi_0$ or $\Phi_t$ or any pl.s.h. function close to
$\Phi_0$ in $C^2.$
\end{lemma}

\underline{\bf Proof:} In the proof $h=1/\lambda$ is real, extension
to $\lambda\in\Lambda_{A,B}$ is straightforward.

We write  $\rho =(x,y,\theta),$
$\phi(\rho)=-(\varphi(x,\theta)-y\theta).$ Stationary point of
$\phi$ \wrt{} $y,\theta$ is given by $(y,\theta) =(0,0).$
 If $M_0 u=v$ we have,
\begin{equation}\label{vu}
v(x)=h^{-n}\int\!\!\int_{\Gamma(x)}
e^{i(\phi_r(\rho)+\psi_r(\rho))/h}
\left(b^r(\rho,\lambda)+c^r(\rho,\lambda)\right)u(y)dyd\theta,
\end{equation} with a good integration contour $\Gamma(x)$ passing through the critical point and $\psi_r(\rho)={\mathcal
O} (|\rho |^{2r+1}),$  $c^r\sim\sum_{j=0}^\infty\lambda^{-j}
c_j^r$ with, for $0\leq j\leq r,$ $c_j^r={\mathcal O} (|\rho
|^{2(r-j)+1}).$ If $|\rho |\leq {\mathcal O}(s^{1/2})$ we get
$$\psi_r(\rho)={\mathcal
O} (s^{r+1/2}),\,\,c^r(\rho)={\mathcal O} (s^{r+1/2}). $$

Let
$\Omega_2\subset B(x_0, c_2s^{1/2}),$ be a \neigh{} of $x_0=0,$    and
let
$$\Omega_1=\pi_y\kappa^{-1}(\pi_x^{-1}\Omega_2\cap\Lambda_\Phi)\subset\subset B(y_0, c_1s^{1/2}),\,\,y_0=0.$$

By modifying $M_0u$ with ${\mathcal O}(h^{N_0})$ in
$L^2_{\Phi}(\Omega_2),$ with some $N_0>0,$ we can restrict the
integration contour $\Gamma(x)$ in (\ref{vu}) to $|y| +|\theta|\leq
\epsilon s^{1/2},$ with $\epsilon$ small enough.

 It is clear that $ e^{-i\lambda F^r(I;1/\lambda)}$ is well defined
as a bounded operator
$H_{\Phi}(\Omega_1)\mapsto
H_{\Phi}(\Omega_2).$

We want to estimate the difference $M_0u-e^{-i\lambda
F^r(I;1/\lambda)}u$ on $\Omega_2.$

 We make substitution
$x=s^{1/2}\tilde{x},\,\,y=s^{1/2}\tilde{y},\,\,\theta=s^{1/2}\tilde{\theta}$
and write $u(y)=\tilde{u}(\tilde{y}),$ $v(x)=\tilde{v}(\tilde{x}),$
$\tilde{\rho}=(\tilde{x},\tilde{y},\tilde{\theta})$ and
$s^{1/2}\tilde{\rho}=(s^{1/2}\tilde{x},s^{1/2}\tilde{y},s^{1/2}\tilde{\theta}).$
Then, if $M_0  u=v$ we have
\begin{equation*}\tilde{v}(\tilde{x})=\left( \frac{h}{s}\right)^{-n}\int\!\!\int
e^{\frac{i}{h}(\phi_r(s^{1/2}\tilde{\rho})+{\mathcal O}(s^{r+1/2}))}
\left( b^r(s^{1/2}\tilde{\rho},\lambda)+{\mathcal
O}(s^{r+1/2})\right)\tilde{u}(\tilde{y})d\tilde{y}d\tilde{\theta}
\end{equation*}

 Let
$\tilde{\Gamma}$ be the image of $\Gamma(x)$ by substitution. We
have then $\tilde{x}\in\tilde{\Omega}_2\subset B(0,c_2),$ and for
$(\tilde{y},\tilde{\theta})\in\tilde{\Gamma}$ we have
$|\tilde{y}|+|\tilde{\theta}|
\leq \epsilon.$
 We write
$\phi_r^s:=\frac{1}{s}\phi_r(s^{1/2}\rho).$ Then on $\tilde{\Gamma}$
$$-s^{-1}\Phi(s^{1/2}\tilde{x})-\Im\phi_r^s
+s^{-1}\Phi(s^{1/2}\tilde{x})\sim-\left(|\tilde{x}|^2+|\tilde{y}|^2+|\tilde{\theta}|^2\right)
+{\mathcal O}(s^{r+1/2}h^{-1}).$$ Then, if
\begin{equation}\label{2.12}
v_r(\tilde{x})=\left( \frac{h}{s}\right)^{-n}\int_{\tilde{\Gamma}}
e^{\frac{s}{h}\phi_r^h}b_0^r(s^{1/2}\tilde{x},s^{1/2}\tilde{y},s^{1/2}\tilde{\theta})\tilde{u}(\tilde{y})
d\tilde{y}d\tilde{\theta},
\end{equation} we get
$$\|\tilde{v}-v_r\|_{\Phi,\tilde{\Omega}_2}\leq c(\ln(1/h))^{n}
s^{r+1/2}h^{-1}\|\tilde{u}\|_{\Phi,\tilde{\Omega}_1+B(0, \epsilon)}=
{\mathcal O}(h^{N})\|\tilde{u}\|_{\Phi,\tilde{\Omega}_1+B(0,
\epsilon)},$$  where $\|.\|_{\Phi,\Omega}$ denotes the norm in
$H_\Phi(\Omega).$ Here $N >0$ can be chosen arbitrary large by
choosing $r$ large enough.

 \hfill\qed

\section{Deformation of the escape function in a {\neigh} of $(0,0)$ and final definition of  the
space.
}\label{s-deformation} We write $h=1/\lambda_1,$
$\lambda_1=\Re\lambda.$ Let $W_0$ be the \neigh{} of $0$ of the size
$c_0\sqrt{h\ln (1/h)}$ as in (\ref{W0}) (with appropriate
identifications of the domains). Let $W_1\supset\supset W_0$ be the
$\lambda-$independent domain as in (\ref{W1}).

Operators  $M_0=B_rM_1B_r^{-1},\,\,e^{-i\lambda F^r(I; 1/\lambda
)}=M_0-R^r$ can be extended to the whole $\lambda$-independent
\neigh{} of $0.$ Let $U=B_rF T,$ where $T$ is the Bargman transform  with a global choice of phase (\ref{A.23}), $F$ is the
analytic \fourior{} quantizing  the symplectic change of coordinates
near $a_1=0,$ $B_r$ - the Birkhoff transform up to order $r.$  On
$W_1\setminus W_0$ we can estimate the norm of $R^r$ as in (\ref{est}).
 Together with Lemma \ref{l-R} it implies
\begin{lemma} For any $N\in\nz$  and $\lambda\in\Lambda_{A,B}$ there is
$r\in\nz$ such that
\begin{align*} & U M u= e^{-i\lambda 2 d} e^{-i\lambda F^r(I;1/\lambda)}
 U u
+e^{-i\lambda 2 d}R^ru ,\,\, \mbox{in}\,\, L_{\Phi}^2(V_1),\,\,\mbox{and}\\
&\|e^{-i\lambda 2 d}R^ru\|_{L^2_\Phi(V_1)}\leq {\mathcal
O}_r(h^N)\|u\|_{H(\Lambda_{tG})},\,\,\mbox{where}\,\,  \pi(W_1)=V_1.
\end{align*}

\end{lemma}

In Section \ref{s-G} we constructed   a global  escape function $G$ in
$T^*(\partial\Omega_1),$ such that $G(\kappa(\rho))-G(\rho) \geq 0$
 on $ T^*\partial(\Omega_1) $ and
such that we have strict inequality, outside an arbitrary small
neighborhood of  $(a_1,0)\in T^*\partial\Omega_1.$

Let  $\Omega_{\rm int}$ and $W$ be a $\lambda$-independent \neigh{}s
of $(a_1,0)$ (after usual identification of the domains),  such that
$W_{0}\subset\subset W\subset\subset\Omega_{\rm int}\subset\subset
W_1.$
 We define a  preliminary space ${\mathcal H}_{\rm pre}$ associated to the  Lagrangian space $\Lambda_{tG}$ as follows
$$u\in{\mathcal H}_{\rm
pre}\,\,\Leftrightarrow\,\,\left\{\begin{array}{c}
  T^{HS} u\in L^2_{ t G}(T^*\partial\Omega_1\setminus W) \\
  U u\in H_{\Phi_t}(\Omega_{\rm int}),\,\,\kappa_{U}(\Lambda_{tG})=\Lambda_{\Phi_t}. \\
\end{array}\right. $$

 In some small
\neigh{} of $(0,0),$ \begin{equation}\label{W00}
W_{00}=\{\rho,\,\,|\rho|\leq c_{00}\sqrt{h\ln(1/h)}\} \subset
\subset W_0,\end{equation} we deform $\Lambda_{\Phi_t}$ in the
following way.

We have
$\kappa_U:=\kappa_{B_r}\circ\kappa_F\circ\kappa_T:
\,\,\Lambda_{tG}\mapsto\Lambda_{t \hat{G}}.$ We denote again
$\hat{G}$ by $G.$
With $G$ we associate an IR-lagrangian manifold which is essentially
$$\Lambda_{tG}=\{(x,\xi)=\exp{(it H_G)}(y,\eta);\,\,(y,\eta)\in
T^*({\rm neigh}(0))\}
$$ for $0 < t\leq 1,$ small.

Locally we have $\Lambda _+=\{\xi =0\},$ $\Lambda_- =\{ x=0\}.$
 For $x^2+\xi^2\leq  {\mathcal O} (1) h\ln
(1/h)$ we have $G\sim (x^2-\xi^2)/2,$ and $\Lambda_{tG}$ is given by
$$\left\{\begin{array}{cc}
  x_j & =(\cos t) y_j-i(\sin t)\eta_j \\
  \xi_j & =-i(\sin t)y_j +(\cos t)\eta_j \\
\end{array}\right. ,\,\, y_j,\eta_j\in\rz,\,\,1\leq j\leq n.$$
Then
$\kappa_U(\Lambda_{tG})=\Lambda_{\Phi_t}=\{\xi=\frac{2}{i}\frac{\partial\Phi_t}{\partial
x}\},\,\,\Phi_t=\frac12(\cot t)(\Im x)^2 +\frac12(\tan t)(\Re
x)^2.$

 Notice that
$\Phi_{\pi/4}=|x|^2/2$ and that the corresponding IR-manifold is
given by $\xi=-i\overline{x}.$ Following \cite{KaidiKerdelhue2000}
we look for a new IR-manifold $\Lambda$ which coincides with
$\Lambda_{\pi /4}=\Lambda_{\Phi_{\pi/4}}$ in some neighborhood of
$(0,0)$ for $x^2 +\xi^2\leq {\mathcal O}(1) h\ln(1/h)$ and with $\Lambda_{\Phi_t}$
outside.

First  we notice that if $q_1(x)$ and $q_2(x)$ are strictly convex
quadratic forms  then we can find smooth strictly convex function
$\phi (x),$ with $\phi (x)=q_1(x)$ near $0$ and with $\phi (x)
=q_2(x)$ outside $V,$ where $V$ is any given \neigh{} of $0.$ Apply
this with
$V=\{ x\in\cz^n,\,\,|x| <K\}$ for some $K={\mathcal O}(\sqrt{h\ln(1/h)}),$ $q_1(x) =\Phi_{\pi
/4} (x),$ $q_2(x)=\Phi_t (x).$

 Then replace $\phi (x)$ by
$\phi_\sigma=\left(\frac{\sigma}{K}\right)^2\phi (\frac{xK}{\sigma}
),$ $0 <\sigma \ll 1,$ in order to decrease the \neigh{} of $0,$
where $\phi\neq \Phi_t,$ even further while keeping $\phi$ bounded
in $ C^2.$ We have $\sigma\sim \sqrt{h\ln (1/h)}$ and the
derivatives of higher order than $2$ are diverging as $h\rightarrow
0.$ Thus we get that the deformation $\phi_\sigma$ is strictly
convex function which is close to the original function in $C^2.$

The deformation of $\Lambda_{\Phi_t}$ is denoted
$\Lambda_{\tilde{\Phi}_t}$ and the new weight function satisfies
$$\tilde{\Phi}_t(x)=|x|^2/2,\,\,|x|^2\leq\epsilon h\ln (1/h),\,\,\tilde{\Phi}_t(x)=\Phi_t,\,
\,|x|^2\geq c_{00} h\ln (1/h), $$ where $0<\epsilon\ll c_{00}$ are
some constants.

The final version of the space is given by
$$u\in{\mathcal H}\,\,\Leftrightarrow\,\,\left\{\begin{array}{c}
  T^{HS} u\in L^2_{ t G}(T^*\partial\Omega_1\setminus W) \\
  U u\in H_{\tilde{\Phi}_t}(\Omega_{\rm int}). \\
\end{array}\right. $$

\section{Model problem.}\label{s-model}

In this section we derive expressions  for the coefficients in the
expansion in powers  $k^{-j}$ of the solution of equation
(\ref{modelproblem}).

We suppose that $\mu_j$ satisfy the non resonance condition
(\ref{non-res}). Thus all $\mu_j$ and $\nu_j=e^{\mu_j}$ are
different. We also have $\mu_j >0$ $\Rightarrow$ $\nu_j >1.$

For $\lambda\in\Lambda_{A,B},$ as in (\ref{domain}),  consider the
equation
$$2d\lambda +\lambda F^r \left(\frac{2\alpha
+1}{2i\lambda};1/\lambda\right)=2\pi
k,\,\,\alpha\in\zz_+^n,\,\,k\in\zz.$$ We omit $r$ in the notations
and proceed in increasing the order of generality.

{\bf A.} Assume $F(\imath;1/\lambda )=G(\imath,\mu)=\sum_1^n\mu_i\imath_i,$ $\imath=x\xi,$ is independent
of $\lambda$ and linear in $\imath.$
\begin{align*}
 &2d\lambda+\lambda G\left(\frac{1}{2i\lambda}(2\alpha_1
+1),\ldots,\frac{1}{2i\lambda}(2\alpha_n+1),\mu\right)=2\pi
k,\,\,k\in\zz\,\,\Rightarrow\\
&2d\lambda +\frac{1}{2i}\sum_1^n\mu_i(2\alpha_i +1)=2\pi
k\,\,\Rightarrow\,\,\lambda (k,\alpha)= k\frac{\pi}{d}
+\frac{i}{4d}\sum_1^n\mu_i(2\alpha_i+1),\,\,
(k,\alpha)\in\zz\times\nz^n.
\end{align*}
We have $${\rm Re}\,\lambda =\frac{\pi k}{d},\,\,{\rm Im}\,\lambda
=\frac{1}{4d}\sum_{j=1}^n\mu_j(2\alpha_j
+1)\,\,\Rightarrow\,\,|\alpha|\sim{\rm Im}\,\lambda,$$ as
$\Re\lambda\rightarrow\infty,$ $|\alpha|=O(|\Im\lambda |),$ $|\alpha
|=O(k^{1-\epsilon}),$ $|\alpha |= O(\log |k|),$ $k\sim\Re\lambda .$

{\bf B.}
Assume that $F(\imath)=G(\imath,\mu) +H(\imath),$ polynomial of
degree $r,$  is independent of $\lambda,$ where $G$ is as above and
$H(\imath)=O(\imath^2).$ We have
\begin{equation}\label{expan}
H\left( \frac{y}{\lambda}\right) =\sum_{j=2}^{r}
\lambda^{-j}h_j(y),\,\,|1/\lambda|\rightarrow 0,\,\,|y|\leq
|\lambda|^{1-\epsilon }, \end{equation} for some $\epsilon
>0,$ where $|h_j(y)|=O(|y|^{j}),$
$h_j(y)=\sum_{|\alpha|=j}\frac{1}{\alpha !}\partial^\alpha
H(0)y^\alpha $ is a homogeneous polynomial of degree $j.$ Then we
have
\begin{align}\label{maineq0}
&\frac{\lambda}{k}\left( 1
+\frac{1}{4id\lambda}\sum_{i=1}^n\mu_i(2\alpha_i
+1)+\frac{1}{2d}\sum_{j=2}^{r} \lambda^{-j} h_j(\alpha)\right)
=\frac{\pi}{d} ,
\end{align}
where $h_j(\alpha)=h_j\left(\frac{2\alpha +1}{2i}\right)={\mathcal
O}\left( |\alpha |^{j}\right)$ is homogeneous polynomial in
$\alpha.$

  It is clear that equation (\ref{maineq0}) has solution in the
  form
\begin{equation}\label{ansatz}\frac{\lambda}{k}=\frac{\pi}{d}
+\frac{a_1}{k}
+\frac{a_2}{k^2}+\ldots\in\Lambda_{A,B}.\end{equation}

Then (\ref{maineq0}) transforms to \begin{align*}
&\frac{a_1}{k}+\frac{a_2}{k^2}+\ldots
+\frac{1}{4id}\sum_{i=1}^n\mu_i(2\alpha_i
+1)+\frac{\lambda}{2dk}\sum_{j=2}^{r} \lambda^{-j} h_j(\alpha) =0 ,
\end{align*}
where the last term is of order ${\mathcal O}(|\alpha|^2/k^2).$

Thus we get
$$ a_1=-\frac{1}{i4d}\sum_{j=1}^n\mu_j
 (2\alpha_j
+1)={\mathcal O}(|\alpha |)={\mathcal O}(\ln k).$$

The next equation is

\begin{align*}
&\frac{a_2}{k^2}+\frac{a_3}{k^3}\ldots
+\frac{1}{2d}\left(\frac{h_2(\alpha)}{\lambda k}
+\frac{h_3(\alpha)}{\lambda^2
k}+\ldots+\frac{h_r(\alpha)}{\lambda^{r-1} k}\right)
=0\,\,\Leftrightarrow\\
&\left(\frac{\pi}{d}k+a_1+\frac{a_2}{k}+\ldots\right)\left(\frac{a_2}{k^2}+\frac{a_3}{k^3}+\ldots\right)
+\frac{1}{2d}\left(\frac{h_2(\alpha)}{ k}
+\frac{h_3(\alpha)}{\lambda
k}+\ldots+\frac{h_r(\alpha)}{\lambda^{r-2} k}\right) =0
\end{align*}
which implies
\begin{equation}\label{a2} a_2=-\frac{1}{2\pi } h_2 \left(\frac{2\alpha
+1}{2i}\right)={\mathcal O}(|\alpha|^2).\end{equation} With this
choice of coefficients we get
\begin{align*}
&\frac{\pi}{d}k\left(\frac{a_3}{k^3}+\frac{a_4}{k^4}+\ldots\right)+\left(a_1+\frac{a_2}{k}+\ldots\right)
\left(\frac{a_2}{k^2}+\frac{a_3}{k^3}+\ldots\right)
+\frac{1}{2d}\left(\frac{h_3(\alpha)}{\lambda
k}+\ldots+\frac{h_r(\alpha)}{\lambda^{r-2} k}\right) =0.
\end{align*}
Considering the coefficients for $k^{-2}$ we get
$$\frac{\pi}{d}a_3
+a_1a_2+\frac{h_3(\alpha)}{d\pi}=0\,\,\Rightarrow\,\,a_3=-\frac{d}{\pi}a_1a_2
-\frac{h_3(\alpha)}{\pi^2}.$$ We get
$$a_3=-\frac{1}{i8\pi^2}\cdot\sum_{i=1}^n\mu_i(2\alpha_i +1)
\cdot
h_2\left(\frac{2\alpha+1}{2i}\right)-\frac{1}{\pi^2}h_3\left(\frac{2\alpha+1}{2i}\right).$$

It is clear that in this way we get all $a_j,$ $j\leq r$ dependent
only on $\{ h_k\}_{k\leq j}$ and with this choice $\lambda_0/k=\pi/d
+k^{-1}a_1 +k^{-2}a_2 +\ldots+k^{-r}a_r$ satisfies

\begin{align}\label{maineq}
&\frac{\lambda_0}{k}\left( 1
+\frac{1}{4id\lambda_0}\sum_{i=1}^n\mu_i(2\alpha_i
+1)+\frac{1}{2d}\sum_{j=2}^{r} \lambda_0^{-j} h_j(\alpha)\right)
-\frac{\pi}{d}={\mathcal
O}\left(\frac{|\alpha|^{r+1}}{k^{r+1}}\right).
\end{align}

{\bf C.} We write $h=1/\lambda.$
 In the general case $F^r(\imath ;h)= F_0(\imath) +h
F_1(\imath )+\ldots +h^r F_r(\imath) =:F_0(\imath) + K(\imath ;h),$
with $F_0(\imath)$ as before, $F_j$ polynomial in $\imath$ of degree
$(r-j).$

Then $ F(0;h)/ h=K(0;h)/h = F_1(0) +h F_2(0) +\ldots +h^{r-1} F_r$
and we can decompose \begin{align*}\frac{F(hy;h)}{h}= &
\frac{K(0;h)}{h}
+ G(y,\mu) +\frac{H(hy)}{h} +\frac{K(hy;h) - K(0;h)}{h}=\\
=& \frac{K(0;h)}{h} + G(y,\mu) +\frac{H(hy)}{h} +\sum_{j=1}^{r-1}
h^j k_j(y),
\end{align*} where $k_j$ is a polynomial of degree $\leq j$ and $k_j(0)=0.$
We get
$$F\left(\frac{y}{\lambda};1/\lambda\right)=F(0;1/\lambda)
+\frac{1}{\lambda}G(y,\mu)+\sum_{j=2}^{r}\lambda^{-j}k_{j-1}(y).$$

  Combining this with (\ref{expan}) where $h_{j+1}$ is a
homogeneous polynomial of degree $j+1,$ we get

$$F\left( \frac{y}{\lambda} ;1/\lambda\right)= F(0;1/\lambda) +\frac{1}{\lambda}G(y,\mu)
+\sum_{j=2}^{r} \lambda^{-j} f_j (y),$$ where
 for $2\leq
j\leq r,$ $f_j(y)=h_j(y)+k_{j-1}(y)$ is a polynomial of degree $\leq
j.$
  We have for example
$$f_2 (y)= F_1' (0) y +\frac12 F_0''(0) y\cdot y.$$

We need to solve the equation
$$\frac{\lambda}{k}\left(
1+\frac{1}{4id\lambda}\sum_{i=1}^n\mu_i(2\alpha_i
+1)+\frac{F_1(0)}{2d\lambda}
+\frac{1}{2d}\sum_{j=2}^r\lambda^{-j}\left( F_j(0)
+f_j(\alpha)\right) \right)=\frac{\pi}{d},$$ where we write
$f_j(\alpha)=f_j((2\alpha+1)/2i)={\mathcal O}(|\alpha|^j).$

Then we can get all coefficients in the expansion
$\frac{\lambda_0}{k}=\frac{\pi}{d} +\frac{a_1}{k}
+\frac{a_2}{k^2}+\ldots+\frac{a_r}{k^r}$ such that
$$\frac{\lambda_0}{k}\left(
1+\frac{1}{4id\lambda_0}\sum_{i=1}^n\mu_i(2\alpha_i
+1)+\frac{F_1(0)}{2d\lambda_0}
+\frac{1}{2d}\sum_{j=2}^r\lambda_0^{-j}\left( F_j(0)
+f_j(\alpha)\right) \right)-\frac{\pi}{d}={\mathcal
O}\left(\frac{|\alpha|^{r+1}}{k^{r+1}}\right).$$ First we get
$$a_1=-\frac{1}{i4d}\sum_{i=1}^n\mu_i(2\alpha_i+1)-\frac{F_1(0)}{2d}.$$
Let $F_j(0)+f_j(\alpha)=q_j(\alpha)$ and fix $2\leq  m \leq r.$
Suppose that we have already chosen $a_1,\ldots,a_{m-1}$ such that
\begin{equation}\label{eqm-1}
\frac{\lambda_0}{k}\left(
1+\frac{1}{4id\lambda_0}\sum_{i=1}^n\mu_i(2\alpha_i
+1)+\frac{F_1(0)}{2d\lambda_0}
+\frac{1}{2d}\sum_{j=2}^r\lambda_0^{-j}q_j(\alpha)-\frac{\pi}{d}\right)={\mathcal
O}\left(\frac{|\alpha|^{m}}{k^{m}}\right).\end{equation}

The \lhs{} of (\ref{eqm-1}) is then equal to
$$\frac{a_2}{k^2}+\frac{a_3}{k^3}+\ldots+\frac{a_r}{k^r}+\frac{1}{2d}\left(\frac{q_2(\alpha)}{\lambda k}+
\frac{q_3(\alpha)}{\lambda_0^2k}+ \ldots
+\frac{q_r(\alpha)}{\lambda^{r-1}k}\right).$$

Denote
$$E_{m-1}(\lambda_0)=\frac{a_2}{k^2}+\frac{a_3}{k^3}+\ldots+\frac{a_{m-1}}{k^{m-1}}
+\frac{1}{2d}\left(\frac{q_2(\alpha)}{\lambda_0
k}+\frac{q_3(\alpha)}{\lambda_0^2k}+ \ldots
+\frac{q_{m-1}(\alpha)}{\lambda_0^{m-2}k}\right).$$

Then equation (\ref{eqm-1}) writes
\begin{equation}\label{123}
E_{m-1}(\lambda_0) +\sum_{j=m}^r\frac{a_j}{k^j}
+\frac{1}{2d}\sum_{j=m}^r\frac{q_j(\alpha)}{\lambda_0^{j-1}k}={\mathcal
O}\left(\frac{|\alpha|^m}{k^m}\right).\end{equation} Let
$f(a_1,a_2,\ldots,a_{m-1})$ be defined by
$$E_{m-1}(\lambda_0)\sim\frac{f(a_1,a_2,\ldots,a_{m-1})}{k^m}
+{\mathcal
O}\left(\frac{|\alpha|^{m+1}}{k^{m+1}}\right)\,\,\mbox{as}\,\,|\alpha|/k\rightarrow
0.$$ Then  equation (\ref{123}) writes
$$\frac{f(a_1,a_2,\ldots,a_{m-1})}{k^m}+{\mathcal
O}\left(\frac{|\alpha|^{m+1}}{k^{m+1}}\right)+\frac{a_m}{k^m}+{\mathcal
O}\left(\frac{|\alpha|^{m+1}}{k^{m+1}}\right)+\frac{1}{2d}\frac{q_m}{\lambda_0^{m-1}k}=
-\sum_{j=m+1}^r\frac{q_j(\alpha)}{\lambda_0^{j-1}k}.$$ We chose
$a_m$ such that
$$f+a_m+\frac{1}{2d}\left(\frac{d}{\pi}\right)^{m-1}q_m(\alpha)=0\,\,\Leftrightarrow\,\,
a_{m}=-f(a_1,a_2,\ldots,a_{m-1})-\frac{1}{2d}\left(\frac{d}{\pi}\right)^{m-1}q_m(\alpha).$$
With this choice $\lambda_0$ satisfies equation (\ref{eqm-1}) with
the \rhs{} ${\mathcal
O}\left(\frac{|\alpha|^{m+1}}{k^{m+1}}\right).$

We summarize in the following theorem.
\begin{theorem}\label{asymptotics} Let $\Lambda_{A,B}$ be defined in
(\ref{domain}). For  $r\in\zz_+$ let $F^r$ be defined in (\ref{F}).
 Then for any $(k,\alpha)\in\nz\times\nz^n$ such that
$|\alpha|={\mathcal O}(\ln k)$ and $k$ large, there exist functions
$a_j=a_j(\alpha )={\mathcal O}(|\alpha |^{j}),$ polynomial in
$\alpha$, $j=0,1,\ldots,r,$ such that if we denote
 \begin{equation}\label{solution}\lambda_r := k\left( a_{0}  + \frac{a_1}{k}
 +\frac{a_2}{k^2} +\ldots +\frac{a_{r+1}}{k^{r+1}}\right)\in\Lambda_{A,B},\,\,a_0=\frac{\pi}{d},\,\,a_1=-\frac{1}{i4d}\sum_{i=1}^n\mu_i(2\alpha_i+1)-\frac{F_1(0)}{2d}
 ,\ldots,\end{equation}
 then we have
$$2d\lambda_r +\lambda_r F^r\left(\frac{2\alpha +1}{2i\lambda_r}
;1/\lambda_r\right) -2\pi k={\mathcal O}\left(\frac{(\ln
k)^{r+2}}{k^{r+1}}\right)$$ and the coefficients $a_j,$ are
independent of $r$  for $j=1,\ldots,r+1.$
\end{theorem}
\begin{note}\label{note}
\begin{enumerate}
\item For any $r,k\in\nz$ large and $\alpha ={\mathcal O}(\ln k),$
$\lambda_0:=\lambda_r(\alpha,k)$ satisfies
$$1-e^{-i2d\lambda_0}K_\alpha^r(\lambda_0)={\mathcal
O}\left(\frac{(\ln k
)^{r+2}}{k^{r+1}}\right),\,\,\mbox{where}\,\,K_\alpha^r(\lambda):=e^{-i\lambda
F^r\left(\frac{2\alpha +1}{2i\lambda} ;\lambda\right)}.$$
  \item Using that operator $M$ is microlocally
unitary we know that $F_1(0)$ is real. \item For large $k$ and
$\alpha ={\mathcal O}(\ln k),$ we have $\Re\lambda (\alpha ,k+1)
-\Re\lambda (\alpha,k)=\pi/2 +{\mathcal O}(1/k).$ For any $\beta
={\mathcal O}(\ln k),$ $\beta\neq\alpha,$ we have
$$ \Im\lambda (\alpha, k)
-\Im\lambda (\beta ,k)=\frac{1}{2d}\mu\cdot (\alpha -\beta)
+{\mathcal O}\left(\frac{\ln^2 k}{k}\right),\,\,k>>1.$$ In order to
have good separation of strings as $k\rightarrow\infty$ (the first
term in the above formula is dominated over the second error term)
we need to impose the Diophantine condition:
\begin{equation}\label{dioph}
\alpha\neq\beta,\,\,|\alpha |,|\beta|\leq
m\,\,\Rightarrow\,\,|\mu\alpha -\mu\beta |\geq \frac{1}{C(D)} e^{-D
m},\,\,D>0.\end{equation} It would imply for any $\alpha\neq\beta,$
$$|\alpha |,|\beta|\leq c\ln k\,\,\Rightarrow\,\,|\mu\alpha -\mu\beta |\geq
\frac{1}{C(D)} e^{-D c\ln k}\geq \frac{1}{C(D)\, k^{cD}}>>{\mathcal
O}\left(\frac{\ln^2 k}{k}\right),\,\,cD <1.$$ As the next term is of
order $k^{-1}$ we have separation.

\end{enumerate}

\end{note}

\section{The first local Grushin problem in $W_{0}.$}\label{s-gr1}

\subsection{Notations}

In $\Omega_{\rm int}$ for any $N\in\nz$ we can chose $r\in\nz$ such that we have $$UMu=e^{-i\lambda 2d}e^{-i\lambda
F^r(I;1/\lambda)}Uu+ {\mathcal O}\left(h^{N}\right)|u|_{\mathcal
H}. $$

Monomials $x^\alpha$ are formal eigenfunctions for $e^{-i\lambda
F^r(I;1/\lambda)}.$ Let $K_\alpha$ denote the corresponding
eigenvalue: $e^{-i\lambda F^r(I;1/\lambda)}x^\alpha =K_\alpha
x^\alpha.$ We have
 \begin{align} &K_\alpha(\lambda)=e^{-i\lambda
F^r(\frac{2\alpha_1 +1}{2i\lambda},\ldots,\frac{2\alpha_n
+1}{2i\lambda};1/\lambda)}
=\nonumber\\
&=\exp\{-\frac{1}{2}\sum_{i=1}^n\mu_i(2\alpha_i +1)-iF_1(0)
-i\sum_{j=2}^r\lambda^{-j+1}\left( F_j(0)
+f_j(\alpha)\right)\},\label{kalpha}\\ &\mbox{ where $f_j(\alpha)$
is polynomial of degree $\leq j.$}\label{polynom}
\end{align}
 We can
have $K_\alpha (\lambda)=K_{\alpha'}(\lambda)$ for
$\alpha\neq\alpha'.$
We chose a value $\lambda_0\in\Lambda_{A,B}$ and a multi-index
$\alpha_0$ such that
\begin{equation}\label{lambda_0}
1-e^{-i2d\lambda_0}K_{\alpha_0}(\lambda_0)=0.\end{equation}
Let $\beta\in\nz^n$ be such that $K_\beta(\lambda_0)\neq
K_{\alpha_0}(\lambda_0).$

We suppose that the
Diophantine condition (\ref{dioph}) in Note \ref{note} in the
previous section be satisfied and write
$|e^{-i2d\lambda_0}
K_{\beta}(\lambda_0)|=e^{|\mu (\alpha_0 -\beta)|+r}.$
Then if $|\beta|,\,\,|\alpha_0| \leq
C\ln \Re\lambda$ the Diophantine condition  implies
\begin{equation}\label{delta}
 |\mu
(\alpha_0 -\beta)|\geq c(\Re\lambda)^{-\delta},\,\,r={ o}
\left((\Re\lambda)^{-\delta}\right),
\end{equation}
 with some $0<\delta <1.$
If $|\mu (\alpha_0 -\beta)| ={\mathcal O}\left(
(\Re\lambda)^{-\delta}\right)$ then  we have $$|1-e^{-i2d\lambda_0}
K_{\beta}(\lambda_0)|={\mathcal O}\left(
(\Re\lambda)^{-\delta}\right),\,\,0<\delta <1.
$$

Let $\lambda$ belong to a small \neigh{} $\Omega_{\lambda_0}$ of
$\lambda_0.$

We chose $\Omega_{\lambda_0}$ sufficiently small such that for all
$\beta\in\nz^n$ with $K_\beta(\lambda_0)\neq
K_{\alpha_0}(\lambda_0)$ we have
$$\lambda\in\Omega_{\lambda_0}\cap\Lambda_{A,B},\,\,\left|1-e^{-i2d\lambda}
K_\beta(\lambda)\right| \geq\epsilon_0\,\,\mbox{and}\,\,\left|
1-e^{-i2d\lambda}K_{\alpha_0}(\lambda_0)\right| \leq\epsilon_0,$$
with \begin{equation}\label{epsilon0}
\epsilon_0={\mathcal O}\left(
(\Re\lambda)^{-\delta}\right),\,\,0<\delta <1.
\end{equation}
For $\lambda\in \Omega_{\lambda_0}\cap\Lambda_{A,B}$ we have
$$\Im\lambda =\frac{1}{4d}\mu (2\alpha_0 +1)-\frac{\Im F_1(0)}{2d}
+{\mathcal O}\left((\Re\lambda)^{-\delta}\right),$$
where we have $\Im F_1(0)=0$ (see Note \ref{note}).
 Let
$$J=\{\alpha\in\nz^n;\,\,|\alpha|\leq C\ln\Re\lambda_0,\,\,K_\alpha
(\lambda_0)= K_{\alpha_0}(\lambda_0)\}\,\,\mbox{and}\,\, a={\rm
Card}\,J.$$ Then we have
\begin{equation}\label{bound2}a\leq{\mathcal
O}(1)\ln\Re\lambda_0.\end{equation} We can chose
$\Lambda_{A,B}(\epsilon_0)\subset\Omega_{\lambda_0}\cap\Lambda_{A,B}$
such that
\begin{equation}\label{bound1}\mbox{for}\,\,\lambda\in\Lambda_{A,B}(\epsilon_0),\,\,\alpha\in J,\,\, \left|
1-e^{-i2d\lambda}K_{\alpha}(\lambda)\right| \leq 2\epsilon_0.
\end{equation}

We will use the following notations:  $\lambda_1=\Re\lambda,$
$h=1/\lambda_1$ and $s=h\ln (1/h).$

\subsection{Monomials form an almost orthonormal base in $H_\Phi$}
For $\Omega\subset\subset\cz^n$ we denote
 $$H_\Phi(\Omega)=\{ u\in{\mathcal Hol}(\Omega)
\,\,\mbox{such that}\,\,\int_{\Omega}|u|^2 e^{-2\Phi/h}
L(dx)<\infty\},$$
 where $\Phi$ is the modified weight defined in Section \ref{s-deformation}  such that $\Phi=|x|^2/2$ for $x\in
B(0,\epsilon\sqrt{h\ln (1/h)}),$ $\epsilon < c_{00}.$

 Since the weight function $\Phi$ is fixed from
now on, we shall suppress it from our notations as much as possible.
 Normalized monomials
$$\varphi_\alpha(x) =c_\alpha h^{-n/2}(h^{-1/2}
x)^\alpha,\,\,c_\alpha:=(\pi^n\alpha !)^{-1/2}$$ form an orthonormal
base in $H_{\Phi_0}(\cz^n)$ for $\Phi_0=|x|^2/2.$

Following \cite{KaidiKerdelhue2000} we show that $\varphi_\alpha(x)$
are orthonormalized in  $H_\Phi(\Omega)$ with an error of order
$h^{N}$ for any $N\in\nz.$

\begin{lemma}\label{l-monomials}
For all $\Omega\subset B(0, c),$  where $c$ is some real number, we
have
$$\langle\varphi_\alpha
|\varphi_\beta\rangle_{H_\Phi(\Omega)}=\delta_{\alpha,\beta}
+{\mathcal O} (h^{N})$$ for $N >0,$ which can be taken arbitrary
large by choosing $c_{00}$ large enough, and ${\mathcal O}$ is
uniform in $\alpha,\beta,\Omega.$
\end{lemma}
\underline{\bf Proof:}
\begin{align*}&\int_\Omega\left(\Pi_{j=1}^n |x_j
|^{2\alpha_j}\right) e^{-2\Phi (x)/h} L(dx) =\Pi_{j=1}^n\int_\cz
|x_j|^{2\alpha_j} e^{-x_j^2/h} L(dx_j) -\\
&-\int_{|x|\geq\epsilon\sqrt{ h\ln(1/h)}}\left(\Pi_{j=1}^n
|x_j|^{2\alpha_j}\right) e^{-|x|^2/h} L(dx) +\int_{\{
x\in\Omega,\,\,|x|\geq\epsilon\sqrt{ h\ln (1/h)}\}}\left(\Pi_{j=1}^n
|x_j|^{2\alpha_j}\right) e^{-2\Phi
(x)/h} L(dx)=\\
& =I_1-I_2 +I_3\end{align*} and $\epsilon$ is such that $\phi
(x)=|x|^2/2$ for $|x|\leq\epsilon\sqrt{h\ln(1/h)}.$ The first term
$I_1$ is equal to $\pi^n\alpha ! h^{|\alpha | +n}.$ We need to show
that
$$\frac{|I_2| +|I_3|}{\pi^n\alpha ! h^{|\alpha | +n}} ={\mathcal
O} (h^{N})$$ for $h\rightarrow 0$ uniformly  in $\alpha.$ As
$\Phi\geq c|x|^2$ for some $c>0,$ it is enough to show the
following:
\begin{align*}
& c_\alpha^2\int_{\{ x\in\cz^n,\,\,|x|\geq\epsilon\sqrt{ h\ln
(1/h)}\}}\left(h^{-1/2} x\right)^{2\alpha} e^{-2c|x|^2/h} L(dx)=
c_\alpha^2\left(\frac{h}{c}\right)^nc^{-|\alpha|}\int_{y\geq\epsilon\sqrt{
c\ln
(1/h)}} y^{2\alpha} e^{-2|y|^2} L(dy)\\
&\leq c_\alpha^2\left(\frac{h}{c}\right)^n
c^{-|\alpha|}e^{-\epsilon^2 c\ln(1/h)}\underbrace{\int_{\cz^n}
e^{-|y|^2} y^{2\alpha} L(dy)}_{\pi^n\alpha !\equiv
c_\alpha^{-2}}={\mathcal O}(h^{N}),
\end{align*}
where we mean $x^\alpha =\Pi_{j=1}^n x_j^{\alpha_j}$ and make
change of variables $y=\sqrt{\frac{c}{h}}x.$

 Next we show approximate
orthogonality:
\begin{align*}&\langle\varphi_\alpha |\varphi_\beta\rangle_{H_\Phi(\Omega)}
=\int_{|x|\leq\epsilon\sqrt{ h\ln(1/h)}}
\varphi_\alpha\overline{\varphi}_\beta e^{-|x|^2/h} L(dx) +\\
&+ \int_{x\in\Omega,\,\,|x|\geq\epsilon\sqrt{ h\ln(1/h)}}
\varphi_\alpha\overline{\varphi}_\beta e^{-2\Phi(x)/h} L(dx) =I_1
+I_2,\end{align*} $I_1=0,$ $I_2={\mathcal O} (h^{N})$ by previous
calculus.\hfill\qed

\subsection{Grushin problem}
 We denote $V_{0}=\pi(W_{0}),$ where $W_{0}$ is the image by $\kappa_U$ of the domain defined in
 (\ref{W0}).

 For
$N\in\nz,$ we introduce the expansions:
$$\tau_N u(x)=\sum_{|\alpha| <N}(\alpha !)^{-1}(\partial_x^\alpha
u)(0) x^\alpha=\sum_{|\alpha |<N} (\alpha
!)^{-1}((h^{1/2}\partial_x)^\alpha u)(0)(h^{-1/2} x)^\alpha.$$ We
denote ${\sl F}^N:={\mathcal Im}\,\tau_N,$ ${\sl O}^N:={\mathcal
Ker}\,\tau_N,$ $M(N)=\dim ({\sl F}^N).$

 Suppose that $N\in\nz$ satisfy
\begin{equation}\label{N}
\max_{\alpha\in J}|\alpha |\leq N={\mathcal
O}(1)\ln\Re\lambda_0,
\end{equation}
where $\lambda_0$ is fixed as before satisfying (\ref{lambda_0}).

 In
addition, we need the expansions
$$\tau_{J} u(x)=
\sum_{\alpha\in J}(\alpha !)^{-1}(\partial_x^\alpha
u)(0) x^\alpha=\sum_{\alpha\in J} (\alpha
!)^{-1}((h^{1/2}\partial_x)^\alpha u)(0)(h^{-1/2} x)^\alpha.$$ We
denote ${\sl F}_{J}:={\mathcal
Im}\,\tau_{J},$ ${\sl O}_{J}:={\mathcal
Ker}\,\tau_{J}.$ We have $a(\epsilon_0)=\dim ({\sl
F}_{J}).$

 We have introduced the approximately orthonormal basis of
 monomials:
$$\varphi_\alpha(x)=c_\alpha
h^{-n/2}(h^{-1/2} x)^\alpha,\,\,c_\alpha=(\pi^n\alpha !)^{-1/2}.$$
Then
$$(\alpha !)^{-1}((h^{1/2}\partial_x)^\alpha
u)(0)(h^{-1/2}x)^\alpha=\varphi_\alpha\cdot (\pi
h)^{n/2}(\alpha !)^{-1/2}((h^{1/2}\partial_x)^\alpha u)(0).$$

Let
 ${\displaystyle
 R_-:\,\,\cz^a\mapsto
F_{J},\,\, R_+:\,\,F_{J} \mapsto\cz^a}$ be
defined by the formulas
\begin{align*}
 &u\in F_{J},\,\,
 u_-=\{u_{-,\alpha}\}_{|\alpha\in J}\in\cz^a,\,\,
 R_- u_-=\sum_{\alpha\in
J} u_{-,\alpha} \varphi_\alpha ,\\
&R_+u=\left\{  (\pi h)^{n/2}(\alpha
!)^{-1/2}(h^{1/2}\partial_x)^\alpha u(0)\right\}_{\alpha\in
J}.
\end{align*}
By direct calculation we get $R_+u=\langle
u,\varphi_\alpha\rangle_{\alpha\in J}(1+{\mathcal
O}(h^{N})).$

Let $$
\|u\|_\Omega:=\|u\|_{H_\Phi(\Omega)},\,\,|x|_{\cz^a}:=\sqrt{\sum_{\alpha\in
J}|x_\alpha|^2}.$$ We write $\|.\|$ for natural operator
norm. We have
\begin{align*}
&R_+R_-=Id_{|\cz^a},\,\,R_-R_+=\tau_{J}=Id_{|F_{J}},\,\,
\tau_N R_-=R_-,\,\,R_+\tau_N=R_+,\\
& \| R_-\|=1+{\mathcal O}(h^{N}),\,\,\| R_+\|=1+{\mathcal
O}(h^{N}).\end{align*}
The estimates on the norms of operators $R_\pm$ follow
 from  Lemma \ref{l-monomials}.

 We pose the first
Grushin problem for $\lambda\in\Lambda_{A,B}(\epsilon_0)$:
\begin{equation}\label{GrPr1}
 \left\{\begin{array}{ll}
 \tau_N( I-e^{i2d\lambda}e^{-i\lambda F^r(I;1/\lambda)})\tau_N u +R_-u_-=v, \\
  R_+u=v_+,\,\,u,v\in{\sl F}^N.
\end{array}\right.
\end{equation}
\begin{theorem}\label{thGrPr1}
 For a given $(v,v_+) \in {\sl F}^N\times
\cz^a $ and all $\lambda\in \Lambda_{A,B}(\epsilon_0),$  the Grushin
problem (\ref{GrPr1}) has unique exact solution $(u,u_-)\in {\sl
F}^N\times\cz^a,$ which satisfies
\begin{equation}\label{GrPr1estimate}
\| u\|_{H_{\Phi}(\Omega_1)} +|u_-|_{\cz^a}\leq {\mathcal O}\left(
(\Re\lambda)^{\delta}\right) (\|v\|_{H_{\Phi}(\Omega_2)}
+|v_+|_{\cz^a}), \,\,0<\delta <1,
\end{equation}
where $\delta$ is as in (\ref{delta}), $\Omega_2\subset\subset\Omega_1\subset\subset V_{0}$ and
$|x|_{\cz^a}:=\sqrt{\sum_{\alpha\in J}|x_\alpha|^2}.$
\end{theorem}
\underline{\bf Proof:}
 We denote  $h=1/\Re\lambda,$
$\lambda\in\Lambda_{A,B}(\epsilon_0),$
$z=e^{i2d\lambda}.$
 The solution operator for (\ref{GrPr1}) is
$${\mathcal E}^0=\begin{pmatrix}
  E^0 & E_+^0 \\
  E_-^0 & E_{+-}^0\end{pmatrix},$$ where
\begin{align*}
  E^0(v)&=\sum_{\alpha\not\in
  J, |\alpha| <N}(1-zK_\alpha)^{-1} (\alpha !)^{-1}
((h^{1/2}\partial_x)^\alpha v)(0)(h^{-1/2}x)^\alpha\\
&=\sum_{\alpha\not\in
  J, |\alpha| <N}(1-zK_\alpha)^{-1}\langle
  v,\varphi_\alpha\rangle\varphi_\alpha\left( 1 +{\mathcal O}(h^{N})\right),\\
 E_+^0(v_+)& =\sum_{\alpha\in J} v_{+,\alpha} (\pi^n\alpha !)^{-1/2}h^{-n/2}(h^{-1/2}x)^\alpha =\sum_{\alpha\in J}
v_{+,\alpha} \varphi_\alpha(x)=R_-v_+,\\
  E_-^0(v) &= \left((\pi
h)^{n/2}(\alpha !)^{-1/2} (h^{1/2} \partial_x)^\alpha
v(0)\right)_{\alpha\in J}=\langle
v,\varphi_\alpha\rangle_{|\alpha\in J}\left( 1
+{\mathcal O}(h^{N})\right)=R_+v,\\
  E_{+-}^0 (v_+)&=-\{(1-zK_{\alpha}(\lambda))\}_{|\alpha\in J}v_+,
\end{align*}
where $\{(1-zK_{\alpha}(\lambda))\}_{|\alpha\in J}$ is diagonal
matrix with entrees $(1-zK_{\alpha_1}),\ldots,(1-zK_{\alpha_{N}}).$
Some of entrees can coincide.
  Solutions to (\ref{GrPr1}) are
$u=E^0(v)+E_+^0(v_+),\,\, u_-=E_-^0(v)+E_{+-}^0(v_+).$

Let
$B_{\epsilon\sqrt{s}}\subset\subset\Omega_2\subset\subset\Omega_1\subset\subset
V_{0},$ $s=h\ln(1/h).$ We have the following estimates:
\begin{align} &\| E^0(v)\|_{\Omega_1}^2\leq
\epsilon_0^{-2}\sum_{\alpha\not\in J,|\alpha|
<N}(\pi h)^{n}(\alpha !)^{-1} \left|((h^{1/2}\partial_x)^\alpha
v)(0)\right|^2(1+ {\mathcal O}(h^{N}))\leq\label{estimateE}\\
& \leq\epsilon_0^{-2}\|v\|_{\Omega_2}^2(1+
{\mathcal O}(h^{N})),\,\,
\|E_+^0(v_+)\|_{\Omega_1}^2=\sum_{\alpha\in
J}|v_{+,\alpha}|^2(1+{\mathcal O}(h^{N}))\nonumber,
\end{align}
as $E_+^0(v_+)=R_-v_+.$

We get
 \begin{equation}\label{eqq1}
 \| u\|_{\Omega_1}=\|E^0(v)+E_+^0(v_+)\|_{\Omega_1}
 \leq \epsilon_0^{-1}\left( 1+{\mathcal O}(h^{N})\right)(\| v\|_{\Omega_2}+|v_+|_{\cz^a}).
 \end{equation}
As $E_-^0(v)=R_+v,$ we have already proven that
$$|E_-^0(v)|_{\cz^a}^2\leq \left( 1+{\mathcal O}(h^{N})\right)\|v\|_{\Omega_2}^2.$$
By  (\ref{bound2}), (\ref{bound1}) and (\ref{epsilon0}) we have also
$$|E_{+-}^0(v_+)|_{\cz^a}\leq Ca\epsilon_0
|v_+|_{\cz^a}\leq {\mathcal
O}(1)\ln\Re\lambda\cdot(\Re\lambda)^{-\delta}|v_+|_{\cz^a},\,\,0<\delta
<1.$$
It implies
\begin{equation}\label{eqq2}|u_-|_{\cz^a} =|E_-^0(v)|_{\cz^a} +
|E_{+-}^0(v_+)|_{\cz^a}\leq (1+{\mathcal O}(h^{N}))
(\|v\|_{\Omega_2} +|v_+|_{\cz^a}).\end{equation}
   Adding inequalities (\ref{eqq1}) and (\ref{eqq2}) and using (\ref{epsilon0}) we get
the estimate.\hfill\qed\\

As the Grushin problem (\ref{GrPr1}) is invertible then it is well
posed. Let
$$P_N(\lambda):=\tau_N(I-e^{i2d\lambda}e^{-i\lambda F^r(I;1/\lambda)}(\lambda))\tau_N.$$ By
general consideration it follows that
$E_+^0:{\mathcal Ker}\,E_{-+}^0\mapsto{\mathcal
Ker}\,P_N(\lambda)$ is a bijection. As $P_N(\lambda)$ is a family
of Fredholm operators depending holomorphically on $\lambda
\in\Lambda_{A,B}(\epsilon_0)$ then
$$\Lambda_{A,B}(\epsilon_0)\ni\lambda\,\,\mapsto\,\,
\left(P_N(\lambda)\right)^{-1}$$ is a meromorphic family of
operators.

We have the formula
$$\left(P_N(\lambda)\right)^{-1}=E^0(\lambda)
-E_+^0(\lambda)\left( E_{-+}^0(\lambda)\right)^{-1} E_-(\lambda).$$

 Let $\gamma$ be simple loop such that $E_{-+}^0(\lambda)$
is invertible for $\lambda\in\gamma.$

The number of poles $n(\gamma)$ of $P_N(\lambda)^{-1}$ inside
$\gamma$ counted with there multiplicity is given by
$$n(\gamma)=\frac{1}{2\pi i}\tr\int_{\gamma}\left(P_N(\lambda)\right)^{-1}d\lambda=\frac{1}{2\pi i}\int_{\gamma}\tr
\left(E_+^0\left( E_{-+}^0\right)^{-1}
E_-\right)d\lambda.$$ From the other side, the number of
roots of $\det E_{-+}^0$ inside $\gamma$ is equal to
$$m(\gamma)=\frac{1}{2\pi i}\int_\gamma\frac{\left( \det E_{-+}^0(\lambda)\right)'}{\det
E_{-+}^0}d\lambda=\frac{1}{2\pi i}\int_\gamma\tr\left\{\left(
E_{-+}^0\right)'\left( E_{-+}^0\right)^{-1}\right\} d\lambda .$$

 As
in \cite{GerardSjostrand1987} we have then $n(\gamma)=m(\gamma).$

\section{Some useful estimates.}

\subsection{Estimate on $u_N\in{\mathcal Ker}\,\tau_N$}

Let $s=h\ln (1/h),$ $h=1/\lambda_1,$ $\lambda_1=\Re\lambda,$
$\lambda\in\Lambda_{A,B},$ $B_r=B(0,r)=\{x,\,\,|x|\leq r\}.$

\begin{lemma}\label{l-0}
Let $N\in\nz,$ $N={\mathcal O}(\ln\Re\lambda),$ satisfying
(\ref{N}). Let $u_N\in H_\Phi(B_{\epsilon\sqrt{s}})\cap{\mathcal
Ker}\,\tau_N.$ Put
$$r_m:=\frac{\sqrt{N+1/2}}{2\sqrt{s\lambda_1}}={\mathcal O}(1).$$
Let $\epsilon\geq r_m$ be
such that $\Phi=|x|^2/2$ for $x\in B_{\epsilon\sqrt{s}}.$ Then for any $r_0\leq r_m,$
$$\|u_N\|_{B_{r_0\sqrt{s}}}\leq\|u_N\|_{B_{\epsilon\sqrt{s}}\setminus
B_{r_0\sqrt{s}}}.$$
\end{lemma}

\underline{\bf Proof:} We need to prove $$\int_{B_{r_0\sqrt{s}}}
e^{-|x|^2/h}|u_N(x)|^2 L(dx)\leq \int_{\complement
B_{r_0\sqrt{s}}\cap B_{\epsilon\sqrt{s}}} e^{-|x|^2 /h}|u_N (x)|^2
L(dx).$$  We change the variables $x=\tilde{x} \sqrt{s}$ in
order to have the domain of integration independent of $h:$
$$x\in
B_{\epsilon\sqrt{s}}\,\,\Leftrightarrow\,\,\tilde{x}\in
B_{\epsilon}.$$

 Denote
$\Phi_s=\frac{1}{s}\Phi(\sqrt{s} x)$ and  $$ H_{\Phi_s}(B_\epsilon)={\mathcal Hol}\,(B_\epsilon)\cap L_{\Phi_s}
    (B_\epsilon)=
  \{ v\in{\mathcal Hol}\,(B_\epsilon),\,\,\int_{B_\epsilon} |v|^2
e^{-s2\Phi_s/h}L(dx) <\infty\}.$$

Let $u\in H_{\Phi}(B_\epsilon),$ ${\displaystyle
u_N(x)=u(x)-\sum_{|\alpha| <N}\frac{\partial^\alpha_xu(0)}{\alpha !}
x^\alpha},$  $B_\epsilon :=B(0,\epsilon ).$

Suppose that $\epsilon >0$ is small enough such that $\Phi =|x|^2/2$
on the domain of integration. Thus $\Phi_s=|x|^2/2.$ We must prove
that for ${\displaystyle (s\lambda_1)^{-1/2}\sqrt{N+1/2}\leq
\epsilon }$ and ${\displaystyle
r_m=1/2(s\lambda_1)^{-1/2}\sqrt{N+1/2}}$ we have
$$\int_{B_{r_m} } e^{-s\lambda_1 |x|^2}|u_N(x)|^2L(dx)\leq
\int_{\complement B_{r_m}\cap B_\epsilon} e^{-s\lambda_1
|x|^2}|u_N(x)|^2L(dx).$$

  The monomials
constitute an orthogonal base in $H_\Phi(B_\epsilon)$ and it is
enough to show the inequality for $x^\alpha,$ $|\alpha |\geq N.$ Let
$\alpha_j\in\nz$ be such that $\alpha_j\geq N.$ For $u_N(x)=
x_j^{\alpha_j}$ we have $$\int_{|x|\leq\epsilon} e^{-s\lambda_1
|x|^2} |f_N(x) |^2L(dx)=C_n(s\lambda_1)^{-n/2}\int_{r\leq \epsilon
\sqrt{s\lambda_1}} e^{-r^2} r^{2\alpha_j+1} dr.$$ As $$\left(
e^{-r^2} r^{2\alpha_j+1}\right)'_r=(-2r^2+2\alpha_j +1)
e^{-r^2}r^{2\alpha_j}=0\,\,\Leftrightarrow\,\,
r=r_m(\alpha):=\sqrt{\alpha_j +1/2}$$ we get that the function
${\displaystyle e^{-r^2}r^{2\alpha+1}}$ is increasing for $0\leq
r\leq r_m.$ Then
$$\int_0^{r_m/2} e^{-r^2} r^{2\alpha_j +1} dr\leq
\int_0^{r_m/2}e^{-(r+r_m/2)^2}\left(
r+\frac{r_m}{2}\right)^{2\alpha_j +1} dr=\int_{r_m/2}^{r_m}e^{-r^2}
r^{2\alpha_j +1} dr.$$ For all $\alpha$ such that $|\alpha |>N$ the
maximum $r_m(\alpha)
>r_m(\alpha_N)$ with $|\alpha_N|=N.$ Thus,  if the
above estimate is valid  for $|\alpha_N|=N$  then it implies that it
is valid for any $|\alpha| >N.$ \hfill$\blacksquare$

\subsection{Model estimate}
\begin{lemma}\label{l-preparatory} As before we write $s=h\ln (1/h),$  $h=1/\lambda_1,$ $\lambda_1 =\Re\lambda.$
Let $\Omega$ be any $\lambda$-independent \neigh{} of $0$ and
 $r_0>0$   such that
$B_{r_0\sqrt{s}}\subset\subset\Omega.$

Let $A={\rm diag}(e^{\mu_1},\ldots, e^{\mu_n})$ be as in Theorem
\ref{Th.2.6}.  Denote $A^{-1}(\Omega):=\{ A^{-1}x,\,\,x\in\Omega\}.$

For any $c_0 >0$ and $|z|\leq e^{c_0\ln\lambda_1}$ and $r_0$  large
enough there exists $m\leq 1/2$ such that
 \begin{equation*}\int_{ \Omega\setminus
B_{r_0\sqrt{s}}} e^{-2\Phi(x) /h} |z
e^{-\sum_{j=1}^n\mu_j(x_j\partial_j +1/2)} u |^2 L(dx)\leq
m^2\int_{A^{-1}(\Omega)\setminus A^{-1}(B_{r_0\sqrt{s}})}
e^{-2\Phi(x) /h} |u|^2 L(dx).\end{equation*}
\end{lemma}
\underline{\bf Proof:} We have $$e^{-\sum_{j=1}^n\mu_j(x_j\partial_j
+1/2)} u=e^{-\sum_{j=1}^n\mu_j/2} u(A^{-1} x),\,\,A={\rm
diag}(e^{\mu_1}\ldots e^{\mu_n}).$$ Let first $\Phi=|x|^2/2.$  We
put for simplicity $n=1,$ $\mu_j=\mu.$ Then $ e^{-\mu(x\partial
+1/2)} u =e^{-\mu/2} u(e^{-\mu}x).$ The general case is
straightforward.  We change the variables $
\tilde{x}=e^{-\mu}x,\,\,L(dx)=e^{2\mu}L(d\tilde{x}),$ and denote the
new variable $\tilde{x}$ again by $x.$ Let $ A^{-1}(\Omega):=\{
e^{-\mu}x,\,\,x\in\Omega\}.$ Then we have, with $s=h\ln\lambda_1,$
\begin{align}
&e^{-\mu}\cdot\int_{\Omega\setminus B_{r_0\sqrt{s}}} e^{-|x|^2/h} |z
u (e^{-\mu}x)|^2 L(dx)  \leq
e^{2c_0\ln\lambda_1+\mu}\int_{A^{-1}(\Omega)\setminus
B_{r_0\sqrt{s}e^{-\mu}}} e^{-
 e^{2\mu}|x|^2/h}|u(x)|^2L(dx) =\nonumber\\
 &=
e^{2c_0\ln\lambda_1+\mu}\int_{A^{-1}(\Omega)\setminus
B_{r_0\sqrt{s}e^{-\mu}}} e^{- (e^{2\mu}-1)|x|^2/h}
e^{-|x|^2/h}|u(x)|^2 L(dx) \leq\nonumber\\&\leq
 e^{2c_0\ln\lambda_1+\mu}e^{-
(e^{2\mu}-1)e^{-2\mu}r_0^2s/h}\int_{A^{-1}(\Omega)\setminus
B_{r_0\sqrt{s}e^{-\mu}}} e^{- |x|^2/h}|u(x)|^2
L(dx)=\label{here}\\
&= m^2\int_{A^{-1}(\Omega)\setminus B_{r_0 e^{-\mu}\sqrt{s}}} e^{-
|x|^2/h}|u(x)|^2 L(dx)\nonumber
\end{align}
with
$$m^2:=  e^{2c_0\ln\lambda_1+\mu- \ln\lambda_1
(1-e^{-2\mu})r_0^2}.$$

 We have $m\leq 1/2$ if
\begin{align*}
& e^{2c_0\ln\lambda_1+\mu- \ln\lambda_1 (1-e^{-2\mu})r_0^2}\leq
2^{-2}\,\,\Leftrightarrow\,\, 2c_0\ln\lambda_1+\mu- \ln\lambda_1
(1-e^{-2\mu})r_0^2\leq -2\ln
2\,\,\Leftrightarrow\\
&r_0\geq\sqrt{\frac{2c_0\ln\lambda_1+\mu +2\ln 2}{\ln\lambda_1\cdot
(1-e^{-2\mu})}}.\end{align*}

For general weights  $\Phi$ we use the strict convexity and get
(after change of variables $\tilde{x} =e^{-\mu} x$)
$$\Phi(e^\mu\tilde{x}) -\Phi (\tilde{x})\geq
c|\tilde{x}|^2,\,\,c>0,\,\, x\in V_1,$$ and in the above formula
(\ref{here}) exchange $e^{2\mu}-1$ by $c.$ \hfill\qed

\subsection{General estimate} In this section we will write
$h=1/\lambda$ and suppose that $\lambda$ is real. Generalization to
 $\lambda\in\Lambda_{A,B}$ is straightforward.

  From (\ref{expan}) it follows that
$$e^{-i H^r(hy)/h}\sim 1+h p_1(y)
+h^{2}p_2(y)+\ldots,\,\, |y| \leq |h|^{\delta -1/2},\,\,\delta
>0,$$ where the polynomial $p_j(y)=\sum_{|\alpha|\in [j+1 ,2j]}
q_{\alpha,j} y^\alpha$ is a linear combination of monomials of
degree in $[j+1, 2j].$ Moreover, $p_1(y)=-\frac{i}{2} H''(0)y\cdot
y.$ Then
\begin{align*}e^{-i F(h y)/h}&=  e^{-i\mu y}\left(
1+h p_1 (y) +h^{2} p_2(y)+\ldots\right) =\left( 1+h
p_1(i\partial_\mu) +h^{2}p_2(i\partial_\mu) +\ldots\right) e^{-i\mu
y}.
\end{align*}
Taking $y=h^{-1}\cdot I=\left(x\partial_x +1/2\right) /i,$ we have
the following representation:
\begin{align*}
& e^{-i F_0^r(I)/h}= e^{-i  H^r(i h\partial_\mu)/h}
e^{-\sum_{i=1}^n\mu_i(x_i\partial_i+1/2)},
\,\,F_0(\imath)=G(\imath,\mu)+H^r(\imath)=\sum_{i=1}^n
   \mu_i\cdot\imath_i +{\mathcal O}_r(\imath^2),
   \end{align*} polynomial of degree $r$
   independent of $h.$

 Then,
\begin{align*}
&  e^{-i  H^r(ih\partial_\mu)/h} e^{-\sum_{i=1}^n\mu_i
(x_i\partial_i +1/2)} u=  e^{-i H(i h\partial_\mu)/h}
e^{-\sum_{j=1}^n\mu_j/2} u(A^{-1} x),\,\, A={\rm
diag}(e^{\mu_1}\ldots e^{\mu_n}).
\end{align*}

For simplicity we put $n=1.$ We have
\begin{align*}
e^{-iF_0^r(I)/h}u&= e^{-i H(i h\partial_\mu )/h}e^{-\mu
(x\partial+1/2)}
u(x)=\\
&=e^{-\mu /2} u(e^{-\mu}x) + h\sum_{j=1}^{r-1}h^{j-1} \sum_{l\in
[j+1, 2j]} q_{l,j}\partial_\mu^l(e^{-\mu /2} u(e^{-\mu}x))=\\
& =e^{-\mu /2} u(e^{-\mu}x) +e^{-\mu /2}T_0^r(e^{-\mu}x),\\
\mbox{where}\,\,&T_0^ru(x):= h\sum_{j=1}^{r-1} h^{j-1}\sum_{l\in
[j+1, 2j]} q_{l,j}\sum_{k=0}^l c_{k,l} u^{(k)}(x)\cdot x^k={\mathcal
O}(h).
\end{align*}

In general case
$$e^{-iF(hy; h)/h} = e^{-i  K(0;h)/h} e^{-i G(y,\mu)} (1+h q_1 (y)
+h^{2} q_2 (y) +\ldots),$$ where $q_j(y)$ is a polynomial of degree
at most $2j,$ and as before, this can be written
$$e^{-i F(h y;h)/h}= e^{-i K(0;h )/h} (1 +h q_1 (i\partial_\mu )
+h^{2} q_2 (i\partial_\mu)+\ldots) e^{-i G(y,\mu)}.$$ We will write
it in the form
$$\frac{F(hy;h )}{h} =\frac{F(0;h )}{h} + G(y,\mu)
+\frac{J(hy;h )}{h}$$ and then (see also
\cite{IantchenkoSjostrandZworski2002}, page 347)
$$e^{-iF (I(h); h)/h} =e^{-i F(0;h)/h} e^{-i  J( i h\partial_\mu
;h )/h } e^{-\sum_{i=1}^n\mu_i ( x_i\partial_i +1/2)}$$ and the
difference from the case $F=G(I,\mu)+H(I)$ is only the factor $e^{-i
F(0;h)/h},$ which is bounded as $h\rightarrow 0,$ and that $q_j (y)$
is a polynomial of degree at most $2j$ and not just in $[j+1, 2j].$

We have
\begin{align}&  e^{-i F (I(h); h)/h}u
 =e^{-i F(0;h)/h} e^{-i J( i h\partial_\mu ;h) /h}
e^{-\sum_{i=1}^n\mu_i ( x_i\partial_i +1/2)}u=\nonumber\\
&e^{-i F(0;h)/h} \left(e^{-\mu /2} u(e^{-\mu}x)
+h\sum_{j=1}^{r-1}h^{j-1} \sum_{l=0}^{2j}
 q_{l,j}\partial_\mu^l(e^{-\mu /2} u(e^{-\mu}x))\right)=\nonumber\\
& =e^{-i F(0;h)/h} \left(e^{-\mu /2} u(e^{-\mu}x) +e^{-\mu /2}T^r(e^{-\mu}x)\right),\label{formula1}\\
&\mbox{where}\,\,T^ru(x):= h\sum_{j=1}^{r-1} h^{j-1}\sum_{l=0}^{2j}
q_{l,j}\sum_{k=0}^l c_{k,l} u^{(k)}(x)\cdot x^k.\nonumber
\end{align}

 We put $z=e^{-i2d\lambda}.$ Then for $\lambda\in\Lambda_{A,B}$ we
  have $|z|<e^{2dA\ln\Re\lambda},$
 with $A$ as in the definition of $\Lambda_{A,B}.$

We can apply Lemma \ref{l-preparatory} with $c_0=2dA.$ Then, for
$r_0$ large enough we have
\begin{align}
&e^{-\mu}\cdot\int_{\Omega\setminus B_{r_0\sqrt{s}}}e^{-2\Phi(x)/h}
|z |^2 \left| u(e^{-\mu} x)+T^ru(e^{-\mu}
x)\right|^2L(dx))\leq \nonumber\\
&m^2\int_{A^{-1}(\Omega)\setminus
A^{-1}(B_{r_0\sqrt{s}})}e^{-2\Phi(x)/h} \left|
u(x)+T^ru(x)\right|^2L(dx).\label{23}\end{align}

For any $\Omega_2\subset\subset\Omega_1\subset\subset\Omega,$ we
have for any $u\in H_\Phi(\Omega),$
$$\| (h^{-1/2} x)^\alpha (h^{1/2} D_x)^\beta u\|_{\Omega_2}\leq \|
\left( 1+h^{-1/2}|x|\right)^{|\alpha|+|\beta |} u\|_{\Omega_1},$$
see Corollary 4.2 in \cite{GerardSjostrand1987}.

If $x\in\Omega\subset V_0$ then $|x|\leq {\mathcal O}(1) h\ln
(1/h).$ Taking $u\in {\mathcal Ker}(\tau_N)$ and  using Lemma
\ref{l-0} we get
$$\| (h^{-1/2} x)^\alpha (h^{1/2} D_x)^\beta u\|_{\Omega_2}\leq
{\mathcal O}\left( \left( \ln(1/h)\right)^{|\alpha|+|\beta |}\right)
\|u\|_{\Omega_1\setminus B_{r_0\sqrt{s}}}. $$  Thus
 the
partial differential operator $T^r$ is bounded of norm ${\mathcal
O}(h):$ $H_\Phi(\Omega_1)\mapsto H_\Phi(\Omega_2)$ and as $u\in
{\mathcal Ker}(\tau_N)$ we have,
\begin{equation}\label{eq-coco}
\| T^r u\|_{\Omega_2}\leq{\mathcal
O}(h)\|u\|_{\Omega_1}\leq{\mathcal O}(h)\| u\|_{\Omega_1\setminus
B_{r_0\sqrt{s}}}.
\end{equation}

Then, using (\ref{eq-coco}) with $\Omega_1=\Omega$ and
$\Omega_2=A^{-1}(\Omega),$  we get
$$ \int_{A^{-1}(\Omega)\setminus A^{-1}(B_{r_0\sqrt{s}})}e^{-2\Phi(x)/h}  \left|
T^ru(x)\right|^2L(dx)\leq {\mathcal O(h^2)} \int_{\Omega\setminus
B_{r_0\sqrt{s}}}e^{-2\Phi(x)/h} \left| u(x)\right|^2L(dx).$$

We can  choose $r_0$ large enough such that
$$\mbox{the \rhs{} of\,\,}(\ref{23})\leq\frac14 \int_{\Omega\setminus
B_{r_0e^{-\mu}\sqrt{s}}}e^{-2\Phi(x)/h} | u(x)|^2L(dx).$$

This implies the following lemma:
\begin{lemma}\label{l-main} Suppose
$\lambda\in\Lambda_{A,B}$ and $|\lambda|>C$ large enough. Let
$h=1/\Re\lambda,$ $r_0,$ $\epsilon$ be as in Lemma \ref{l-0}. Let
$\Omega$ be any $\lambda$-independent \neigh{} of $0$ and
$$B_{r_0\sqrt{s}}\subset\subset
B_{\epsilon\sqrt{s}}\subset\subset\Omega.$$

Let $A={\rm diag}(e^{\mu_1},\ldots, e^{\mu_n})$ be as in Theorem
\ref{Th.2.6}.  Let $A^{-1}(\Omega):=\{ A^{-1}x,\,\,x\in\Omega\}.$

  If $r_0$ and $\epsilon$   are sufficiently large then for some
$m\leq 1/2$ \begin{equation}\label{wanted}\int_{ \Omega\setminus
B_{r_0\sqrt{s}}} e^{-2\Phi /h} |e^{-i2d\lambda} e^{-i\lambda
F^r(I;1/\lambda)} u |^2 L(dx)\leq m^2\int_{\Omega\setminus
A^{-1}(B_{r_0\sqrt{s}})} e^{-2\Phi /h} |u|^2 L(dx)\end{equation} for
any $u\in H_\Phi(\Omega)\cap{\mathcal Ker}(\tau_N).$
\end{lemma}

\subsection{Bounds useful for the second local Grushin problem, Theorem \ref{thGrPr2}}
We suppose $\lambda\in\Lambda_{A,B}.$  Vi fix $N\in\nz$ and let $u_N\in
H_\Phi(V_{0})\cap{\mathcal Ker}(\tau_N).$ Let $r_0>0$ be large enough.
In all lemmas below we use the following convention for the
domains:
\begin{equation}\label{domm}
B_{r_0\sqrt{s}}\subset\subset
B_{\epsilon_0\sqrt{s}}\subset\subset\Omega_2\subset\subset\Omega_1\subset\subset
V_{0}.\end{equation}

Lemma \ref{l-R} implies that the norm of $R_r=M_0-e^{-i\lambda
F^r(I;1/\lambda)}$ is small in $H_\Phi(V_{0}).$ We have
\begin{lemma}\label{co-1}
For any $\tilde{N}\in\nz$ there is $r=r(N,\tilde{N}, r_0)\in\nz$    sufficiently large
such that
$$\|e^{-i2d\lambda} R^r u_N \|_{\Omega_2\setminus
B_{r\sqrt{s}}}^2 \leq  {\mathcal O}_r(h^{\tilde{N}})\|
u_N\|_{\Omega_1\setminus B_{r_0\sqrt{s}}}^2$$ with domains
satisfying (\ref{domm}).
\end{lemma}
\underline{\bf Proof:} Let $z=e^{-i2d\lambda}.$
 We use Lemma \ref{l-R} (with $\tilde{N}$ instead of $N$) and Lemma \ref{l-0}:
\begin{align*}&\| zR^ru_N\|^2_{\Omega_2\setminus B_{r_0\sqrt{s}}}\leq
{\mathcal O} (h^{\tilde{N}})\|u_N\|_{\Omega_1}^2={\mathcal O}
(h^{\tilde{N}})\left(\|u_N\|_{B_{r_0\sqrt{s}}}^2
+\|u_N\|_{\Omega_1\setminus B_{r_0\sqrt{s}}}^2\right)\leq\\
&\leq {\mathcal O}
(h^{\tilde{N}})\left(\|u_N\|^2_{B_{\epsilon\sqrt{s}}\setminus
B_{r_0\sqrt{s}}} +\|u_N\|^2_{\Omega_1\setminus
B_{r_0\sqrt{s}}}\right)\leq{\mathcal O}(h^{\tilde{N}})\|
u_N\|_{\Omega_1\setminus B_{r_0\sqrt{s}}}^2.
\end{align*}
\hfill\qed

 Combining Lemma \ref{l-main} and Lemma \ref{co-1} we get
\begin{lemma} Let  $M_0=e^{-i\lambda
F^r(I;1/\lambda)}+R^r.$  Suppose $\lambda\in\Lambda_{A,B},$   (\ref{domm}) and $A^{-1}(\Omega)=\{ A^{-1}x,\,\,x\in\Omega\},$ $A={\rm
diag}(e^{\mu_1},\ldots, e^{\mu_n}).$

For any $\tilde{N}\in\nz$ there is $r\in\nz$ such that
if $r_0$   is sufficiently large then for some $m\leq 1/2$
\begin{equation}\label{p-2}
\|e^{-i2d\lambda} M_0 u_N \|^2_{\Omega_2\setminus
B_{r_0\sqrt{s}}}\leq m^2\|u_N\|^2_{\Omega_2\setminus
A^{-1}(B_{r_0\sqrt{s}})}+ {\mathcal O}_r(h^{\tilde{N}})\|
u_N\|^2_{\Omega_1\setminus B_{r_0\sqrt{s}}}\end{equation} and
\begin{equation}\label{cor-p-main}\| u_N\|^2_{\Omega_2\setminus B_{r_0\sqrt{s}}}\leq {\mathcal O}(1)\|
 (I-e^{-i2d\lambda}M_0)u_N\|^2_{\Omega_2\setminus B_{r_0\sqrt{s}}}
 +{\mathcal O}_r(h^{\tilde{N}})\|u_N\|^2_{\Omega_1\setminus
 B_{r_0\sqrt{s}}}.\end{equation}
 \end{lemma}
\underline{\bf Proof:} We show (\ref{cor-p-main}). Choose $\epsilon$ such that $\Phi =|x|^2/2$
for $x\in B_{\epsilon\sqrt{s}}$ and $r_0\leq
\frac{\sqrt{N+1/2}}{2\sqrt{s\lambda_1}}$ as before. Then, with
$z=e^{-i2d\lambda},$ and some constant $m \leq 1/2,$  bound
\ref{p-2} implies
\begin{align}
 &\| (I-zM_0) u_N\|^2_{\Omega_2\setminus B_{r_0\sqrt{s}}}\geq
\| u_N\|^2_{\Omega_2\setminus B_{r_0\sqrt{s}}} -\| zM_0
u_N\|^2_{\Omega_2\setminus
B_{r_0\sqrt{s}}}\geq\nonumber\\
&\geq\| u_N\|^2_{\Omega_2\setminus B_{r_0\sqrt{s}}}
-m^2\|u_N\|_{\Omega_2\setminus A^{-1}(B_{r_0\sqrt{s}})}^2- {\mathcal
O}_r(h^{\tilde{N}})\| u_N\|_{\Omega_1\setminus B_{r_0\sqrt{s}}}^2=\nonumber\\
&=(1-m^2)\|u_N\|^2_{\Omega_2\setminus B_{r_0\sqrt{s}}}
-m^2\|u_N\|^2_{B_{r_0\sqrt{s}}\setminus A^{-1}(B_{r_0\sqrt{s}})} -
{\mathcal O}_r(h^{\tilde{N}})\| u_N\|_{\Omega_1\setminus
B_{r_0\sqrt{s}}}^2.\label{heer2}
\end{align}
Lemma \ref{l-0} implies
$$\|u_N\|^2_{B_{r_0\sqrt{s}}\setminus
A^{-1}(B_{r_0\sqrt{s}})}\leq
\|u_N\|^2_{B_{r_0\sqrt{s}}}\leq\|u_N\|^2_{B_{\epsilon\sqrt{s}}\setminus
B_{r_0\sqrt{s}}}\leq \|u_N\|^2_{\Omega_2\setminus
B_{r_0\sqrt{s}}}.$$

Then we get that
$$\mbox{the \rhs{} of (\ref{heer2})}\geq(1-2m^2)\|u_N\|^2_{\Omega_2\setminus
B_{r_0\sqrt{s}}}  - {\mathcal O}_r(h^{\tilde{N}})\|
u_N\|_{\Omega_1\setminus B_{r_0\sqrt{s}}}^2.$$

 Then we  have
  \begin{align*} (1-2m^2)\|u_N\|^2_{\Omega_2\setminus
B_{r_0\sqrt{s}}} \leq & \| (I-zM_0) u_N\|^2_{\Omega_2\setminus
B_{r_0\sqrt{s}}} +{\mathcal O}_r(h^{\tilde{N}})\|
u_N\|^2_{\Omega_1\setminus B_{r_0\sqrt{s}}}.\end{align*} \hfill\qed

Lemma \ref{l-0} implies that $\|u_N\|^2_{B_{r_0\sqrt{s}}}\leq\|
u_N\|^2_{\Omega_2\setminus B_{r_0\sqrt{s}}}$ which together with estimate (\ref{cor-p-main}) implies
 \begin{lemma}\label{locel}
Let $\lambda\in\Lambda_{A,B}.$ For any $\tilde{N}\in\nz$ there is $r\in\nz$ such that
if $r_0$   is sufficiently large then
\begin{equation}\label{VIP2}\| u_N\|_{\Omega_2}^2\leq
{\mathcal O}(1)\|(I-e^{-i2d\lambda}M_0)u_N\|_{\Omega_2}^2+ {\mathcal
O}_r(h^{\tilde{N}})\|u_N\|_{\Omega_1\setminus
B_{r_0\sqrt{s}}}^2\end{equation} with domains satisfying
(\ref{domm}).
 \end{lemma}

\section{The second local Grushin problem.}
With the same notations as in (\ref{GrPr1}), $M_0=e^{-i\lambda
F^r(I;1/\lambda)}+R^r,$ $\lambda\in\Lambda_{A,B}(\epsilon_0)$
verifying (\ref{bound1}), for any $\Omega_2\subset\subset
\Omega_1\subset\subset V_{0},$
 we pose the second Grushin problem:
\begin{equation}\label{GrPr2}
 \left\{\begin{array}{ll}
 ( I-e^{-i2d\lambda}M_0) u +R_-u_-=v, \\
  R_+u=v_+,\,\,u\in H_\Phi(\Omega_1),\,\,v\in H_\Phi(\Omega_2),\,\,u_-,v_+\in\cz^a.
\end{array}\right.
\end{equation}
\begin{theorem}\label{thGrPr2}
Let $\lambda\in \Lambda_{A,B}(\epsilon_0).$ We consider the Grushin
problem (\ref{GrPr2}).

 Then for any $\tilde{N}\in\nz$ there is  $r\in\nz$ large enough such that we have
\begin{equation}\label{GrPr2estimate}
\| u\|_{H_\Phi(\Omega_1)} +|u_-|_{\cz^a}\leq {\mathcal O}\left(
(\Re\lambda)^{\delta}\right)\left( \|v\|_{H_\Phi(\Omega_2)}
+|v_+|_{\cz^a}\right)+{\mathcal O}_r \left( |\Re\lambda|^{-\tilde{N}}\right) \|
u\|_{H_\Phi (\Omega_1)},\,\,0 <\delta <1.
\end{equation}
Here $\delta$ is  as in (\ref{delta}).
\end{theorem}
\underline{\bf Proof:} Let $z=e^{-i2d\lambda},\,\,h=1/\Re\lambda,$ and $N$ sufficiently large.   Applying $(1-\tau_N)$
to the first equation in (\ref{GrPr2}) and using $(1-\tau_N)R_-=0,$ we
get
$(1-\tau_N)(I-zM_0)u=(1-\tau_N)v$ which implies
\begin{equation}\label{e-1}(I-zM_0)(1-\tau_N)u=(1-\tau_N)v+\tau_N(I-zM_0)(1-\tau_N)u-(1-\tau_N)(I-zM_0)\tau_N
u.\end{equation}
Then we apply estimate (\ref{VIP2}) (Lemma \ref{locel}). Let
$B_{r_0\sqrt{s}}\subset\subset\Omega_2\subset\subset\Omega_1\subset\subset
V_{0}.$ Then
\begin{align}
&\| (1-\tau_N)u\|_{\Omega_2}\leq
c\|(1-\tau_N)v\|_{\Omega_2}+\nonumber\\
&+{\mathcal O}(1)\|\tau_N(I-zM_0)(1-\tau_N)u-(1-\tau_N)(I-zM_0)\tau_Nu\|_{\Omega_2}+\label{vopros}\\
&+{\mathcal O}_r(h^{N_0})\|(1-\tau_N) u\|_{\Omega_1\setminus
B_{r_0\sqrt{s}}}.\nonumber
\end{align}
We need to bound (\ref{vopros}).

We have
 $$\|\tau_Nze^{-i\lambda F^r(I;1/\lambda)}(1-\tau_N)u\|_{\Omega_2}\leq {\mathcal
 O}(h^{\tilde{N}})\|u\|_{\tilde{\Omega}_2},
 $$ where
$\Omega_2\subset\subset\tilde{\Omega}_2\subset\subset\Omega_1\subset
\subset V_{0},$
 and similar, interchanging $\tau_N$ and $1-\tau_N.$
It follows from the following facts:
$$\mbox{if}\,\, u\in{\sl O}_N={\mathcal Ker}\,\tau_N,\,\,\mbox{then}\,\,
u^{(\beta)}(e^{-\mu} x)(e^{-\mu}x)^\beta\in{\sl
O}_N\,\,\forall\beta,$$ and formula (\ref{formula1}), which implies that
$\tau_Nze^{-i\lambda F^r(I;1/\lambda)}(1-\tau_N)=0$ for any fixed
order $r.$

  Then it is enough to estimate $\|\tau_N
zR(1-\tau_N)\|_{\Omega_2}$ and $\|(1-\tau_N) zR \tau_N\|_{\Omega_2}$
which can be done directly, but from more general result,
 Lemma
\ref{l-R} on $V_{0},$ we know that  for any $\tilde{N}$ we can find $r$
such that $(\ref{vopros})\leq {\mathcal
O}_r(h^{\tilde{N}})\|u\|_{\tilde{\Omega}_2}.$

We get the estimate
\begin{equation}\label{VIP3} \| (1-\tau_N)u\|_{\Omega_2}\leq
{\mathcal O}(1)\|(1-\tau_N)v\|_{\Omega_2} +{\mathcal
O}_r(h^{\tilde{N}})\|u\|_{\tilde{\Omega}_2}.
\end{equation}
Apply $\tau_N$ to (\ref{GrPr2}). Then
$$\left\{\begin{array}{rl}
  \tau_N(I-zM_0)u+R_-u_- & =\tau_Nv \\
  R_+\tau_Nu & =v_+ \\
\end{array}
\right.
$$ since $\tau_N R_-=R_-$ and $R_+=R_+\tau_N.$
Then
$$\left\{\begin{array}{rl}
  \tau_N(I-ze^{-i\lambda F^r(I;1/\lambda)})\tau_N u+R_-u_- & =\tau_Nv -\tau_N (I-zM_0)(1-\tau_N)u-\\&\tau_N z
  (\underbrace{e^{-i\lambda F^r(I;1/\lambda)}-M_0}_{R^r})
  \tau_Nu\\
  R_+\tau_Nu & =v_+. \\
\end{array}
\right.
$$
The first term in  $-\tau_N (I-zM_0)(1-\tau_N)u-\tau_N zR^r \tau_Nu$ which is equal to  $\tau_NzM_0(1-\tau_N)u=\tau_Nze^{-i\lambda
F^r(I;1/\lambda)}(1-\tau_N)u
 +\tau_NzR^r(1-\tau_N)u$ is already estimated. The second term can be bounded as before using   Lemma \ref{l-R} in Section \ref{s-qbnf}.

Then estimate (\ref{GrPr1estimate}) implies
\begin{equation}\label{VIP4}
\|\tau_N u\|_{\Omega_1} +| u_-|_{\cz^a}\leq {\mathcal O}\left(
h^{-\delta}\right)\left(
\|\tau_Nv\|_{\Omega_2}+|v_+|_{\cz^a}\right) +{\mathcal
O}_r(h^{\tilde{N}})\|u\|_{\tilde{\Omega}_2}.
\end{equation}
Adding (\ref{VIP3}) and (\ref{VIP4})
 we get
\begin{equation}\label{e-2}\| u\|_{\Omega_1}+|u_-|_{\cz^a}\leq
{\mathcal O}(1)\|(1-\tau_N)v\|_{\Omega_2} +{\mathcal O}\left(
h^{-\delta}\right)\left(
\|\tau_Nv\|_{\Omega_2}+|v_+|_{\cz^a}\right)  +{\mathcal
O}_r(h^{\tilde{N}})\| u\|_{\tilde{\Omega}_2} .\end{equation}
Then using Proposition 4.3 in
  \cite{GerardSjostrand1987},
 we get the desired estimate (\ref{GrPr2estimate}).\hfill\qed


\section{The global Grushin problem}\label{gr-global}

We denote $\Omega_{\rm int}\in T^*\partial\Omega_1$  the
$\lambda$-independent \neigh{} of $a_1$ such that $W_0\subset\subset
\Omega_{\rm int}\subset\subset W_1 ,$ where $W_0,$ $W_1$ are defined
in (\ref{W0}), (\ref{W1}).
 Let $\Omega_{\rm ext}$ be  such that $\complement\Omega_{\rm ext}\subset\subset\Omega_{\rm
int}.$

We consider the original operator $I-M(\lambda)$ in ${\mathcal H}$
 equipped with the norm
$$|u|_{\mathcal H}:=\|(1-\chi _2)T_{\Lambda_{tG}} u\|_{L^2_{tG}(T^*(\partial\Omega_1))}
+\| \chi_1 U u\|_{L^2(e^{-2\Phi/h} L(dx))}=:|u|_{\Omega_{\rm ext}}
+|u|_{\Omega_{\rm int}},\,\, u\in {\mathcal H},
$$ where $\chi_1\in C_0^\infty(\cz^n)$ is equal to $1$ in a
\neigh{} of $a_1=0,$ with the support independent of $\lambda$ and
$W_0\subset\subset\supp\chi_1\subset\subset \Omega_{\rm
int}\subset\subset W_1,$ and $\chi_2\in C_0(T^*X)$ equal to $1$ near
$(a_1,0)$ and essentially the same function as $\chi_1$ after
suitable identification of domains. Here
$L^2_{tG}(T^*(\partial\Omega_1))$ stands for $L^2(\Lambda;
e^{-2H/h}|\langle\alpha_\xi\rangle |^{2m}d\alpha).$ Let
$$U:\,\,{\mathcal H}\mapsto H_\Phi(\Omega_{\rm int}),\,\, U=B_rF
T_{\rm Bargman},\,\,V_1\supset{\rm neigh}(a_1)\mapsto{\rm
neigh}(0)\subset \hat{V}_1 ,$$
 $\pi(W_1)=V_1,$ be   as before (we omit $\hat{}$ ) such that
$$U M(\lambda) u= e^{-2id\lambda} M_0  U u +{\mathcal
O} \left( h^{N}\right)\| u\|_{\mathcal H},\,\, \mbox{in}\,\,
H_{\Phi}(\Omega_{\rm int}),\,\,M_0-e^{-i\lambda F^r}\equiv{\mathcal
O}(\rho^{2r+1},h^{2r+1}),$$ where
$h=1/\Re\lambda,\,\,\lambda\in\Lambda_{A,B}.$

  Let $V:\,\,H_\Phi(\Omega_{\rm int})\mapsto
H(\Lambda_{tG})$ be  the approximate microlocal inverse of $U$ such
that if $\chi (x,hD_x)$ is \pseudor{} adapted to $H(\Lambda_{tG})$
with compact symbol and with $\supp\chi\subset\subset\Omega_{\rm
int},$ then $\chi(VU-I)$ and $(UV-I)\chi$ are neglectible,
$VUu-u={\mathcal O}(h^{N})\|u\|_{\mathcal H}.$

\begin{theorem}\label{GrushinM1}

 For a given $(v,v_+)\in
{\mathcal H}\times \cz^a$ consider the Grushin problem
\begin{equation}\label{GrPr}
\left\{\begin{array}{rl}
  (I-M)u + V R_-u_-&=v,  \\
  R_+ U u&=v_+.
\end{array}\right.
\end{equation}

 For  $\lambda\in \Lambda_{A,B}(\epsilon_0),$ verifying (\ref{bound1}),
and $\Re\lambda$ sufficiently large the Grushin problem (\ref{GrPr})
has a unique solution $(u,u_-)\in {\mathcal H}\times\cz^a.$
Moreover, we have the a priori estimate
\begin{align}\label{finalestimate} &|u|_{{\mathcal H}}
+|u_-|_{\cz^a}\leq  {\mathcal O}((\Re\lambda)^{\delta})\left(|v|_{{\mathcal
H}}  +| v_+|_{\cz^a}\right)
\end{align}
with some $\delta,$ $0<\delta <1.$
\end{theorem}
\underline{\bf Proof:} The existence of a solution follows as in \cite{Gerard1988}. We need to show estimate (\ref{finalestimate}). We apply $U$ from the left in the first equation
in (\ref{GrPr}) and denote $Uu=\tilde{u}.$ Then
$$\left\{\begin{array}{rl}
  (I-e^{-2id\lambda}M_0)\tilde{u} +  R_-u_-&=U v +\tilde{w},\,\,\tilde{w}:=(UM-e^{-2id\lambda} M_0U) u+(I-U V)R_-u_-,  \\
  R_+ \tilde{u}&=v_+.
\end{array}\right.
$$ Denote $\tilde{v}=Uv+\tilde{w}.$
Introduce as before $\Omega_{\rm ext},$ $\Omega_{\rm int},$ $W_{0},$
$V_{0}=\pi(W_{0}).$ Let
$$\Omega_2\subset\subset\Omega_1\subset\subset V_{0}\subset\subset
\Omega_3'\subset\subset \Omega_2'\subset\subset
\Omega_1'\subset\subset\Omega_{\rm int}.$$
{\bf Estimate  on $V_{0}:$}\\
By (\ref{GrPr2estimate}) we have:
\begin{equation}\label{1int}
\| \tilde{u}\|_{H_\Phi(\Omega_1)} +|u_-|_{\cz^a}\leq {\mathcal O}
\left(h^{-\delta}\right)\left(
\|\tilde{v}\|_{H_\Phi(\Omega_2)} +|v_+|_{\cz^a}\right)+{\mathcal
O}_r \left( h^N\right)\| \tilde{u}\|_{H_\Phi (\Omega_1)}
,\,\,0<\delta <1.
\end{equation}
 We have
\begin{equation}\label{2int}
\|\tilde{w}\|_{H_\Phi(\Omega_2)}\leq {\mathcal O}_r(h^{N})\left(\|
\tilde{u}\|_{H_\Phi(\Omega_1)}+|u_-|_{\cz^a}\right).
\end{equation}
 This gives
\begin{equation}\label{1intt}
\| \tilde{u}\|_{H_\Phi(\Omega_1)} +|u_-|_{\cz^a}\leq  {\mathcal O}
\left(h^{-\delta}\right)\left( \|Uv\|_{H_\Phi(\Omega_2)}
+|v_+|_{\cz^a}\right)+{\mathcal O}_r \left( h^{N}\right)\left( \|
\tilde{u}\|_{H_\Phi (\Omega_1)} + | u_-|_{\cz^a}\right).
\end{equation}
{\bf Estimate  on $\Omega_{\rm int}\setminus V_{0}:$}\\
In  $(I-zM_0)\tilde{u}=-R_-u_-+\tilde{v}$ we use (a variant of)
Theorem \ref{Th.2.6}, estimate (\ref{est}), which implies that there is $N_1>0$ such that
$$\|z\chi M_0\tilde{u}\|_{L^2_\Phi(\Omega_2')}\leq{\mathcal O}(h^{N_1})\|\tilde{u}\|_{H_\Phi(\Omega_1')},\,\,\supp\chi\subset \Omega_{\rm
int}\setminus V_{0}.$$ We have then
$$\|\chi
(1-zM_0)\tilde{u}\|_{L^2_\Phi(\Omega_2')}\geq\|\chi\tilde{u}\|_{L^2_\Phi(\Omega_1')}
-{\mathcal O}(h^{N_1})\|\tilde{u}\|_{H_\Phi(\Omega_1')}$$ and
$$ \|\chi\tilde{u}\|_{L^2_\Phi(\Omega_1')}\leq \| \chi R_-u_-\|_{L^2_\Phi(\Omega_2')}
+\|\chi \tilde{v}\|_{L^2_\Phi(\Omega_2')}+{\mathcal
O}(h^{N_1})\|\tilde{u}\|_{H_\Phi(\Omega_1')}.$$
 This implies for any $N\in\nz$
\begin{equation}\label{3int}
\|\chi\tilde{u}\|_{L^2_\Phi(\Omega_1')}\leq {\mathcal
O}(h^{N})|u_-|_{\cz^a} +\|\chi \tilde{v}\|_{L^2_\Phi(\Omega_2')}
+{\mathcal O}(h^{N_1})\| \tilde{u}\|_{H_\Phi(\Omega_1')}.
\end{equation}

 Similar, applying (\ref{est}) to  each term in $\tilde{w}$ separately, we get
\begin{equation}\label{20int}
\|\chi\tilde{w}\|_{L^2_\Phi(\Omega_2')}\leq {\mathcal
O}(h^{N_1})\left(\|
\chi\tilde{u}\|_{L^2_\Phi(\Omega_1')}+|u_-|_{\cz^a}\right)
\end{equation} and we have
\begin{equation}\label{3intt}
\|\chi\tilde{u}\|_{L^2_\Phi(\Omega_1')}\leq {\mathcal
O}(h^{N_1})|u_-|_{\cz^a} +\|\chi Uv\|_{L^2_\Phi(\Omega_2')}
+{\mathcal O}(h^{N_1})\| \tilde{u}\|_{H_\Phi(\Omega_1')}.
\end{equation}

Applying  (\ref{1intt}) and (\ref{3intt}) with $\tilde{u}$ replaced
by $Uu$  we get
\begin{equation}\label{12int}
| u|_{\Omega_1'} +|u_-|_{\cz^a}\leq {\mathcal O} (h^{-\delta})\left(
|v|_{\Omega_2'} +|v_+|_{\cz^a}\right)+{\mathcal O} \left(
h^{N_1}\right)\left( | u|_{\Omega_1'} + | u_-|_{\cz^a}\right),
\end{equation}
where
$\Omega_2'\subset\subset\Omega_1'\subset\subset\Omega_{\rm int}.$\\
{\bf Estimate on $\Omega_{\rm ext}:$}\\
In order to estimate $u=Mu-VR_-u_-+v$ we use  bound (\ref{boundCW}) with some $N_2 >0:$
$$|Mu|_{\Omega_{\rm ext}}\leq {\mathcal
O}(h^{N_2})|u|_{\Omega_{\rm ext}}$$ and get for any $N\in\nz$
$$|u|_{\Omega_{\rm ext}}\leq | Mu|_{\Omega_{\rm ext}}
+|v|_{\Omega_{\rm ext}} +|VR_-u_-|_{\Omega_{\rm ext}}\leq {\mathcal
O}(h^{N_2})| u|_{\Omega_{\rm ext}} +| v|_{\Omega_{\rm ext}}
+{\mathcal O}(h^N)| u_-|.
$$
Summing up we get
\begin{align*} &|u|_{\Omega_{\rm ext}} +|u|_{\Omega_{\rm int}}
+|u_-|_{\cz^a}\leq\\
&\leq {\mathcal O}(h^{-\delta})\left(|v|_{\Omega_{\rm ext}} +|v|_{\Omega_{\rm
int}} +| v_+|_{\cz^a}\right) +{\mathcal O}(h^{N_0})\left(|
u|_{\Omega_{\rm ext}} +| u|_{\Omega_{\rm int}}
+|u_-|_{\cz^a}\right),
\end{align*}
with $N_0=\min (N,N_1,N_2).$
The last term can be absorbed in the left hand side and we get
  estimate (\ref{finalestimate}).

Then it is clear that $${\mathcal P}=\left(%
\begin{array}{cc}
  I-M(\lambda) & VR_- \\
  R_+U &  \\
\end{array}%
\right):\,\,{\mathcal H}\times\cz^a\mapsto{\mathcal H}\times\cz^a$$
is Fredholm of index $0.$ Estimate (\ref{finalestimate}) implies
that ${\mathcal P}$ is injective and thus bijective for
$\lambda\in\Lambda_{A,B}(\epsilon_0).$

 \hfill\qed

Denote $${\mathcal E}=\left(%
\begin{array}{cc}
  E & E_+ \\
  E_- & E_{+-} \\
\end{array}%
\right)$$ the inverse of ${\mathcal P}.$ Let $h=1/\Re\lambda,$
$\lambda\in\Lambda_{A,B}.$ Then it is known
$$0\in\sigma (I-M(\lambda))\,\,\Leftrightarrow\,\,
0\in\sigma (E_{+-}(\lambda)).$$ Let $v_+\in\cz^a.$ If
$(u,u_-)={\mathcal E} (0,v_+)$ then $u=E_+(v_+),$ $u_-=E_{+-}(v_+).$

Let $$ \left\{\begin{array}{rl}
  (I-M)u + V R_-u_-&=0,  \\
  R_+ U u&=v_+.
\end{array}\right.
$$
Then we have $$(I-e^{-2id\lambda}M_0)\tilde{u} +
R_-u_-=\tilde{w}=(UM-e^{-2id\lambda} M_0U) u+(I-U V)R_-u_-$$ and
$\tilde{w}$ satisfies (\ref{2int}):
$$\|\tilde{w}\|_{H_\Phi(\Omega)}\leq {\mathcal
O}_r(h^{N})\left( | u|_{\mathcal H}+|u_-|_{\cz^a}\right)={\mathcal
O}_r(h^{N}) |v_+|_{\cz^a}.$$ We have then for $r$ large enough
$$ \left\{\begin{array}{rl}
  \tau_J\left( I-e^{-2id\lambda}e^{-i\lambda F^r(I( 1/\lambda);1/\lambda)}\right)\tau_J\tilde{u} +
R_-u_- &={\mathcal O}_r(h^{N}) |v_+|_{\cz^a},  \\
   \tau_J\tilde{u}&=R_-v_+
\end{array}\right.
$$ and \begin{align*} &\left( I-e^{-2id\lambda}e^{-i\lambda F^r(I(1/\lambda);1/\lambda)}\right) R_-v_+ +
R_-u_- ={\mathcal O}_r(h^{N})
|v_+|_{\cz^a}\,\,\Leftrightarrow\\
&u_-=-R_+\left( I-e^{-2id\lambda}e^{-i\lambda
F^r(I(1/\lambda);1/\lambda)}\right) R_-v_++{\mathcal O}_r(h^{N})
|v_+|_{\cz^a}.\end{align*} We have then for
$\lambda\in\Lambda_{A,B}(\epsilon_0):$
$$E_{+-}(\lambda)=E_{+-}^0(\lambda) +{\mathcal O}_r\left(
(\Re\lambda_{\alpha,k}^r)^{-N}\right),\,\,\mbox{where}\,\,E_{+-}^0=-\left\{\left(
1-e^{i2d\lambda}K_{\alpha}(\lambda)\right)\right\}_{|\alpha\in
J}.$$

 Let $\gamma=\gamma(\lambda_{\alpha,k}^r)$
  be simple loop such that
$E_{-+}^0(\lambda)$ is invertible for $\lambda\in\gamma.$ Suppose
$$\lambda\in\gamma\,\,\leftrightarrow\,\,\lambda=\lambda_\alpha^r+{\mathcal
O}_r\left( (\Re\lambda_{\alpha ,k}^r)^{-N}\right),$$ where
$\lambda_{\alpha,k}^r$ is solution of (\ref{lambda_0}).
 Then
 the number of
poles of $(I-M(\lambda))^{-1}$ inside $\gamma$ counted with
multiplicity is equal to the number of roots
$m(\gamma)=m(\lambda_{\alpha,k}^r)$ inside $\gamma$ of $\det
E_{-+}^0(\lambda)=0,$
\begin{equation*}
m(\gamma)=\frac{1}{2\pi i}\int_\gamma\frac{\left( \det
E_{-+}^0(\lambda)\right)'}{\det E_{-+}^0}d\lambda .\end{equation*}
We have $0\in\sigma (I-M(\lambda))$ if
$\lambda=\lambda_\alpha^r+{\mathcal O}_r\left( (\Re\lambda_{\alpha
,k}^r)^{-N}\right).$

 This accomplishes the proof of Theorem \ref{th-main-bis}.\hfill\qed
\appendix
\section{Some facts about the quantum billiard operator $M.$}\label{s-M}

\subsection{Boundness}

 Let
$\Omega_1$ be an obstacle with analytic boundary,
non-trapping.
 Let
$\Omega_1\subset B(0,R).$ Denote \begin{equation}\label{3.1}
U_{A,B}:=\{ \lambda;\,\,\Im\lambda\leq A |\lambda|^{1/3} -B\}.
\end{equation}
Let $R(\lambda)$ be the outgoing Dirichlet resolvent in
$\complement{\Omega_1}.$

\begin{theorem}
Let $\chi\in C_0^\infty (\rz^n)$ be such that $\chi\equiv 1$ in
$B(0,R).$ There exist $A,B>0$ such that $\chi R\chi (\lambda ),$
defined for $\Im\lambda <0$ has an analytic extension to the domain
(\ref{3.1}) as a bounded operator $
L^2(\complement{\Omega_1})\mapsto H_0^1(\complement{\Omega_1}),$
satisfying the estimate
$$\exists C, D >0;\,\,\|\chi R\chi\|_{\mathcal{L} (L^2(\complement{\Omega_1}),
H_0^1(\complement{\Omega_1}))}\leq Ce^{D\Im\lambda^+},$$  where
$\Im\lambda^+=\max (\Im\lambda ,0).$
\end{theorem}

For $i=1,2$ let $$H_{i,+}(\lambda):\,\, C^\infty(\partial\Omega_i)\mapsto
C^\infty(\overline{\rz^{n+1}\setminus\Omega_i})$$ be the outgoing
resolvent of the problem
$$\left\{\begin{array}{l}
  (\Delta +\lambda^2) H_{i,+}(\lambda) u=0\,\,\mbox{in}\,\,\rz^{n+1}\setminus\Omega_i \\
  H_{i,+}(\lambda) u_{|\partial\Omega_i}=u
\end{array}\right.$$ extended
as an  operator ${\displaystyle H^{1/2}(\partial\Omega_i)\mapsto
H_{\rm loc}^1 (\overline{\complement\Omega_i})}.$ In
\cite{BardosLebeauRauch1987} it is proven that this resolvent,
analytical for $\Im\lambda <0$ has an analytical extension  to the
domain of the form (\ref{3.1}) as a bounded operator ${\displaystyle
H^{1/2}(\partial\Omega_i)\mapsto  H_{0,{\rm loc}}^1
(\overline{\complement\Omega_i})}$
 and satisfy the following estimate:
$$\forall R >0,\,\,\exists C>0,\,\,\exists D
>0,\,\,\forall\lambda\in U_{A,B},\,\,\|H_{i,+}\|_{{\mathcal L}(H^{1/2}(\partial\Omega_1), H^1(\complement{\Omega_1}\cap B(0,R)))}\leq C
e^{D\Im\lambda^+}.$$

We denote $H_{21}(\lambda)=H_{2,+}\circ\gamma_2\circ H_{1,+},$ where
$\gamma_i$ is the operator of restriction to $\partial\Omega_i$, the
outgoing resolvent of the problem
$$\left\{\begin{array}{rl}
  (\Delta +\lambda^2) H_{21} u &  =0,\,\,\mbox{in}\,\,\complement{\Omega}_2 \\
  H_{21}u_{|_{\partial\Omega_2}} & =H_{1,+} u_{|_{\partial\Omega_2}}.  \\
\end{array}
\right.$$

 We
define $H_i(\lambda)u=H_{i,+}(\lambda) u_{|\partial\Omega_{i+1}},$
where $\partial\Omega_3=\partial\Omega_1$ and
$$M(\lambda)=H_2(\lambda) H_1(\lambda)=\gamma_1H_{2,+}\gamma_2 H_{1,+}.$$

\begin{lemma}[Burq] Operator $M(\lambda)$ defined on ${\displaystyle H^1(\partial{\Omega}_1)\mapsto H^1(\partial{\Omega}_1)}$ for $\Im\lambda <0$ has
an analytic extension in the domain of the form (\ref{3.1}) and
there satisfies the following estimate
$$\exists D>0,\,\,\exists C>0,\,\,\forall\lambda\in
U_{A,B},\,\,\|M(\lambda)\|_{{\mathcal
L}(H^{1/2}(\partial{\Omega}_1))}\leq C|\lambda|^2
e^{D\Im\lambda^+}$$
\end{lemma}

Moreover, as in a \neigh{} of $\Omega_1,$ $Mu$ satisfies $(\Delta
+\lambda^2)Mu=0,$ we get
$$\exists D>0,\,\,\forall s\in\rz,\,\,\exists
C_s>0,\,\,\forall\lambda\in U_{A,B},\,\,\| M(\lambda)\|_{{\mathcal
L}(H^s(\partial\Omega_1))} \leq C_s e^{D\Im\lambda^+}
|\lambda|^{s+2}.$$

\subsection{Microlocal unitarity of $M$ with respect to the flux
norm.}\label{ss-unitary} We suppose that $\lambda$ is real. Denote
$h=1/\lambda.$  Let $v$ satisfy
\begin{equation}\label{eqv}
\left\{\begin{array}{rl}
  Pu:=(-h^2\Delta -1)u & =0,\,\, u\,\,\mbox{outgoing} \\
  \gamma_1 u & =v. \\
\end{array}\right. \end{equation}

We define $H_{1,+}:\,\,{\mathcal
D}'(\partial\Omega_1)\mapsto{\mathcal
D}'(\overline{\complement\Omega_1})$ the solution operator. In the
similar way we define $H_{2,+}:\,\,{\mathcal
D}'(\partial\Omega_2)\mapsto{\mathcal
D}'(\overline{\complement\Omega_2}),$ which satisfies $PH_{2,+} v=0$
and $\gamma_2H_{2,+} v=v.$ Let $H_1=\gamma_2H_{1,+} $ and
$H_2=\gamma_1H_{2,+}. $  Then
$M=H_2H_1:\,\,H^1(\partial\Omega_1)\mapsto H^1(\partial\Omega_1).$

The billiard operator $M$ can be identified with the monodromy
operator $M$ as in \cite{SjostrandZworski2001} in the form presented
in (\cite{IantchenkoSjostrandZworski2002}, p. 360). Let ${\rm
ker}_{\rho_1}(P)$ be local kernel $P$ near $\rho_1=(a_1,0).$ The
solution operator $K$ in (\cite{IantchenkoSjostrandZworski2002}) is
the operator $H_{1,+},$
$$H_{1,+}v\in{\rm ker}_{\rho_1} (P),\,\,
H_{1,+}v_{|\partial\Omega_1}=v.$$ We have
$${\mathcal M}=H_{1,+}\gamma_1 H_{2,+}\gamma_2:\,\,{\rm
ker}_{\rho_1}(P)\mapsto {\rm ker}_{\rho_1}(P).$$ As in
(\cite{IantchenkoSjostrandZworski2002}) we have identification ${\rm
ker}_{\rho_1}(P)\simeq {\mathcal D}'(\partial\Omega_1)$ via
$K=H_{1,+}.$ Then the monodromy operator $M$ on ${\mathcal
D}'(\partial\Omega_1)$ satisfies $$KMv={\mathcal M}Kv,\,\,v\in
{\mathcal D}'(\partial\Omega_1).$$


Let $\chi\in C^\infty(\overline{\complement\Omega_1})$ be a
microlocal cut-off function supported in a \neigh{} $W_2$ of
$\Omega_1$ such that $\chi =1$ in $\complement\Omega_1\cap W_1,$
where $W_1\subset\subset W_2.$

We define the {\em quantum flux} norm (see
\cite{IantchenkoSjostrandZworski2002}, p. 360) on the outgoing
solutions $u$ of (\ref{eqv}) as follows
$$\|u\|_{\rm QF}^2:=\langle\frac{i}{h}[P,1-\chi ]
u|u\rangle.$$

It is easy to see that $\|u\|_{\rm QF}$ is independent of $\chi$
which implies that $M$ is microlocally unitary with respect to
$\|.\|_{\rm QF}$ for real $h:$
\begin{lemma} The billiard operator $M$ is microlocally unitary for all real $\lambda$ with
respect to  $\|.\|_{\rm QF}:$ $\|{\mathcal M} Kv\|_{\rm QF}=\|
H_{1,+}M v\|_{\rm QF}=\| H_{1,+}v\|_{\rm QF}+{\mathcal
O}(h^\infty),$ where $u=H_{1,+}v$ is the outgoing solution of
(\ref{eqv}) and $WF_h(v)\subset{\rm neigh}(\rho_1).$ Here the wave
front set $WF_h(v)$ is defined as in \cite{Gerard1988}.
\end{lemma}

Using the Green's formula in a \neigh{} of $\Omega_1$ bounded by
$\partial\Omega_1$ on one side, we get
\begin{align*}&\|u\|_{\rm QF}^2=\langle\frac{i}{h}[P,1-\chi ]
u|u\rangle=\frac{i}{h}\langle (h^2\Delta +1)\chi u -\chi(h^2\Delta
+1)u|u\rangle =\\
&=\frac{h}{i}\int_{\partial\Omega_1} \{(\partial_\nu
H_{1,+}v)\overline{H_{1,+}v}
-H_{1,+}v(\partial_\nu\overline{H_{1,+}v})\}
S(dx)=\\&=2\Re\langle\frac{h}{i}\partial_\nu
H_{1,+}v|v\rangle_{L^2(\partial\Omega_1)}=2\Re\langle
Av|v\rangle_{L^2(\partial\Omega_1)},\,\,u=H_{1,+} v,
\end{align*}
where $\partial_\nu$ is the normal derivative,
$Av:=\gamma_1\frac{h}{i}\partial_\nu H_{1,+}v.$  As $A$ is
\pseudor{} of order $0,$ elliptic near $(a_1,0)$ (we use the
explicit WKB construction of $H_{1,+}$ in the hyperbolic zone) then
we have for $v$ with $WF_h(v)\subset{\rm neigh}(\rho_1)$
$$\|H_{1,+}v\|_{\rm QF}\sim \| v\|_{L^2(\partial\Omega_1\cap\, {\rm neigh}(\rho_1))}.$$

\begin{lemma} There exists \pseudor{} $B$ of order $0$ elliptic near $a_1$ such that
$BMB^{-1}$ is  microlocally unitary for all real $\lambda$ with
respect to $\|.\|_{L^2(\partial\Omega_1)}:$ $$\| B M B^{-1}
v\|_{L^2(\partial\Omega_1)}=\| v\|_{L^2(\partial\Omega_1)}+{\mathcal
O}(h^\infty),\,\,WF(v)\subset{\rm neigh}(a_1,0).$$
\end{lemma}


\end{document}